\documentclass[12pt]{amsart}
\usepackage{amsbsy}
\usepackage{graphicx,epsfig,subfigure,psfrag}

\newtheorem{theorem}{Theorem}%[section]
\newtheorem{proposition}{Proposition}%[section]
\newtheorem{lemma}{Lemma}%[section]
\newtheorem{definition}{Definition}%[section]
\newcommand{\beq}{\begin{equation}}
\newcommand{\eeq}{\end{equation}}
\numberwithin{equation}{section}
\numberwithin{theorem}{section}
\numberwithin{proposition}{section}
\numberwithin{lemma}{section}
\numberwithin{corollary}{section}
\numberwithin{definition}{section}
\numberwithin{remark}{section}
\newcommand{\re}{{\mathbb R}}
\newcommand{\n}{\nabla}
\newcommand{\ren}{{\mathbb R}^N}
\newcommand{\iy}{\infty}

\renewcommand{\a}{\alpha}
\renewcommand{\b}{\beta}
\newcommand{\g}{\gamma}

\newcommand{\D}{\Delta}
\newcommand{\e}{\varepsilon}
\newcommand{\var}{\varphi}
\renewcommand{\l}{\lambda}

\renewcommand{\O}{\Omega}
\newcommand{\s}{\sigma}

\newcommand{\noi}{\noindent}
\newcommand{\inA}{\quad \mbox{in} \quad \ren \times \re_+}
\newcommand{\inB}{\quad \mbox{in} \quad}

\newcommand{\forA}{\quad \mbox{for} \quad}
\newcommand{\whereA}{,\quad \mbox{where} \quad}

\newcommand{\andA}{\quad \mbox{and} \quad}
\newcommand{\ssk}{\smallskip}

\title %%%[Blow-up localization]
{\bf Variational approach to complicated  similarity solutions of
higher-order
 nonlinear  %%%%%evolution
 PDEs. II}
\begin{document}

%%{\bf Refined multiple structure of blow-up\\
%% and compacton
%%similarity solutions\\
%%% of some higher-order nonlinear PDEs. I}

%%{\bf On refined multiple structure of blow-up and compacton
%%similarity solutions of some higher-order nonlinear PDEs. II}

%% {\bf Variational approach to complicated similarity solutions
%%of higher-order
%% nonlinear evolution equations of parabolic, hyperbolic, and nonlinear dispersion
%% types. II}
 %%%%. Regular operators and homotopy classification}

\author{
V.A.~Galaktionov, E. Mitidieri,  and S.I.~Pohozaev}

\address{Department of Mathematical Sciences, University of Bath,
 Bath BA2 7AY, UK}
 %and Keldysh Institute of Applied Mathematics,
 %Miusskaya Sq. 4, 125047 Moscow, RUSSIA}
\email{vag@maths.bath.ac.uk}

\address{Dipartimento di Matematica e Informatica, Universit\`a  di Trieste, Via A. Valerio 12/1 34127  Trieste, ITALY}
 \email{mitidier@units.it}

\address{Steklov Mathematical Institute,
 Gubkina St. 8, 119991 Moscow, RUSSIA}
\email{pokhozhaev@mi.ras.ru}

 \keywords{Quasilinear ODEs, non-Lipschitz terms,
similarity solutions, blow-up, extinction, compactons,
Lusternik--Schnirel'man category, fibering method.}
%% multiplicity.
 %% \\
%%{\bf Submitted to}: ???.}
 \subjclass{35K55, 35K40, 35K65}
\date{\today}
%%%% \quad {\bf PohF2m/F2mII.tex}}

%%%%%%%%%%%%%%%%%%%%%%%%%%%%%%%%%
% \newcommand{\begin{equation}}{\begin{equation}}
%%%%%%%%%%%%%%%%%%%%%%%%%%%%%%%%%%%%%%%%%%%%%%%
% \newcommand{\end{equation}}{\end{equation}}
%%%%%%%%%%%%%%%%%%%%%%%%%%%%%%%%%%%%%%%%%%%%%%%%%%%%%
%%%%%%

%%%%%%%%%%%%%%%%%%%%%%%%%%%%%%%%%%%%%%%%%%%%%%%%%%%%%%%%

% \vskip 1.5cm

\begin{abstract}

%%%We study
 %%Some common features of
  %%blow-up and travelling wave
 %%compacton-type
 %%self-similar
 This paper continues the study began in \cite{GMPSob, GMPSobIarX}
 of
   the Cauchy problem for  $(x,t) \in \ren
  \times \re_+$ for three higher-order degenerate  quasilinear partial differential equations (PDEs), as
  basic models,
 $$
  \begin{matrix}
 u_t = (-1)^{m+1}\D^m(|u|^n u)+|u|^n u,\qquad\,\quad\ssk\\
 u_{tt} =(-1)^{m+1} \D^m(|u|^n u)+|u|^n u,\qquad\quad\ssk\\
 u_{t}= (-1)^{m+1}[\D^m(|u|^n u)]_{x_1}+(|u|^n
 u)_{x_1},
 \end{matrix}
 $$
  where $n>0$ is a fixed exponent and  $\D^m$ is the $(m \ge 2)$th iteration of the
  Laplacian.
  %%is considered.
   A
  diverse  class of degenerate PDEs from various areas of applications
 of three types:
  %%$D_x=\frac{ \partial}{\partial x}$, $n>0$, and  $m \ge
  %%2$,
  parabolic, hyperbolic, and nonlinear dispersion, is dealt with.
General local, global, and blow-up features of such PDEs are
studied on the basis of
  their {\em blow-up} similarity or
{\em travelling wave} (for the last one) solutions.

In \cite{GMPSob,GMPSobIarX},  Lusternik--Schnirel'man category
theory of variational calculus  and fibering methods were applied.
The case $m=2$ and $n>0$ was studied in greater detail
analytically and numerically. Here,  more attention is paid to a
combination of a Cartesian approximation
 and fibering to get new compactly supported similarity patterns.
 Using numerics,
 such
 compactly supported
 solutions constructed   for $m = 3$ and for higher orders. The ``smother" case of negative
%% and also include the smooth case
$n <0$ is included, with a typical ``fast diffusion-absorption"
parabolic PDE:
 $$
u_t = (-1)^{m+1}\D^m(|u|^n u)-|u|^n u \whereA n \in (-1,0),
 $$
 which admits {\em finite-time extinction} rather than blow-up.
%%%%for the above PDEs reflects the case extinction patterns.
 Finally, a homotopy approach is developed for some kind of
 classification of various patterns obtained by variational and
 other methods.
Using a variety of analytic, variational,
 qualitative, and numerical methods
allows
 to justify that  the above PDEs admit an   infinite
 countable set  of
countable families of compactly supported blow-up (extinction) or
travelling wave solutions.
%% that are
%%oscillatory near finite interfaces. In a whole, this solution set
%%%%m exhibit typical features of being of a chaotic structure.

\end{abstract}

%%%%%%%%%%%%%%%%%%%%%%%%%%%
\maketitle

%%%\ms

%%\begin{center}
%%{\em Dedicated to the memory of Professor S.L.~Sobolev on his
%%Centenary}
%% on the
%%occasion of his 100th birthday anniversary}
%%\end{center}

%%%%%%%%%%%%%%%%%%%%%%%%%%%

%%%%%%%%%%%%%%%%%%%%%%%%%%%%%%%%%%%%%%%%%%%%%%%%%%%%%%%%%%%%%%%%%%%%%%%%
% INTRODUCTION
%%%%%%%%%%%%%%%%%%%%%%%%%%%%%%%%%%%%%%%%%%%%%%%%%%%%%%%%%%%%%%%%%%%%%%%%
\setcounter{equation}{0}
\section{Introduction:
 higher-order blow-up and compacton models}
%% and  similarity
%%%solutions}
%%via an
%% analytic-numerical approach}
%%%Brezis--Friedman results  and extensions to
%%%%%higher-order PDEs}
 \label{Sect1}
  \setcounter{equation}{0}

A general physical and PDE motivation of the present research can
be found in \cite{GMPSob, GMPSobIarX}, together with basic history
and related key references, so we just briefly comment on where
%%%% this kind of
quasilinear elliptic problems under consideration
%%%with various nonlinearities
 are coming from.

%%%%%%%%%%%%%%%%%%%%%%%%%%%%%%%%%%%%%%%%%%
\subsection{(I) Combustion-type models with  blow-up}
%% and main
%%goals}
%%and first discussion}

Our first model is a
 quasilinear degenerate $2m$th-order parabolic
equation of the reaction-diffusion (combustion) type:
 \beq
    \label{S1}
 u_t = (-1)^{m+1}\D^{m}(|u|^n u)+|u|^n u \inA,
 \eeq
 where $n>0$ is a fixed exponent, $m \ge 2$ is integer, and $\D$ denotes the Laplace operator in
 $\ren$.
 %%For $m=1$, this is the classical semilinear heat equation from combustion theory
 %% $$
 %% u_t=\D u + u^{p} \quad (u \ge 0)
 %%$$
  Physical, mathematical, and blow-up history of
 (\ref{S1}) for the standard classic
 case $m=1$ and $m \ge 2$ is explained in \cite[\S~1]{GMPSob,GMPSobIarX}.
 Consider
 %% indicate the
 {\em regional blow-up}  solutions of (\ref{S1})
\beq
 \label{RD.31}
 u_{\rm S}(x,t)=(T-t)^{-\frac 1n} f(x) \quad \mbox{in}
 \quad \ren \times (0,T)
  \eeq
  in separable variables,
  where $T>0$ is the blow-up time. Then
the similarity blow-up profile $f=f(x)$  solves a quasilinear
elliptic equation of the form
 %%%%a degenerate equation %%%an ODE with non-Lipschitz nonlinearity
 \beq
 \label{mm.561}
  \mbox{$
 (-1)^{m+1} \D^m (|f|^n f) + |f|^n f= \frac 1n \, f \inB  \ren.
 $}
  \eeq
This reduces to the following semilinear equation with a
non-Lipschitz nonlinearity:
 $$ %%\beq
 %%\label{S200}
 \mbox{$
 F=|f|^n f \quad \Longrightarrow \quad
 \mbox{$(-1)^{m+1} \D^m F+F- \frac 1n \,
 \big| F \big|^{-\frac n{n+1}} F=0 \quad \mbox{in} \,\,\, \ren.
  $}
 $}
  $$ %%%\eeq
  Scaling out the multiplier $\frac 1n$ in the nonlinear term
  yields
\beq
 \label{S2NN}
 \mbox{$
 F \mapsto n^{-\frac{n+1}n} F \quad \Longrightarrow \quad
 \fbox{$(-1)^{m+1} \D^m F+F-
 \big| F \big|^{-\frac n{n+1}} F=0 \quad \mbox{in} \,\,\, \ren.
  $}
 $}
  \eeq
 For $N=1$, this is a simpler ordinary differential equation (an ODE):
\beq
 \label{S2}
 \mbox{$
 F \mapsto n^{-\frac{n+1}n} F \quad \Longrightarrow \quad
 \fbox{$(-1)^{m+1} F^{(2m)}+F-
 \big| F \big|^{-\frac n{n+1}} F=0 \quad \mbox{in} \,\,\, \re.
  $}
 $}
  \eeq
According to (\ref{RD.31}), the elliptic problems (\ref{S2NN}) and
the ODE (\ref{S2}) for $N=1$ are responsible for the possible
``geometrical shapes" of regional blow-up described by the
higher-order combustion model (\ref{S1}).

\ssk

%%%%%%%%%%%%%%%%%%%%%%%%%%%%%%%%%%%%%%%%%%%%%%%%%%%%%
\noi{\bf Remark: relation to ODEs from extended KPP theory}. There
exists vast mathematical literature, starting essentially from the
1980s, devoted to  the fourth-order ODEs (looking rather
analogously to that in (\ref{S2}) for $m=2$)
 \beq
 \label{KPP1}
  F^{(4)}= \b F'' + F-F^3 \quad \mbox{in} \quad \re,
   \eeq
   where $\b>0$ is a parameter.
   This ODE also admits a complicated set of solutions
 with various classes of patterns and even with
  chaotic features.
 We refer to  Peletier--Troy's book \cite{PelTroy} for the most diverse  account, as well as
 %% and refer to
 %% results and lists of references in
 to papers \cite{KKVV00, VV02}, where a detailed and advanced
 solution description for (\ref{KPP1})  is obtained  by combination of
 variational and  homotopy theory.
 Regardless their rather similar forms,
 %%it turns out that our
 the ODEs (\ref{S2}) belong to a completely
 different class of equations with {\em non-coercive} operators,
 unlike in (\ref{KPP1}). Therefore, direct homotopy approaches and
 several others, that used to be rather effective for (\ref{KPP1}), fail
 in
 principle for (\ref{S2}). In this sense, (\ref{S2}) is  similar to the cubic ODE
 to be studied in Section \ref{SectAn}:
  %%%%with cubic nonlinearity
  \beq
  \label{nn1}
   \fbox{$
  F^{(4)}= - F +F^3 \quad \mbox{in} \quad \re,
  $}
   \eeq
  %%% to be studied in  Section \ref{SectAn}.
of course, excluding complicated oscillatory behaviour at finite
interfaces, which are obviously nonexistent for analytic
nonlinearities. However, we claim that the sets of solutions of
(\ref{S2}), with $m=2$, and  of (\ref{nn1}) are equivalent,
though, not having a rigorous proof, we will devote some efforts
to a homotopy approach connecting solutions of such smooth
(analytic) and non-smooth ODEs. Thus, though going to develop
homotopy approaches for classifying solutions of (\ref{S2})
(Section \ref{SectHom}), our main tool to describe countable
families of solutions $\{F_l\}$ is a combination of
Lusternik--Schnirel'man category-genus theory \cite{KrasZ}
and the fibering method \cite{Poh0, PohFM}.
%% This is
%%done in Sections \ref{SectVar} and \ref{SFFh}.

%%%\noi{\bf (vi)} \underline{\sc Problem ``Pattern classification"}:
%%developing a ``homotopy" approach to classify countable sets of
%%%patterns (Section \ref{???}).

\ssk

 Thus, we show that ODEs (\ref{S2}), as
well as the PDE (\ref{S2NN}),  admit infinitely many countable
families of compactly supported solutions, and the whole solution
set exhibits certain {\em chaotic} properties.
%%% (and moreover an uncountable set) of {\em
%%%chaotic} orbits that are compactly supported in $\re$.
Our analysis will be based on a combination of analytic
(variational and others), numerical, and some more formal
techniques. Explaining existence,
 multiplicity, and asymptotics for the nonlinear problems
 involved, we state and
leave several open difficult mathematical problems.
 %% Some of these for higher-order equations are  extremely
 %%%difficult.
Meantime, let us characterize other models involved.

 %%\ssk

%%%%%%%%%%%%%%%%%%%%%%%%%%%%%%%%%%%%%%%
 \subsection{(II) Regional blow-up in quasilinear hyperbolic equations}

%%%Meantime, we continue to describe main  applications of the
 %%%%%presented mathematics.
  Consider next the  $2m$th-order hyperbolic  counterpart of
  (\ref{S1}),
 \beq
 \label{S3}
 u_{tt} =(-1)^{m+1} \D^{m}(|u|^n u)+|u|^n u \inA.
  \eeq
The blow-up solutions take a similar form with a different
exponent $-\frac 2n$ instead of $-\frac 1n$:
 \beq
 \label{RD.31H}
 \mbox{$
 u_{\rm S}(x,t)=(T-t)^{-\frac 2n}  f(x),
 %%% \quad \Longrightarrow \quad
 %%%\frac 2n \bigl(\frac 2n+1 \bigr) \, \tilde f= (\tilde f^{n+1})'' +  \tilde f^{n+1}.
 $}
  \eeq
%%For the $N$-dimensional PDE (\ref{S3}),
%% \beq
%% \label{S3}
%% u_{tt} =(-1)^{m+1} \D^{m}(|u|^n u)+|u|^n u \inA,
%%  \eeq
where $F=|f|^n f$, after  scaling,
%%% looking for  the same solution (\ref{RD.31H}),
%%%%, (\ref{S5})
%%% after scaling, leads to
solves  the same elliptic equation (\ref{S2NN}).

%%%%%%%%%%%%%%%%%%%%%%%%%%%%%%%%%%
 \subsection{(III) Nonlinear dispersion equations and compactons}

Such rather unusual PDEs in $N$-dimensions (the origin is
integrable PDEs and other
 areas) take the form
 \beq
 \label{NDEN}
 u_t= (-1)^{m+1}[\D^m(|u|^n u)]_{x_1} + (|u|^n u)_{x_1}
 \inA,
  \eeq
  where the right-hand side is the derivative  $D_{x_1}$ of
  that in the parabolic counterpart (\ref{S1}).
Then the elliptic problem (\ref{S2NN}) occurs when studying {\em
travelling wave} (TW) solutions of (\ref{NDEN}). Note that, as
being  PDEs with nonlinear dispersion mechanism, (\ref{NDEN}) and
other NDEs listed below
 admit
shock waves and other discontinuous solutions. Here, we study
smooth solutions and
 do not touch difficult entropy-like approaches
 for such shock and rarefaction waves and refer to
\cite[\S~4.2]{GSVR} and \cite{GPnde} for an account to such
phenomena.

Thus, for the PDE (\ref{NDEN}), looking for a TW {\em compacton}
(i.e., a solution having all the time compact support;
 see key
references in \cite[\S~1]{GMPSob}) moving in the $x_1$-direction
only,
 \beq
 \label{nn1NN}
  \mbox{$
 u_{\rm c}(x,t)=f(y_1,x_2,...,x_N), \,\,\, y_1=x_1-\l t \whereA \l= - \frac
 1n,
  $}
 \eeq
 we obtain on integration in $y_1$ the elliptic problem (\ref{S2NN}).
 Analogously, for the higher-order evolution extension of
nonlinear dispersion PDEs,
 $$
  \mbox{$
 D_t^k u = D^k_{x_1}\bigl[(-1)^{m+1} \D^m (|u|^n u)+|u|^n u \bigr] \inA \quad
 (k \ge 2),
  $}
 $$
 to get the same PDE (\ref{S2NN}),
the compacton (\ref{nn1NN})  demands the following wave speed:
 $$
  \mbox{$
 (-\l)^k= \frac 1n.
  $}
  $$

%%%%%%%%%%%%%%%%%%%%%%%%%%%%%%%%%%%%%%%%%%%%%%%
\subsection{(IV) PDEs with ``fast diffusion" operators:
 parabolic, Schr\"odinger, hyperbolic, and nonlinear dispersion models}

This is about negative exponents:
 \beq
 \label{n1}
 n \in (-1,0),
 \eeq
which generate other types of elliptic equations of interest. To
connect such problems with typical models of diffusion-absorption
type, consider the following parabolic PDE:
 \beq
    \label{S1F}
 u_t = (-1)^{m+1}\D^{m}(|u|^n u)-|u|^n u \inA \quad (-1<n<0),
 \eeq
with the strong non-Lipschitz at $u=0$ {\em absorption}  term. It
is well known that such PDEs describe finite-time extinction
phenomenon instead of blow-up. See \cite[Ch.~4,5]{AMGV} for $m=1$
and \cite{Galp1, Shi2} for $m \ge 2$ for necessary references and
history of strong absorption phenomena. Therefore,  the similarity
solution takes an analogous to  (\ref{RD.31}) form with the
positive exponent $- \frac 1n
>0$, so that $u(x,T^-) \equiv 0$, while $f$ solves a similar
elliptic equation  (cf. (\ref{mm.561}))
 \beq
 \label{mm.561F}
  \mbox{$
 (-1)^{m+1} \D^m (|f|^n f) - |f|^n f= \frac 1n \, f\quad \mbox{in} \,\,\, \ren.
 $}
  \eeq
 By the scaling as in (\ref{S2NN}) (recall that $n<0$),
  we eventually obtain the semilinear elliptic
 problem with a sufficiently smooth nonlinearity:
\beq
 \label{S2NNF}
 \mbox{$
 %%%F \mapsto n^{-\frac{n+1}n} F \quad \Longrightarrow \quad
 \fbox{$(-1)^{m+1} \D^m F - F+
 \big| F \big|^{\a} F=0 \quad \mbox{in} \,\,\, \ren \whereA
 \a=-\frac n{n+1}>0.
  $}
 $}
  \eeq
The nonlinearity is now  $C^1$ at $F=0$, so the solutions are
classic. For instance, for $m=2$ and $n=- \frac 23$, we obtain
equation with a cubic analytic nonlinearity:
 %%The ODE with non-coercive operators
 %%%such as (\ref{nn1})
 %%  and  similar elliptic problems,
 \beq
 \label{nn1E}
 %%% \quad \mbox{and} \quad
   \D^2 F = - F + F^3 \quad \mbox{in} \quad \ren.
 \eeq
   %% occur in the study of
For $N=1$, this is the ODE (\ref{nn1}).
 Indeed, these equations do not admit solutions with finite
 interfaces and exhibit {\em exponentially decaying oscillatory} behaviour at
 infinity. We show that the total set of such
 ``effectively" spatially localized patterns well matches those ones for
 $n>0$ always having finite interfaces.

%%%%%%%%%%%%%%%%%%%%%%%%%%%%%%%%%%%%%%%%
%%\subsection{

\ssk

 As a connection to another classic PDE area and applications,
let us note that, for $m=1$, we obtain  the classic second-order
case
 of the {\em ground state equation} \cite{Coff72}
  \beq
  \label{gr1}
  \D F - F + F^3 =0 \inB \ren.
   \eeq
   This elliptic problem is key in blow-up analysis of
 the critical {\em nonlinear Schr\"odinger equation} (NLSE)
  $$
  \mbox{$
   %%\begin{matrix}
 {\rm i}\,  u_t=- \D u -|u|^{p-1} u, \quad
%% \inB \ren \times \re_+, \,\,\, (N=2)
   p=p_0=1+ \frac 4N, \,\, N=2, \quad
%%%  \ssk\ssk\\
   \mbox{and} \quad u(x,t)={\mathrm e}^{{\rm i}\,
  t}F(x);
  %%\end{matrix}
  $}
  $$
  see Merle--Raphael \cite{Mer05}--\cite{MerR052} as a guide.
 The solution $F$ of (\ref{gr1})
    is {\em strictly positive} (with exponential decay at infinity) and
   is unique up to translations, while the ground state $F_0$ for
   (\ref{nn1E}) is oscillatory
   at infinity, to say nothing about a huge variety of other, Lusternik--Schnirel'man or not,
   solutions.
   Thus, (\ref{nn1E}) is the {\em ground state
   equation} for the fourth-order NLSE
 $$
  \mbox{$
 {\rm i}\,  u_t= \D^2 u -
|u|^{p-1} u
 \inB \ren \times \re_+
  \whereA p=p_0=1+ \frac 8N \,\,\, (N=4).
  $}
  $$

Analogous similarity analysis is performed for the corresponding
``fast diffusion" hyperbolic equation (the extinction patterns are
given by  (\ref{RD.31H}))
 \beq
 \label{S3F}
 u_{tt} =(-1)^{m+1} \D^{m}(|u|^n u)-|u|^n u \inA \quad (-1<n<0),
  \eeq
 and for the nonlinear ``fast dispersion" PDE
 \beq
 \label{NDENF}
 u_t= (-1)^{m+1}[\D^m(|u|^n u)]_{x_1} - (|u|^n u)_{x_1}
 \inA \quad (-1<n<0),
  \eeq
where the moving TW profiles are as in (\ref{nn1NN}).

%%%%%%%%%%%%%%%%%%%%%%%%%%%%%%%%%%%%%%%%%%%
\subsection{Main goals and connections with previous results}

 It turns out that
such profiles $F$ solving (\ref{S2NN}) have rather complicated
local and global structure. The study of equations (\ref{S2NN})
and (\ref{S2}) was began in \cite{GMPSob}, where the following
goals were posed:
%% We are not aware of any rigorous or even formal
%%qualitative results concerning existence, multiplicity, or global
%%structure of solutions of ODEs such as (\ref{S2}).
 %%%Therefore,

 \ssk

\noi{\bf (i)} \underline{\em Problem ``Blow-up"}: proving
finite-time blow-up in the parabolic (and hyperbolic) PDEs under
consideration \cite[\S~2]{GMPSob};

\ssk

 \noi{\bf (ii)} \underline{\em Problem ``Existence and Multiplicity"}:
 existence and
  multiplicity
   %%%, and Sturm-type %%%--Morse
 %% classification
  for elliptic PDEs (\ref{S2NN}) and the
ODEs (\ref{S2}) for $m \ge 2$
%%%% and similar ones with analytic nonlinearities
\cite[\S~3]{GMPSob};
%% \ref{SFFh});

\ssk

 \noi{\bf (iii)} \underline{\em Problem ``Oscillations"}: the
generic structure of  oscillatory
 solutions of (\ref{S2}) near interfaces for arbitrary $m \ge 2$ \cite[\S~4]{GMPSob};
 and

\ssk

\noi{\bf (iv)} \underline{\em Problem ``Numerics"}: numerical
study of all families of $F(x)$ for $m=2$ \cite[\S~5]{GMPSob}.
%%%% \ref{Sectm4}) including those that will not get a
 %%%%rigorous justification;

\ssk

The non-Lipschitz problem (\ref{S1}) possesses so complicated set
of admissible compactly supported solutions (note that Lusternik--Schnirel'man
category theory detects only a single countable subset), that
using effective {\tt MatLab} (or other similar or advanced)
numerical techniques for classifying the critical points becomes
an unavoidable tool of any analytic-numerical approach, which
cannot be dispensed with at all. We recall that in
\cite[\S~3]{GMPSob} the identification of Lusternik--Schnirel'man sequence of
critical values was confirmed numerically only (and for $m=2$
essentially).

\ssk

Therefore,  in the present paper, in Section \ref{Sectm4}, we
begin by continuing achieving the goal {\bf (iv)} for the
sixth-order case $m=3$:

\ssk

 \noi{\bf (iv$'$)} \underline{\em Problem ``Numerics"}: numerical study of all
families of $F(x)$ for $m \ge 3$ (Section \ref{Sectm4}).

\ssk

 In
addition, we also aim new targets:

\ssk

\noi{\bf (ii$'$)} \underline{\em Problem ``Existence and
Multiplicity"}: using variational Lusternik--Schnirel'man approach and fibering with
an auxiliary Cartesian approximation of  critical points (Section
\ref{SFFh});

\ssk

\noi{\bf (v)} \underline{\em Problem ``Fast diffusion"}: $n \in
(-1,0)$, where  smoother elliptic problems (\ref{S2NNF}) occur
(extinction in Sections  \ref{S.Ext} and existence-multiplicity in
Section \ref{SectAn}).

%%Let us explain how these new smoother (or even analytic) elliptic
%%equation occur.

\ssk

Finally, for both the  non-smooth (\ref{S2NN}) and smooth
(\ref{S2NNF}) problems, we pose:
%% the following:

\ssk

 \noi{\bf (vi)} \underline{\em Problem ``Sturm Index"}:  classification of
 various patterns according to their spatial shape for both (\ref{S2NN}) and (\ref{S2NNF})
(Sections \ref{SectHom}, \ref{SectAn}, and \ref{Sect8}). [For
$m=1$,
 this is governed by classic Sturm's First Theorem on zero sets.]

 \ssk

Thus, we are introducing three classes, {\bf (I), (II), (III)}, of
nonlinear higher-order PDEs in $\ren \times \re_+$ including
similar fast diffusion ones {\bf (IV)}. These are  representatives
of PDEs of three different types. However, it  will be shown that
these exhibit quite analogous evolution features (if necessary, up
to replacing blow-up by moving travelling waves  or extinction
behaviour), and  coinciding complicated countable sets of
evolution patterns.
 This
 %% common features
 reveals an exiting feature of
%%concept of a certain
 a certain {\em unified principle of singularity formation
phenomena}
 in general nonlinear PDE theory, which we would like to
believe in, but which is very difficult to justify more
rigorously.

%%% that we seem just begin
%%to touch and study in the twenty first century. Several classic
%%mathematical concepts and techniques successfully developed in the
%%twentieth century including, of course,  Sobolev
%% legacy continue to be  key, but also new ideas from different
%% ranges of various rigorous and qualitative natures are
%% desperately needed for tackling such fundamental difficulties and
%% open problem arising.

%%%%%%%%%%%%%%%%%%%%%%%%%%%%%%%%%%%%%%%%%%%%%%%%%%%%%%%%%%%%%%%
 \subsection{On
  extensions to essentially quasilinear equations}

The three-fold unity of PDE classes {\bf (I)--(III)} is available
for other types of nonlinearities.
 In Section \ref{SectExt1}, we briefly discuss the following
 classes of equations with fourth-order $p$-Laplacian
 operators (here $n=p-2>0$):
  \beq
  \label{Cl1}
  \begin{matrix}
  {\bf (I)} \quad u_t = -\D(|\D u|^n \D u) + |u|^n u \quad
  (\mbox{parabolic}), \,\,\ssk\ssk \\
  {\bf (II)} \quad u_{tt} = -\D(|\D u|^n \D u) + |u|^n u \quad
  (\mbox{hyperbolic}), \ssk\ssk\\
{\bf (III)} \quad u_t = -\big[\D(|\D u|^n \D u)\big]_{x_1} +
(|u|^n u)_{x_1} \quad
  (\mbox{NDE}).
   \end{matrix}
   \eeq
It turns out that these equations admit
 similar blow-up or
compacton (for the NDE {\bf (III)}) solutions that are governed by
variational elliptic problems with similar countable variety of
oscillatory compactly supported solutions. As a first step,  an
approach  to blow-up of solutions of the parabolic equation
(\ref{Cl1})
 for $N=1$ and
 some
 other related results on similarity solutions can be found in
 \cite{GpLap}.

 Thus, these cases are more difficult and are {\em essentially}
quasilinear, since the resulting elliptic problems cannot be
reduced to {\em semilinear} equations such as (\ref{S2}).
%% and
%%(\ref{S2N}).

%%%%%%%%%%%%%%%%%%%%%%%%%%%%%%%%%%%%%%%%%%%%%%%%
\subsection{Towards non-variational
 problems: branching}

Principally more difficult  problems occur under  a slight change
of the lower-order nonlinearity in
 (\ref{Cl1}):
  \beq
  \label{pn12}
|u|^n u \mapsto |u|^{p} u, \quad \mbox{with} \quad p>n.
 \eeq
This leads to {\em non-variational} elliptic problems, to be also
briefly discussed
 in Section \ref{SectExt1} using the idea of {\em branching} of proper
 solutions at $p=n$ from the similarity profiles studied above.
Proving existence of countable sets of solutions for such
non-potential operators reveals  a number of open problems of a
higher level of complexity.

%%\noi{\bf A general description of Lusternik--Schnirel'man critical points.} Thus, the
%%reason that the basic family of patterns $\{F_l\}$ shown in Figure
%%\ref{G4} with simple changing sign oscillatory almost periodic
%%structures, are generated by Lusternik--Schnirel'man critical point theory follows
%%directly from the minimax approach (\ref{ck1}). Namely, the outer
%%minimization $\inf_{{\mathcal F}\subset {\mathcal M}_k}$ demands
%%to use {\em closed subsets ${\mathcal F}$ of genus $\ge k$ of the
%%minimally possible structure and density}. Otherwise, the internal
%% $\sup_{v \in {\mathcal F}$ will be larger.
%%Therefore, the  minimax procedure (\ref{ck1})
%% imposes special restrictions on ${\mathcal F}$, and actually
%% {\em excludes from sets
%%${\mathcal F}$ {all} closed subsets and families that can be
%% embedded to subsets of lower genus}. This leaves in appropriate
%% sets ${\mathcal F}$ only simplest closed families of functions
%%of elementary almost periodic chains of $k$ humps.

%%%%%%%%%%%%%%%%%%%%%%%%%%%%%%%%%%%%%%%%%%%%%%%%%%%%%%%%%%%%%%%%555
%%%%%%%%%%%%%%%%%%%%%%%%%%%%%%%%%%%%%%%%%%%%%%%%%%%
\section{Other families of patterns: Cartesian  approximation and fibering}
 \label{SFFh}

Application of the Lusternik--Schnirel'man  category theory
to constructing a countable family of solutions of (\ref{S2}) is
explained in \cite[\S~3]{GMPSob}.
 This allowed us to detect the so-called  {\em basic family} of
 patterns $\{F_l\}$, which has been  shown for $N=1$ in a number
of figures in \cite[\S~4]{GMPSob}.
 For convenience, we restate the Lusternik--Schnirel'man/fibering result in \cite[\S~3.2]{GMPSob}:

\begin{proposition}
 \label{Pr.MM}
The elliptic problem $(\ref{S2NN})$ has at least a countable set
of different solutions denoted by $\{F_l, \, l \ge 0\}$, each one
$F_l$ obtained as a critical point of the functional \beq
 \label{V1}
  \mbox{$
 {E}(F)= -  {\displaystyle \frac 12\int_{B_R} |\tilde D^m F|^2 + \frac 12  \int_{B_R}
 F^2 -\frac{1}{\b} \,\int_{B_R}  |F|^{\b}, \quad
 %% \mbox{with} \,\,\,
 \b=\frac {n+2}{n+1}\in (1,2)},
 \,\,\, n>0,
  %% \,\,
  %%%\b= \frac 1 n,
  $}
  \eeq
%%$(\ref{V1})$
 in $W^{2,m}_{0}(B_R)$ in a ball $B_R$ with a sufficiently large
 radius
$R=R(l)>0$.
 \end{proposition}

%%%%%%%%%%%%%%%%%%%%%%%%%%%%%%%%%%%%%%%%%%%%%%%%%%%%%%%%%%%%
\subsection{Basic computations}

We next develop approaches for obtaining other patterns, which are
not detected in Proposition \ref{Pr.MM} by Lusternik--Schnirel'man and classical  fibering
 techniques.

%%%%\smallskip

\begin{figure}
 %%\vskip -.3cm
\centering \subfigure[$\s=\{+6,2,+2,2,+6\}$]{
\includegraphics[scale=0.52]{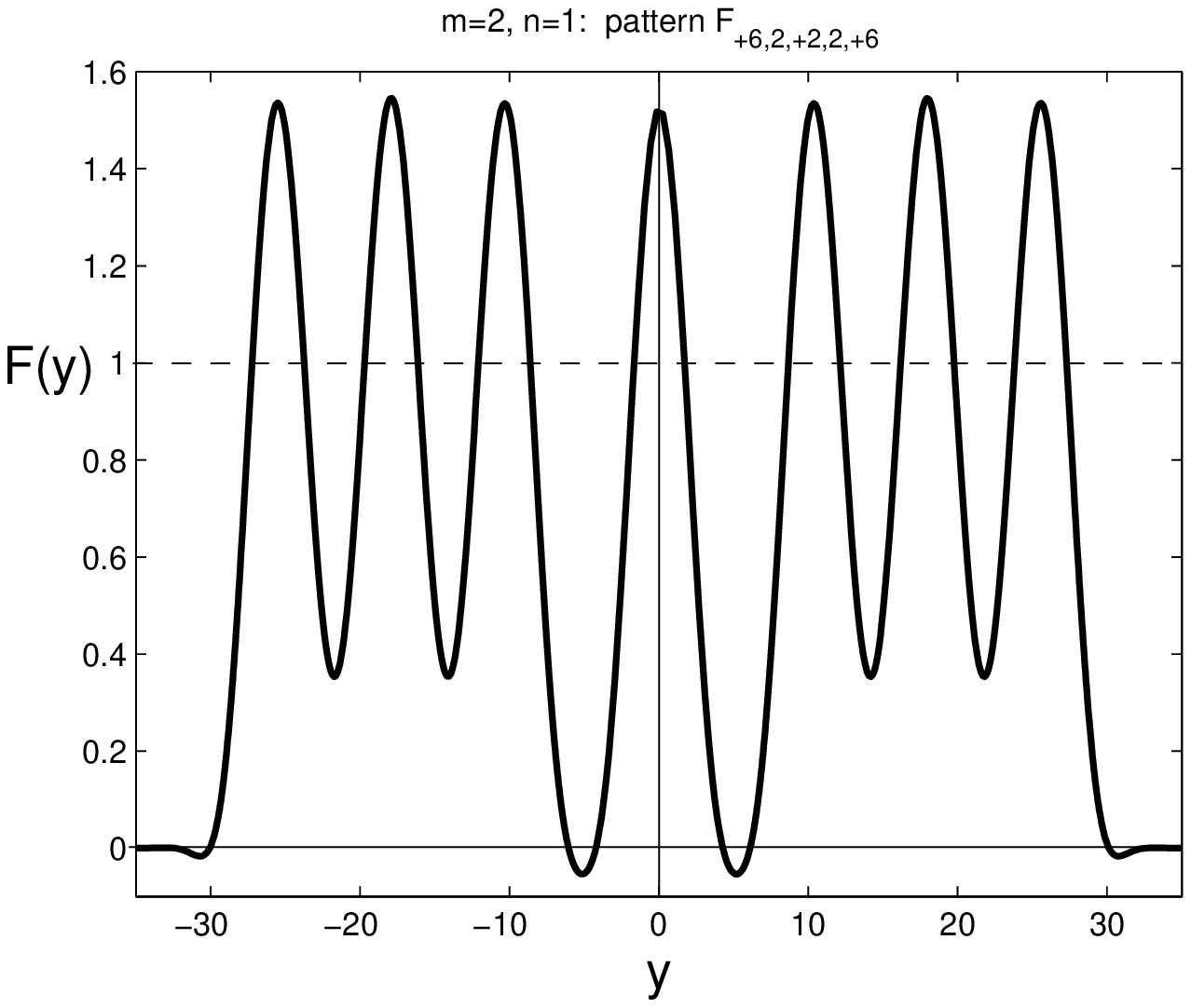}
} \centering \subfigure[$\s=\{+6,2,+4,1,-2,1+2\}$]{
\includegraphics[scale=0.52]{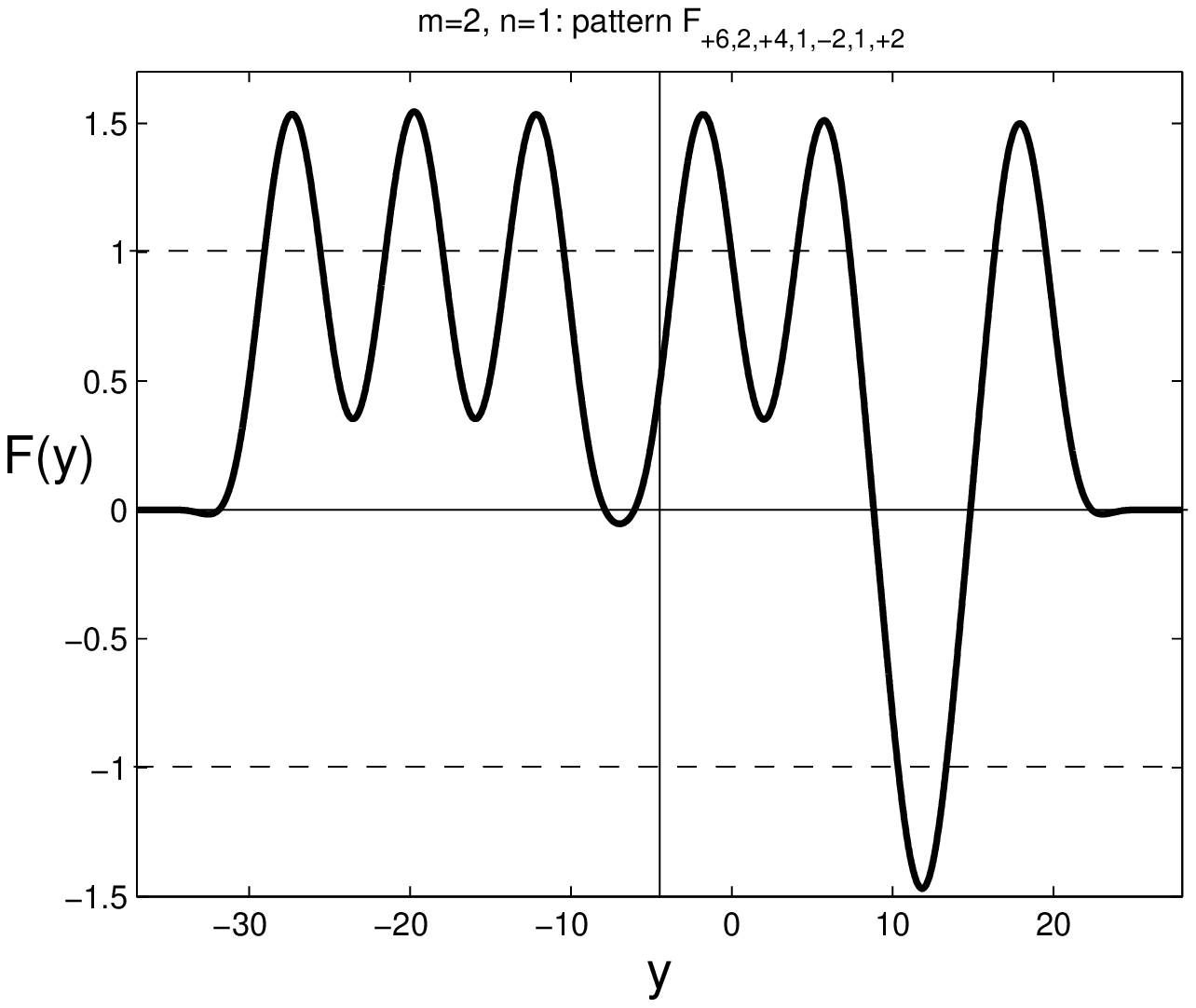}
} \centering \subfigure[$\s=\{+2,2,+4,2,+2,1,-4\}$]{
\includegraphics[scale=0.52]{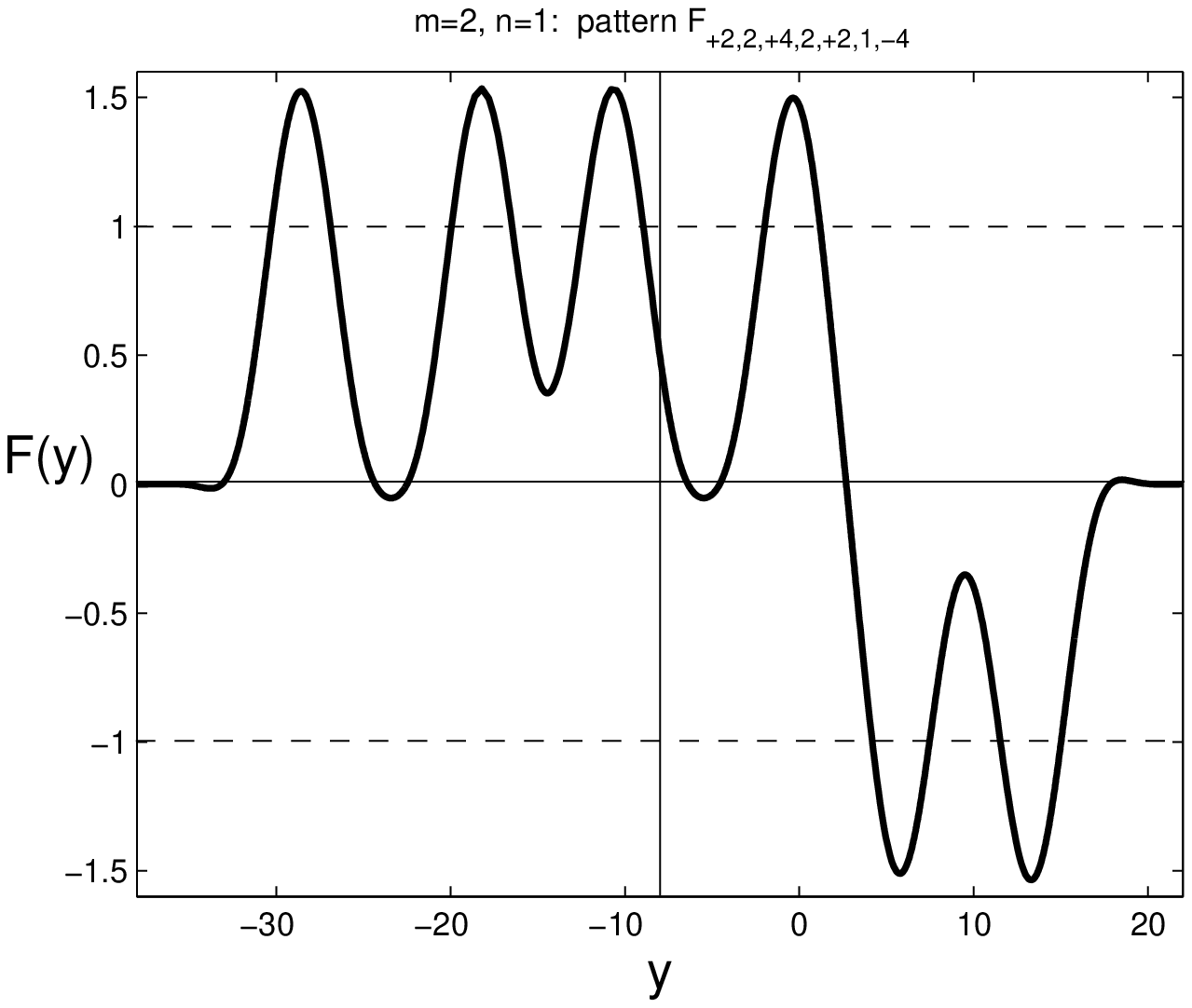}
} \centering \subfigure[$\s=\{+6,3,-4,2,-6\}$]{
\includegraphics[scale=0.52]{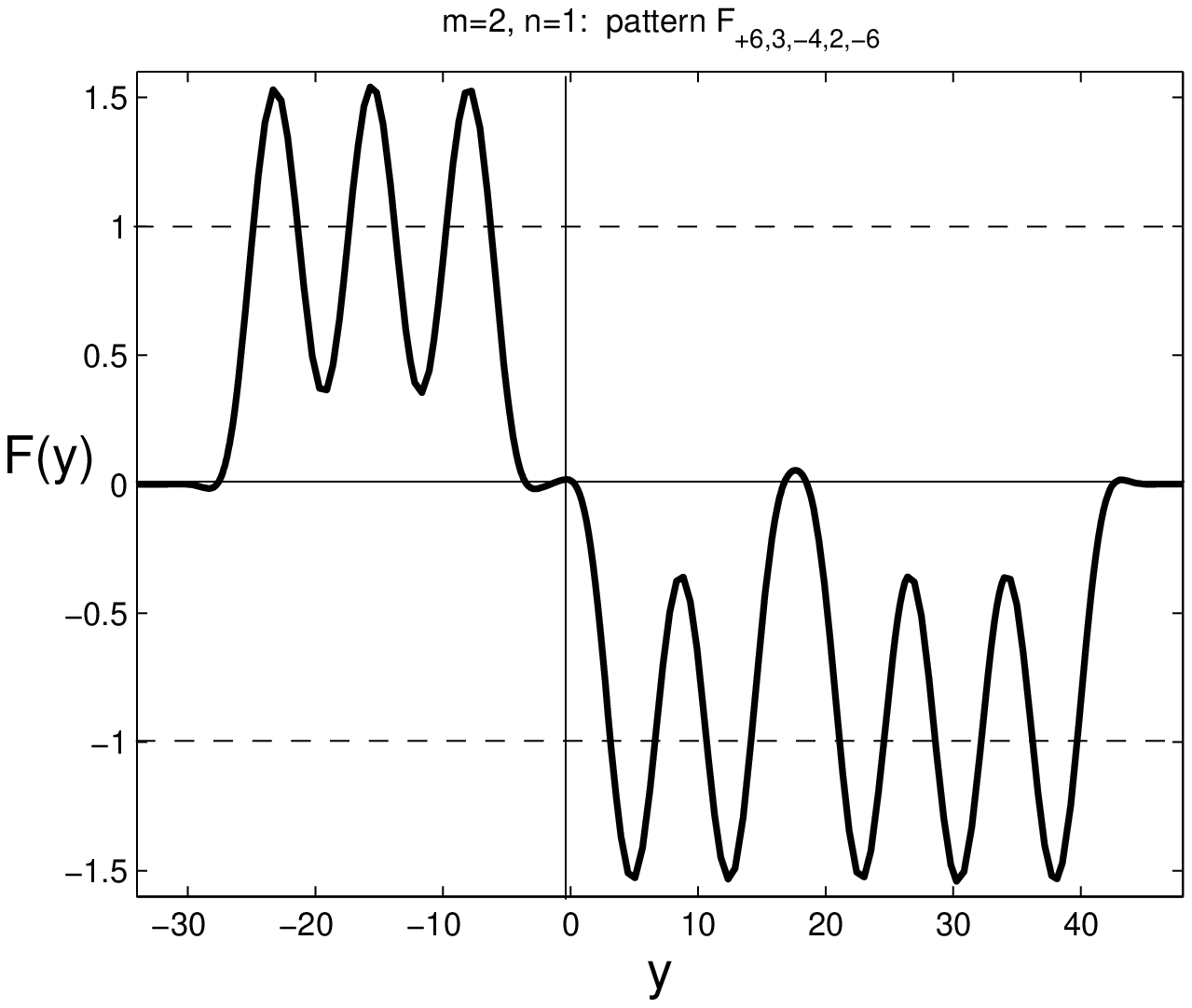}
}
 \vskip -.2cm
\caption{\rm\small  Various patterns for $m=2$ and  $n=1$.}
 \label{FC1}
  %%%%%{FGF.1fig}
\end{figure}

 In order to construct other
families of solutions (see Section \ref{Sectm4} for illustrations
of those for $m=3$),
%% such as shown in Figures \ref{G6},
%%\ref{G7}, \ref{G8}, and  \ref{FC1}, to say nothing about the
%%complicated pattern in Figure \ref{FC2},
 we need an auxiliary
approximation of patterns.
 %%% we are talking about.
Namely, we first perform the Cartesian decomposition
 \beq
 \label{V5}
 F= h+w,
  \eeq
  where $h \in W^{2,m}_{0}(B_R)$ is a  smooth ``step-like
  function" that takes the equilibrium values $\pm1$ and 0
   %%%(we
  %%%%assume $n=1$ for simplicity)
    on some disjoint subsets of $B_R$ (with a smooth connections in between).
  Sufficiently close to the boundary points, we always have $h(y)=0$.
    For instance, in one dimension  for getting the patterns in
   Figure \ref{FC1}, we take $h(y)$ as a smooth approximation of
   the step function, which takes values $\pm 1$ and 0 on the
   intervals of oscillations of the solution about these equilibria.
%%The ordered numbers of oscillations on the successive intervals
%%about equilibria will compose the {\em total Sturm index} ({\em
%%genus}) of the pattern under consideration, as suggested in
%%(\ref{mm1}).

%%\com{To SIP: is it good to say ``genus" for a function?}

%%%\com{Or keep just ``Sturm index"?}

%%%%%%%%%%%%%%%%%%%%%%%%%%%%%%%%%%%%%%%%%%%%%%%%%%%

In other words, we are going to perform the radial fibering not
about the origin  but about the non-trivial point $h$, which plays
the role of an initial approximation of the pattern that we are
interested in. Obviously, the choice of such $h$'s is of principal
importance, which thus should be done very carefully.

Substituting (\ref{V5}) into the functional yields the new one,
 \beq
 \label{V6}
 \begin{matrix}
\hat E(w) = E(h+w)= -  \displaystyle \frac 12  \int_{B_R} |\tilde D^m h|^2 + \frac 12
\displaystyle \int_{B_R}
 h^2+ L_0(h) w
 \qquad \quad
\ssk\ssk\\
 - H_0(w)
  %%%\, \frac 12  \displaystyle \int |\tilde D^m w|^2 + \frac 12 \displaystyle \int
  %%%w^2
  - \displaystyle {\frac{1}{\b} \, \int_{B_R} |h+w|^{\b}},\qquad\quad
   \end{matrix}
    \eeq
    where by $L_0$ we denote the linear functional
     $$
      \mbox{$
     L_0(h) w =- \displaystyle \int_{B_R} \tilde D^m h \cdot \tilde D^m w + \displaystyle \int_{B_R} h w.
      $}
     $$
We next apply the fibering approach by setting, as usual,
 \beq
 \label{V7r}
 w= r(v) v, \quad v \in {\mathcal H}_0, \quad \mbox{whence}
  \eeq
  %%to get the functional
   \beq
   \label{V7}
 \begin{matrix}
\hat H(r,v) = E(h+r(v)v)= -  \displaystyle  \frac 12 \int_{B_R} |\tilde D^m h|^2 +
\frac 12 \displaystyle \int_{B_R}
 h^2+ r \, L_0 v
 %%%\bigl(- \displaystyle \int \tilde D^m h \cdot \tilde D^m v + \displaystyle \int h v \bigr)
 \qquad \quad
\ssk\ssk\\ +\,\displaystyle{ \frac 12 \, r^2
  -\frac{1}{\b} \,} \displaystyle \int_{B_R}  |h+ r v|^{\b}.\qquad\quad
   \end{matrix}
    \eeq
In order to find the absolute minimum point, we need to solve the
scalar equation $\hat H'_r=0$, %%% i.e.,
 \beq
 \label{V8}
  \mbox{$
\hat H_r(r,v) \equiv r-  \displaystyle \int_{B_R} |h+ r v|^{\b-2}(h+ r v)v + L_0(h) v
= 0.
  $}
 \eeq
 For $h=0$, this coincides with the standard equation derived in
 \cite[\S~3]{GMPSob},
 %%(\ref{f31}),
   and
 has three  roots, $r_0(v)=0$ and
  \beq
  \label{V81}
  r_\pm= r_\pm (v),
   \eeq
   which are  positive and negative respectively.
For $h \not = 0$, these roots exist and are slightly deformed for
sufficiently small $h$. For large $h$, one of the roots $r_\pm(v)$
may disappear, and at this instance the resulting functional ({\em
q.v.} below) may loose its smoothness. To distinguish the roots,
%%%Namely,
we observe that
 \beq
 \label{V82}
 \begin{matrix}
  r_-(v)<0 \,\,\, \mbox{always exists and smooth if\,
   $J(h,v)\equiv -\displaystyle \int_{B_R} |h|^{\b-2}h v +L_0(h) v > 0$},\quad\ssk\ssk \\
 r_+(v)>0 \,\,\, \mbox{always exists and smooth if \,
  $J(h,v)\equiv -\displaystyle \int_{B_R} |h|^{\b-2}h v +L_0(h) v < 0$}.\quad
  \end{matrix}
 \eeq

 %%%\com{To SIP: does it work or not ?}

 %%\com{To SIP: can this be applied to $\D u + u^3=h(x)$? Then
 %%existence of smooth branches depends on the sign of $\displaystyle \int hv$
 %%only, and always exist. So Lusternik--Schnirel'man applies for any $p<p_s$ to give a
 %% countable set of solutions for any $h \in L^2$...}

Thus, calculating the extremum point $r_\pm(v)$ (when  exists)
from (\ref{V8}) and substituting into (\ref{V7}) yields the new
functional
 \beq
 \label{V9}
  \tilde H(v)= \hat H(r_\pm(v),v) \quad \mbox{on} \quad {\mathcal
  H}_0, \quad \mbox{being even, since} \quad
  %% \eeq
%%which is even, since
%% $$
  r_\pm(-v)=-r_\pm(v).
 \eeq
  %% $$
Therefore, if (\ref{V9}) is smooth on an appropriate branch
$r_\pm$, this gives a  set of critical point $\{v_k\}$, as above.
    Moreover, in a
neighbourhood of any critical points $v_k$ satisfying
  (\ref{V82}), i.e., $J(h,v_k)>0$ or $<0$, the corresponding branches $r_{\pm}(v)$ are smooth
  for $v \approx v_k$ (and hence along a minimizing sequence $\{v_k^{j}\} \to v_k$), so that
(\ref{V9}) is sufficiently regular. Even in the  delicate case,
when $J(h,v_k)=0$ and $r=0$ is an inflection point of $\hat
H_r(r,v_k)$, i.e.,
 $$
 \hat H_{rr}(r,v_k)=0,
  $$
  one can choose a smooth, existing, and ``stable" branch  $r=r_0(v)$ for $v \approx v_k$ along
  a suitable minimizing sequence.

 %% in a neighbourhood of
%% critical
 %%points under consideration.
 This provides us with a finite number
of critical points  associated with the category of ${\mathcal
H}_0$. As we have seen, all these  critical points of (\ref{V6})
such as
 $$
- {\rm sign} \,\{ L_0(h)w\} =  {\rm sign} \, \{L_0(h)v\} = \pm 1
\quad(\mbox{or} \,\,\, 0)
  $$
  can be obtained on the branch $r_\pm(v)$
  (or $r_0(v)$).
  %% and others, vice versa,
 %%on
 %%% $r_+(v)$.

 %%  \com{To SIP: more explanations how to get smooth branches?}

 %%%%\com{To SIP: a better alternative approach ?  Simpler ?}

Actually, using this mixture (\ref{V5}) of the Cartesian and
spherical fibering decomposition of the functional space, we are
interested, mainly,  in the first critical point, which is defined
via  the absolute infimum of the functional (\ref{V9}) (roughly
speaking, in the case of genus 1). With a  choice of a
sufficiently ``large"
 approximating function $h$, this first pattern will be
different from other basic patterns constructed above for $h=0$.
Indeed, this first pattern is characterized by the condition of
the ``minimal deviation" from $h(y)$, while, e.g., $F_0$
corresponds to the minimal deviation from $h= 0$, so that these
cannot coincide if $h$ is large enough and has a proper shape
concentrating about equilibria $\pm 1$ and 0.

\begin{figure}
 \centering
 \psfrag{F(y)}{$F(y)$}
 \psfrag{F0(y)}{$F_0(y)$}
 \psfrag{h(y)}{$h(y)$}
 \psfrag{F44222(y)}{$F_{+4,-4,+2,-2,+2}(y)$}
 \psfrag{F+4}{$F_{+4}$}
  \psfrag{v(x,t-)}{$v(x,T^-)$}
  \psfrag{final-time profile}{final-time profile}
   \psfrag{tapp1}{$t \approx 1^-$}
\psfrag{x}{$x$}
 \psfrag{0<t1<t2<t3<t4}{$0<t_1<t_2<t_3<t_4$}
 %% \psfrag{0}{$0$}
 \psfrag{y}{$y$}
 %%\psfrag{-1}{$-1$}
 %%\psfrag{1}{$1$}
 %%\psfrag{-1}{$-1$}
\includegraphics[scale=0.5]{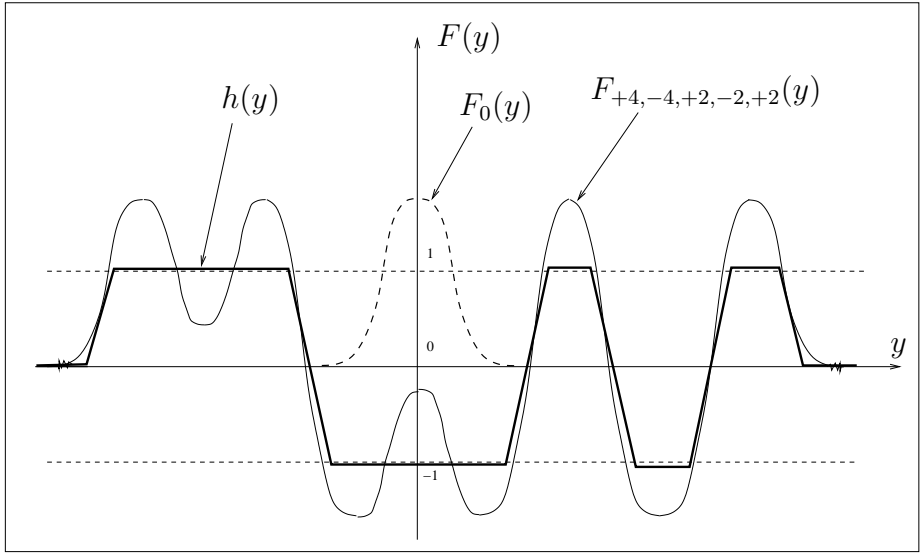}     %%%{F1EE1.pdf}
\caption{\small For getting $F_{+4,-4,+2,-2,+2}$, the Cartesian approximation $h$
in (\ref{V5}) should be chosen properly.}
     \vskip -.3cm
 \label{F42}
\end{figure}

  Figure \ref{F42} illustrates such a statement and
shows a typical Cartesian approximation $h$, which is necessary to
detect the patterns $F_{+4,-4,+2,-2,+2}$. Obviously, then the
absolute extremum of $\tilde H (v)$ cannot be attained at already
known critical point $F_0$ given by the dashed line, which is
characterized by a much larger deviation from the fixed $h$ that
is given by a boldface line (it should be slightly smoothed at
corner points).

%%%%%%%%%%%%%%%%%%%%%%%%%%%%%%%%%%%%%%%%%%%%%%%%%%%%%%%
%%%%%%%%

Therefore, the main result in \cite[\S~3]{GMPSob}, such as
 Proposition
\ref{Pr.MM},
 remains true for any
sufficiently regular initial approximation $h$. Of course, some of
the critical points $F_l(y;h)$ with $l \gg 1$ may coincide with
already known basic patterns $F_l$, but, in fact, we are
interested in the first critical value and point, which thus give
$F$ that has the minimal deviation from $h$, and must be different
from $F_l$'s. Obviously, for approximations $h$ that are far away
from 0, the first pattern $F(y;h)$ obtained by using the sets
${\mathcal H}_0$ of arbitrary category, $\rho \ge 1$, cannot
coincide with the first basic patterns $F_l(y)$, which are
sufficiently small and have a specific and different geometric
structure.

The actual and most general rigorous ``optimal" choice and
characterization of such suitable approximations $h$ (possibly a
sequence of such $\{h_k\}$) remains an open problem, though we
have got a convincing experience in understanding of such
patterns, in particular, using numerical experiments and some
analytic estimates; see  related comments below.

%% \com{To SIP:  branches smooth for a.e. $h$?}

 %%\com{To SIP: how to complete the analysis ?}

 %%%% \com{To SIP: why do we have  sequences $\{h_k\}$?}

 %%%\smallskip

%%%%%%%%%%%%%%%%%%%%%%%%%%%%%%%%%%%%%%%%%%%%%%%%%%%%%%
\subsection{Some asymptotic analysis}

It is easy to show that, asymptotically, for sufficiently
``spatially wide" patterns, the Cartesian-spherical fibering
(\ref{V5}), (\ref{V7r}) provides us with families of patterns that
are different from basic ones $\{F_l\}$.

\begin{figure}
%\vskip -.3cm
 \centering
\includegraphics[scale=0.8]{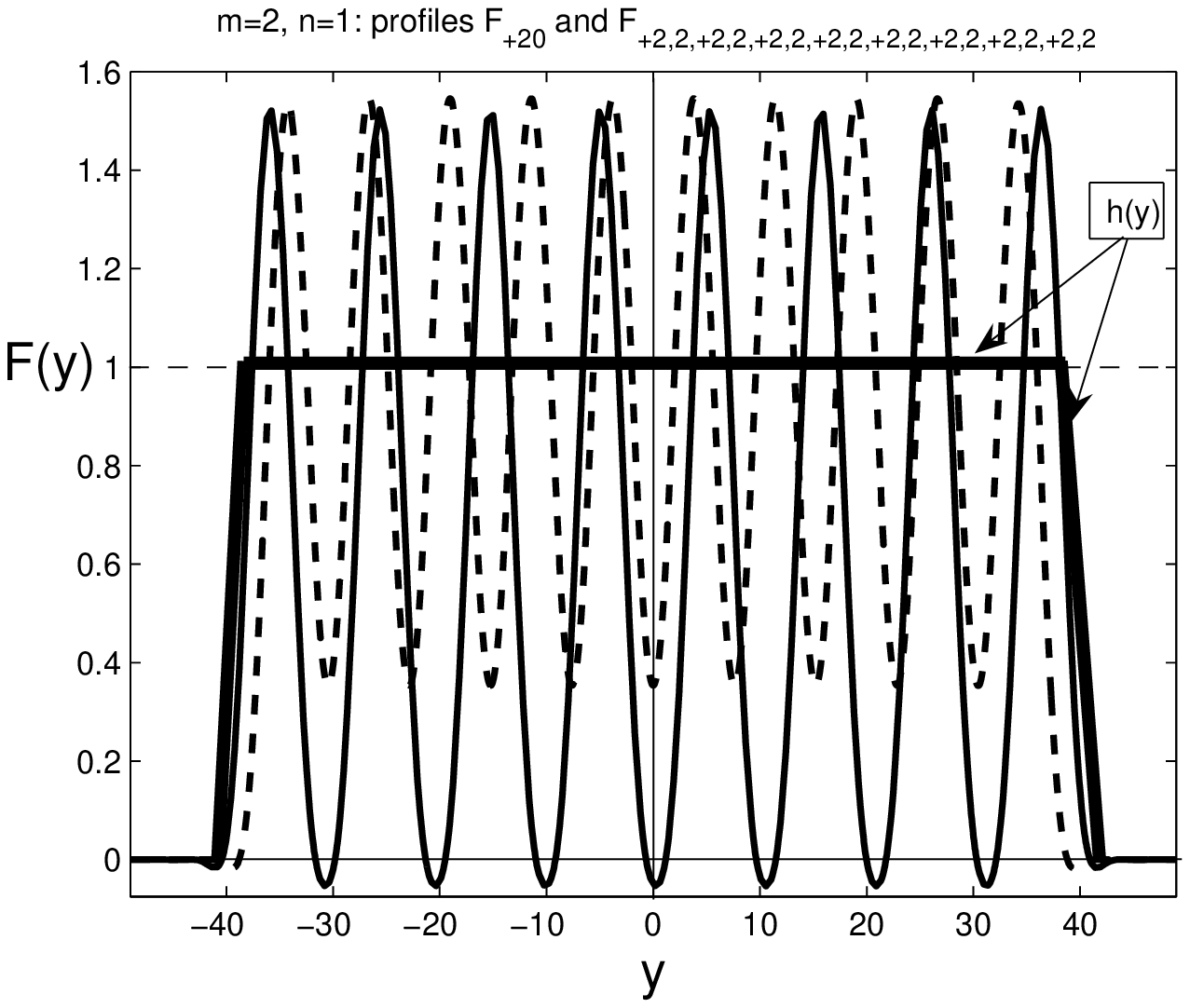}
 \vskip -.4cm
\caption{\rm\small Two patterns, $F_{+20}$ (dashed line) obtained
by Cartesian-spherical fibering and
$F_{+2,2,+2,2,+2,2,+2,2,+2,2,+2,2,+2,2,+2}$ corresponding to
$h=0$; $m=2$, $n=1$.}
   \vskip -.3cm
 \label{F14mm}
\end{figure}

For instance, in Figure \ref{F14mm} we compare the patterns
 $$
 F_{+2k} \,\,(\mbox{the dashed line})  \quad \mbox{and} \quad F_{+2,2,+2,2,...,+2,2,+2}
 \quad
 \mbox{for large $k=10$, i.e.,}
 $$
 %%%for large $k=10$, i.e.,
  $$
   \begin{matrix}
 F_{+20} \whereA \fbox{$c_F=2.9398...$}\, , \ssk\\
  \mbox{and} \quad
 F_{+2,2,+2,2,+2,2,+2,2,+2,2,+2,2,+2,2,+2} \whereA \fbox{$c_F=2.7197...$} \, .
  \end{matrix}
  $$
Here $c_F$ stand for the corresponding critical values of the
functional obtained after fibering \cite[\S~3]{GMPSob,GMPSobIarX}
  \beq
   \label{dd1}
 \mbox{$
 c_F \equiv \tilde H(v) = \frac { \displaystyle \int_{B_R}  |F|^\b}{\left( - \displaystyle \int_{B_R}  |\tilde D^m F|^2 + \displaystyle \int_{B_R}
  F^2 \right)^{ \b/ 2}} \quad \bigl(\b = \frac{n+2}{n+1} \bigr).
   $}
 \eeq

 Thus,  these
two patterns are clearly recognized
 by their different critical values $c_F$ indicated.
For $F_{+20}$, the corresponding functions $h(x) \approx 1$ on
$(-40,40)$ is shown by the boldface line. It is  seen that the
global minimum of the functional (\ref{V9}) for such $h$ {\em
cannot be attained} on any profile from the basic family
$\{F_l\}$, because the total deviation becomes huge in comparison
with the almost periodic deviation achieved via $F_{+20}$. In this
case, the minimum is attained on the profile $F_{+2k}$ having a
completely different geometry.

%%%%%%%%%%%%%%%%%%%%%%%%%%%%%%%%%%%%%%%%%%%%%%%%%%
%%%%%%%%%%%%%%%%%%%%%%%%%%%%%%%%%%%%%%%%%%%%%%%%%%%%%%

For $k \gg 1$, these observations can be fixed in a standard
asymptotically rigorous manner, which we are not going to do here.

%%%%%%%%%%%%%%%%%%%%%%%%%%%%%%%%%%%%%%%%%%%%%%%%%%%%%%%%%%%%%
%%%%%%%%%%%%%%%%%%%%%%%%%%%%%%%%%%%%%%%%%%%%%%%%%%%%%%%%%%%%%%%%%%%%%%
\subsection{The origin of countable sequences of solutions: a formal double
fibering}

 Taking into account both changes (\ref{V5}) and
(\ref{V7r}), we arrive at the functional
 \beq
 \label{l1}
 \hat E(h,r,v) \equiv E(h+r(v)v).
  \eeq
 The relative critical point of (\ref{l1}) are given by the
 system
  \beq
  \label{l2}
  \left\{
  \begin{matrix}
  \hat E_h '(h,r,v)=0, \\
   \hat E_r '(h,r,v)=0, \\
    \hat E_v '(h,r,v)=0.
     \end{matrix}
     \right.
      \eeq
Of course, the first equation is just equivalent to the original
one, since
 $$
 \hat E_h'= E',
 $$
so that (\ref{l2}) is a formal system comprising the spherical
fibering in the $\{r,v\}$-variables, and the original equation.
Let us see what kind of conclusions can be derived from this.

 The second equation  is scalar and gives us necessary smooth
 branches (under certain hypotheses as above)
  \beq
  \label{l3}
  r= r_*(h,v) \quad (r_*(h,-v)=-r_*(h,v)).
   \eeq
   Then we arrive at the system
 \beq
  \label{l31}
  \left\{
  \begin{matrix}
  E_h '(h+r_*(h,v)v)=0, \\
    E_v '(h+r_*(h,v)v)=0.
     \end{matrix}
     \right.
      \eeq

      The first equation is difficult to handle (possibly more difficult
      than the original one). However, assume that this can be
      solved.
 In view of (\ref{l3}), gives  an even dependence,
 \beq
 \label{l4}
 h=h_*(v) \quad (h_*(-v)=h_*(v)).
  \eeq
Finally, we arrive at the even weakly continuous functional for
$v$,
 \beq
 \label{l5}
 E_*(v) \equiv E(h_*(v) + r(h_*(v),v)v) \quad \mbox{in} \quad
 {\mathcal H}_0.
  \eeq
  This has a countable\footnote{As usual, we mean compactly supported solutions in $\ren$,
  i.e., $R=+\iy$.} set of critical point $\{v_k\}$, which by
  fibering method \cite{PohFM}, generate critical point of the
  original functional (\ref{l1}).

Therefore, eventually, as $R \to \infty$, we obtain a countable
set of necessary critical points
 $$
 (h_k,w_k)=((h_*(v_k), r_*(v_k)v_k), \quad k =0,1,2,... \, .
  $$
  The corresponding $F$-patterns are denoted by
 $$
 \{\tilde F_l\}.
 $$
The actual general structure of such special solutions remains
unclear and needs extra analysis. Currently, we know a little
about this and present a few comments only. Using an analogy with
the basic Lusternik--Schnirel'man patterns $\{F_l\}$ obtained in \cite[\S~3]{GMPSob}
%% Section \ref{SectVar}
for $h=0$, it
may be expected that each $\tilde F$ is composed from $l\ge 1$
copies of the ``elementary" profile $F_{+4}$, i.e.,
 $$
 \tilde F_l \sim \{F_{+4},-F_{+4},...,(-1)^{l+1}F_{+4}\},
  %%%%\quad \mbox{for} \quad l \ge 1,
 $$
 with the obvious choice of the corresponding Cartesian
 approximation $h_l(y)$ that is concentrated about equilibria at
 $\pm 1$ and 0 in between.
 %%%$$
 %%%h_l(y) = 1 \quad \mbox{for}
 %%%%\quad |y| \le l y_*
It is more likely that $\{\tilde F_l\}$ includes other profiles of
the
 $\{\pm F_{+4},...,\pm F_{+4}\}$-gluing (see \cite[\S~4]{GMPSob} for definitions), or, in particular,
 can be composed from completely ``non-oscillatory" profiles,
 i.e.,
  $$
  \tilde F_l=F_{+2l}.
   $$

 %%%\com{SIP: to justify a bit more ? Even functional in $v$ ???}

 %%%%\com{SIP: this makes sense or not ? To skip ?}

%%%%%%%%%%%%%%%%%%%%%%%%%%%%%%%%%%%%%%%%%%%%%%%%%%
\section{\underline{\bf Problem ``Numerics"}: patterns in higher-order cases, $m \ge 3$}
\label{Sectm4}

The main features of the pattern  classification by their
structure and computed critical values
 for $m=2$ \cite{GMPSob} can
be extended to arbitrary $m \ge 3$ in the ODEs (\ref{S2}), so we
perform this in less detail.

\begin{figure}
%\vskip -.3cm
 \centering
\includegraphics[scale=0.80]{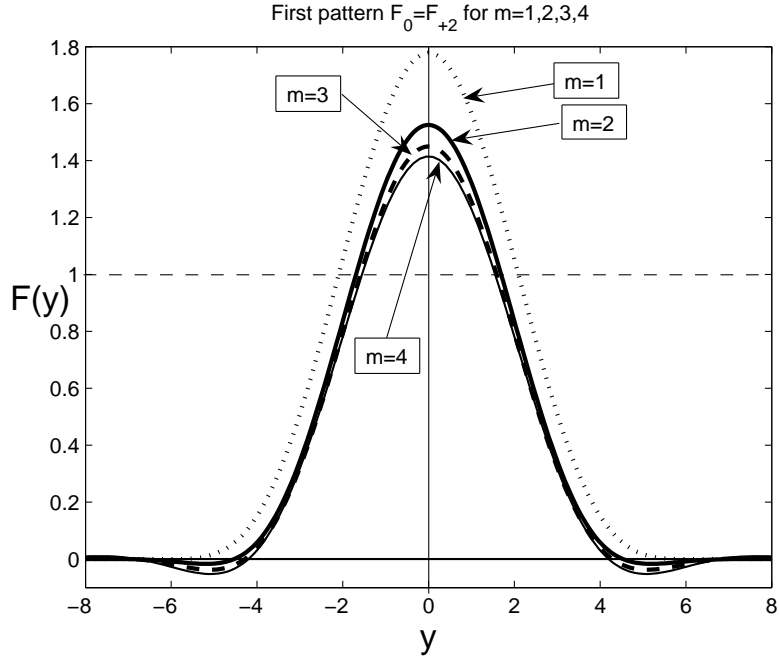}  %%%% 4G1.eps for n=1 only
 \vskip -.4cm
\caption{\rm\small The first solution $F_0(y)$  of  (\ref{S2}),
 $n=1$,
for  $m=1,2,3,4$.}
   \vskip -.3cm
 \label{Gmm1}
\end{figure}

 In Figure
\ref{Gmm1}, for the purpose of comparison, we show the first basic
pattern $F_0(y)$ for $n =1$ in four main cases $m=1$ (the only
non-negative profile by the Maximum Principle known from the 1970s
\cite{SZKM2}, \cite[Ch.~4]{SGKM}),
%% (\ref{RD.4})),
and $m= 2,\, 3,\, 4$. Next
 Figure \ref{Gmm2}  explains
 oscillatory properties of such $F_0(y)$ close to the interface points. It turns out
that, for $m=4$, the solutions are most oscillatory, so  it is
convenient to use this case  for further illustrations.

 %%FIG%%%%%%%%%%%%%%%%%%%%%%%

\begin{figure}
 %%\vskip -.3cm
\centering \subfigure[scale $10^{-2}$]{
\includegraphics[scale=0.52]{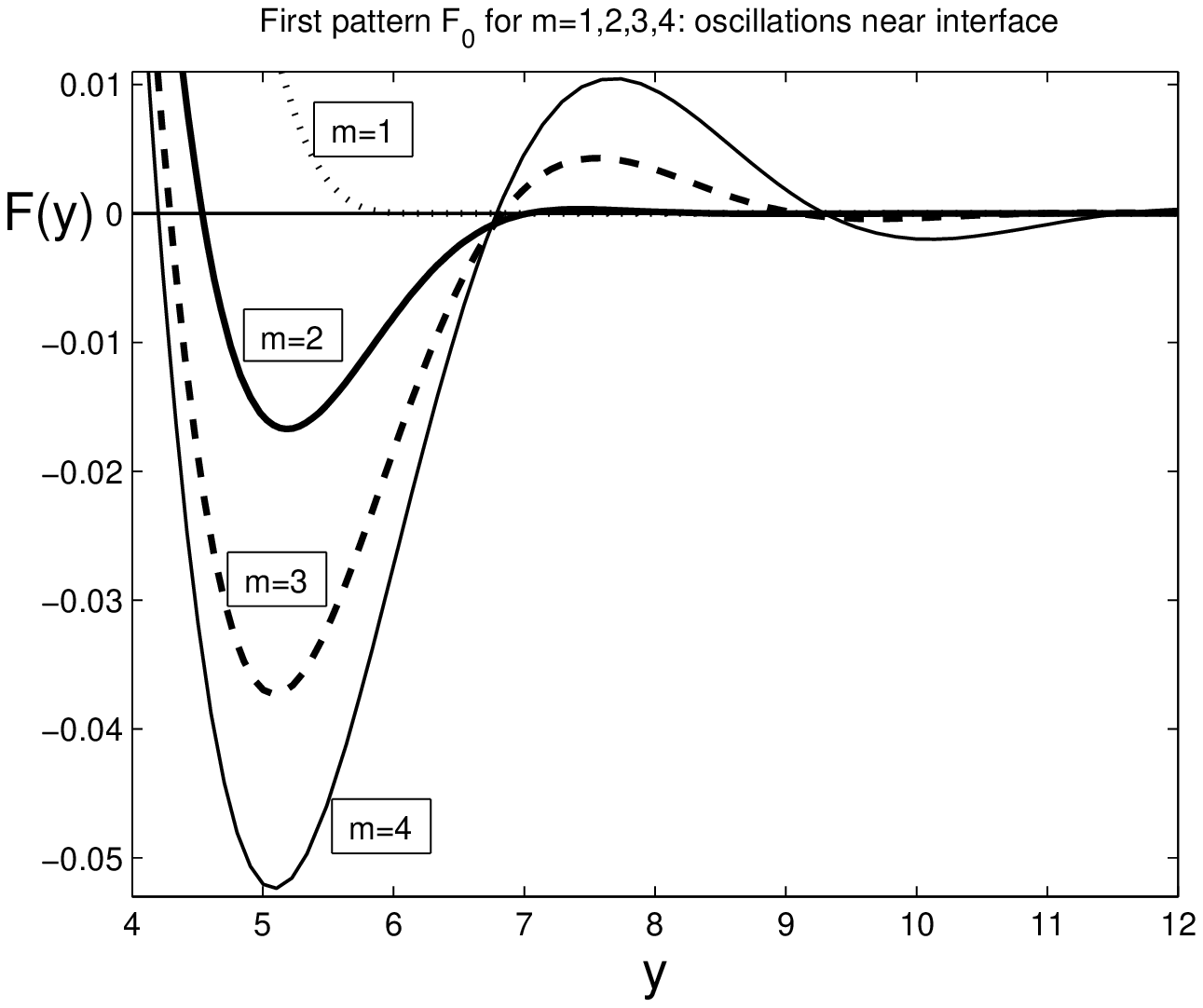}
} \subfigure[scale $10^{-3}$]{
\includegraphics[scale=0.52]{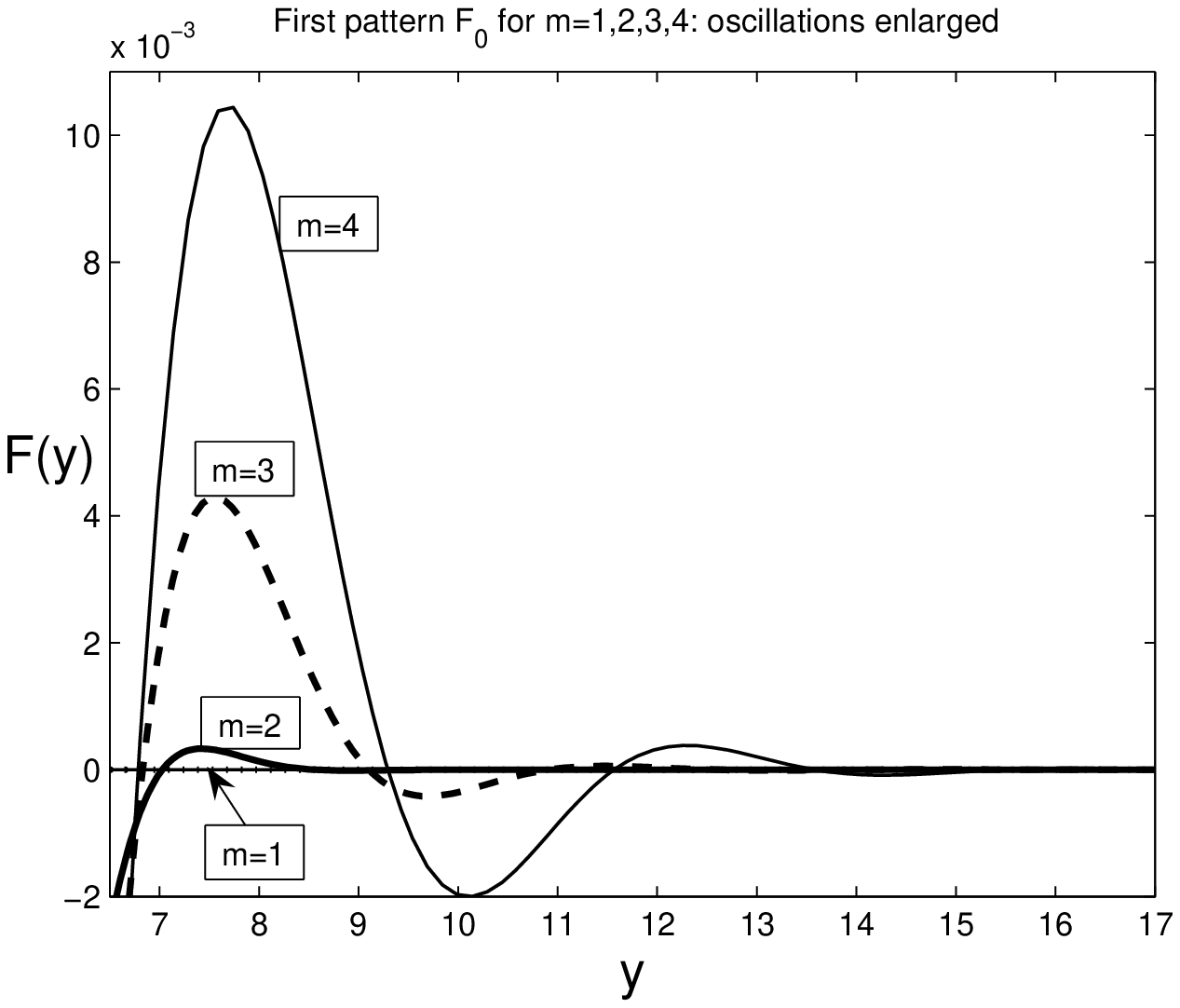}               %%%{FFP01NN.eps}
}
 \vskip -.2cm
\caption{\rm\small  Enlarged zero structure of the profile
$F_0(y)$ for $n=1$ from Figure \ref{Gmm1}; the linear scale.}
 \label{Gmm2}
  %%%%%{FGF.1fig}
\end{figure}

%%\begin{figure}
%\vskip -.3cm
%% \centering
%%\includegraphics[scale=0.65]{4G4.eps}
%% \vskip -.4cm
%%\caption{\rm\small The zero structure of the profile $F_0(y)$ for
%%$n=1$ in the log-scale.}
%%   \vskip -.3cm
%% \label{G3}
%%\end{figure}

%%%%%%%%%%%%%%%%%%%%%%%%%%%%%%%%%%%%%%%%%%%%%%%%%%%%%%%%%%%%%%%%%%%%%

\begin{figure}
%\vskip -.3cm
 \centering
\includegraphics[scale=0.8]{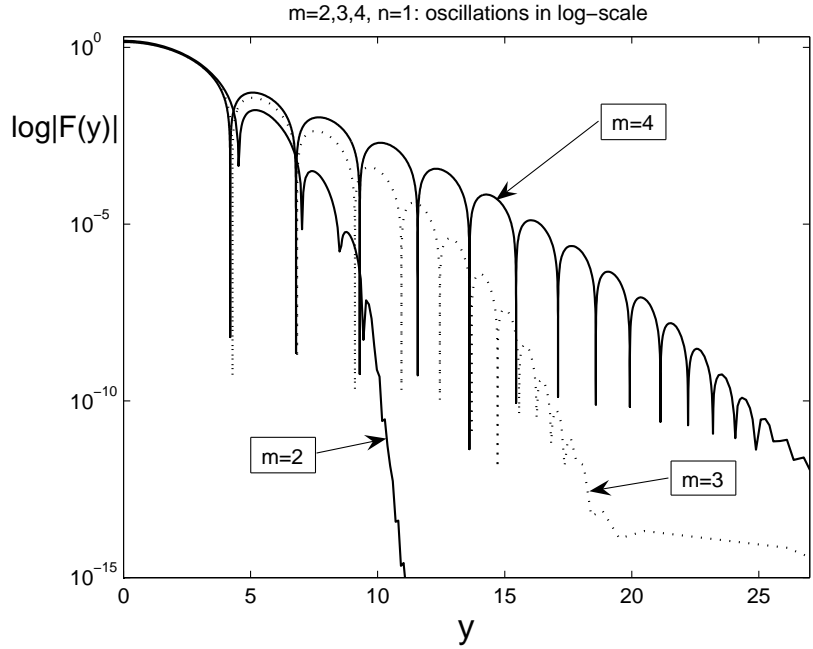}  %%%% 4G1.eps for n=1 only
 \vskip -.4cm
\caption{\rm\small Behaviour of $F_0(y)$  for
 $n=1$,
for  $m=2,3,4$; the log-scale.}
   \vskip -.3cm
 \label{ZZ1}
\end{figure}

In the log-scale, the zero structure of $F_0(y)$ near interfaces
is shown in Figure \ref{ZZ1} for $m=2,3$, and 4 ($n=1$). For $m=4$
and $m=3$, this makes it possible to observe a dozen of
``nonlinear" oscillations that well correspond to the already
known  oscillatory component structure
%%% (see analytic formulae in
%%%\cite[\S~4]{GMPSob}
%% (\ref{2.2})
 close to
interfaces; see \cite[\S~4]{GMPSob}. For the less oscillatory case
$m=2$, we observe 4 reliable oscillations up to $10^{-10}$, which
is our best accuracy achieved.

 %%FIG%%%%%%%%%%%%%%%%%%%%%%%

The basic countable family satisfying approximate Sturm's property
has the same topology as for $m=2$ \cite[\S~5]{GMPSob},
%%in Figure \ref{G4},
 and we do
not present such numerical illustrations.

\begin{figure}
%\vskip -.3cm
 \centering
\includegraphics[scale=0.8]{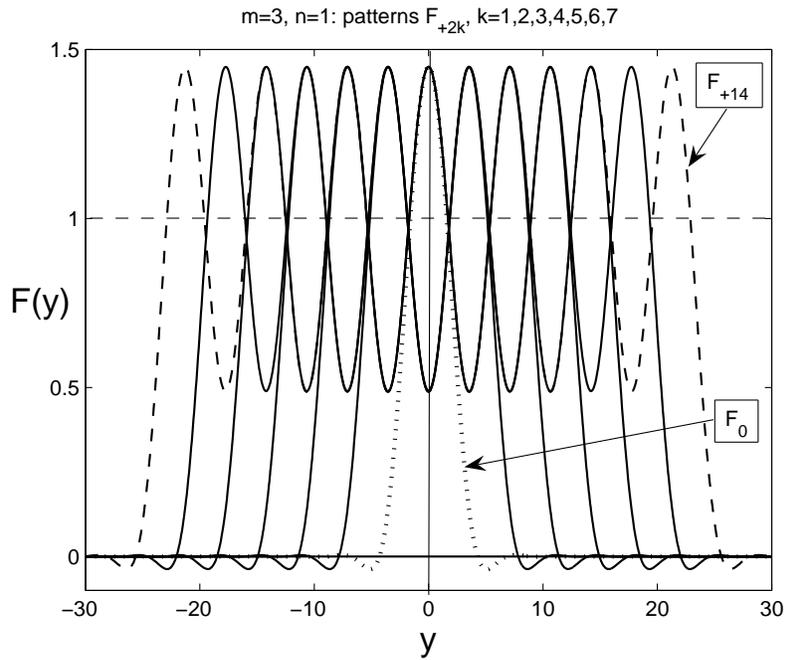}
 \vskip -.4cm
\caption{\rm\small The first seven patterns from the family
$\{F_{+2k}\}$;
 %% of the  $\{F_0,F_0\}$-interactions;
 $m=3$ and
$n=1$.}
   \vskip -.3cm
 \label{Gmm6}
\end{figure}

\begin{figure}
%\vskip -.3cm
 \centering
\includegraphics[scale=0.85]{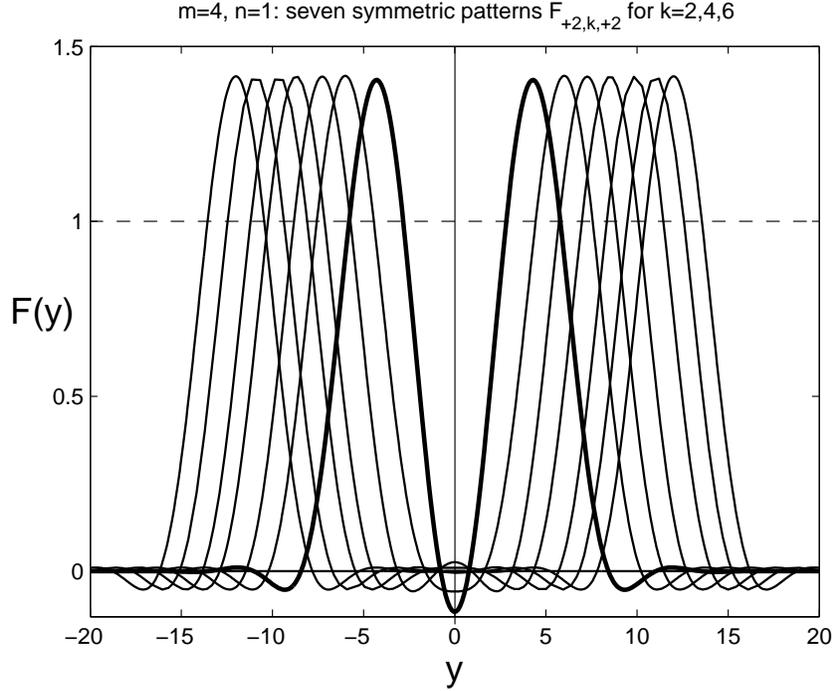}
 \vskip -.4cm
\caption{\rm\small The first  patterns from the family
$\{F_{+2,k,+2}\}$
 of the  $\{F_0,F_0\}$-interactions;
  $m=4$ and $n=1$.}
   \vskip -.3cm
 \label{Gmm68}
\end{figure}

\begin{figure}
 %%\vskip -.3cm
\centering \subfigure[profiles]{
\includegraphics[scale=0.7]{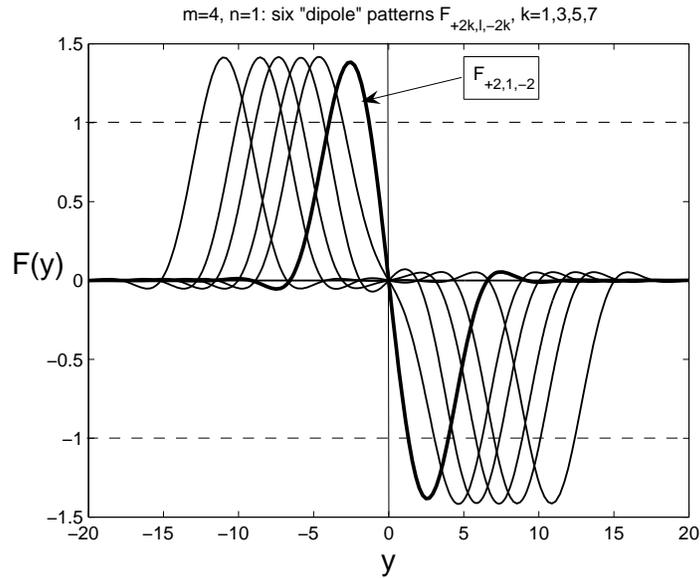}
} \subfigure[zero structure on $(-8,8)$]{
\includegraphics[scale=0.7]{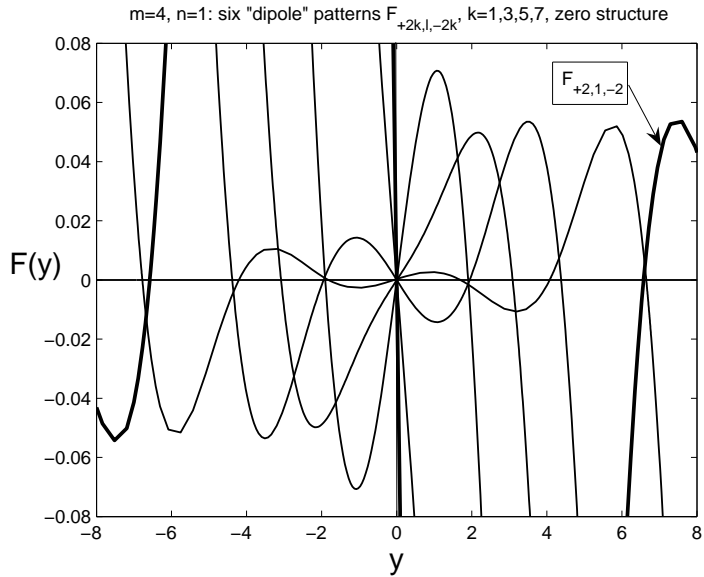}
}
 \vskip -.2cm
\caption{\rm\small The first patterns from the family
$\{F_{+2,k,-2}\}$ of the  $\{-F_0,F_0\}$-interactions, for $m=4$
and $n=1$: profiles (a), and zero structure around $y=0$ (b).}
 \label{Gmm7}
  %%%%%{FGF.1fig}
\end{figure}

In Figure \ref{Gmm6} for $m=3$ and  $n=1$, we show the first
profiles from the family  $\{F_{+2k}\}$, while Figure \ref{Gmm68}
explains typical structures of $F_{+2,k,+2}\}$ for $m=4$, $n=1$.
In Figure \ref{Gmm7} for $m=4$ and $n=1$, we show the first
profiles from the family $\{F_{+2,k,-2}\}$.

%%%%%%%%%%%%%%%%%%%%%%%%%%%%%%%%%%%%%%%%%%%%%%%%%%%%%%%%%%
%%%%%%%%%%%%%%%%%%%%%%%%%%%%%%%%%%%%%%%%%%%%%%%%%%%%%%%%%%%%%%%%

%%%%%%%%%%%%%%%%%%%%%%%%%%%%%%%%%%%%%%%%%%%%%%%%%%%%%%%%%%
%%%%%%%%%%%%%%%%%%%%%%%%%%%%%%%%%%%%%%%%%%%%%%%%%%%%%%%%%%%%%%%%

\begin{figure}
%\vskip -.3cm
 \centering
\includegraphics[scale=0.8]{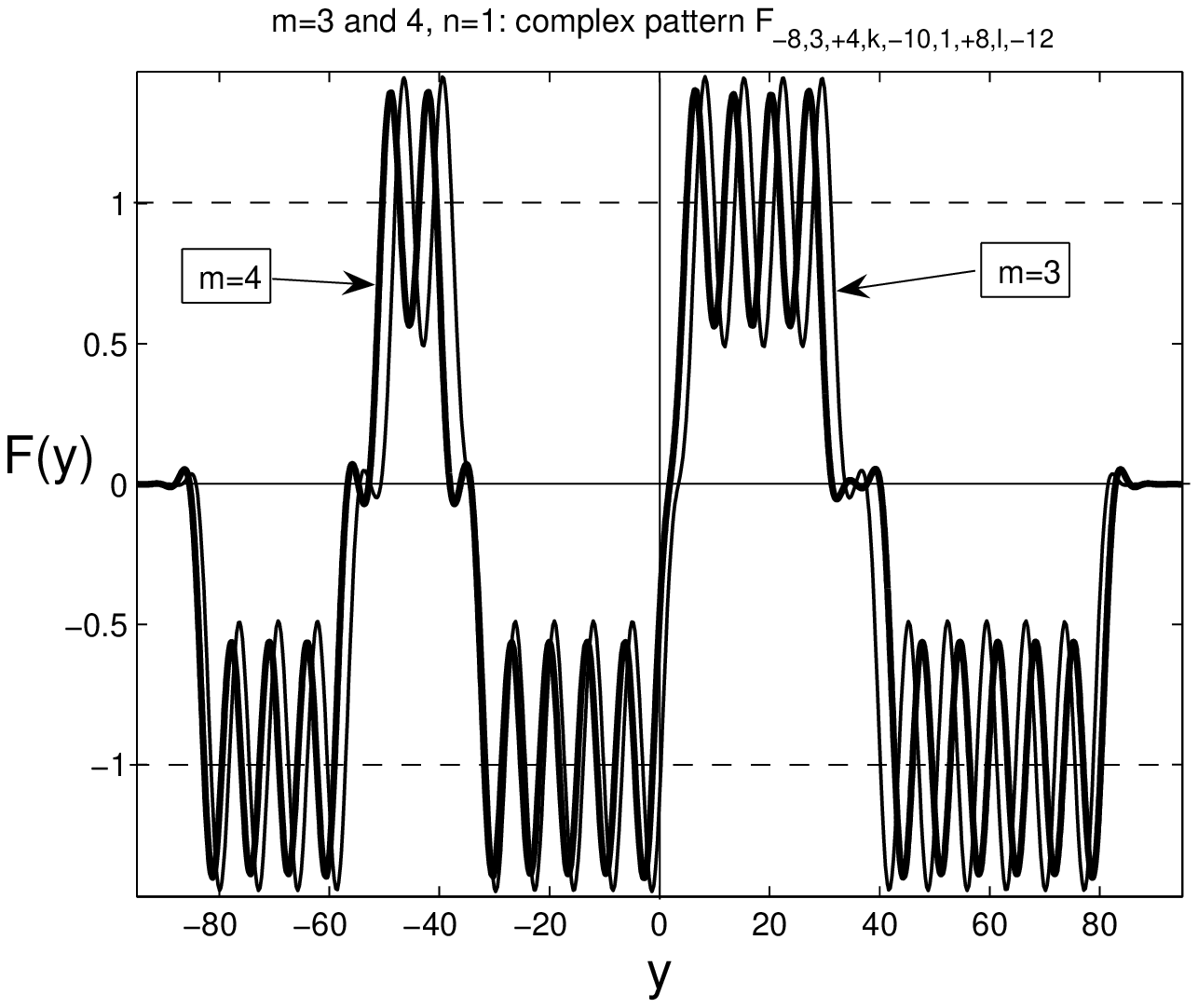}
 \vskip -.4cm
\caption{\rm\small A complicated pattern $F_{\s}(y)$  for $m=3,4$
and $n=1$.}
 %%%%  \vskip -.3cm
 \label{FCmm2}
\end{figure}

 Finally, in
  Figure \ref{FCmm2}, for comparison, we present  a complicated pattern
   for $m=3$ and $4$ (the bold line), $n=1$, with the index
   \beq
   \label{ll1}
   \s=\{-8,3,+4,k,-10,1,+8,l,-12\}.
    %% \whereA k=3,5.
    \eeq
Both numerical experiments were performed starting with the same
initial data. As a result, we obtain quite similar patterns, with
the only difference that, in (\ref{ll1}), $k=1$, $l=3$ for $m=3$,
and for more oscillatory case $m=4$, the number of zeros increase,
so now $k=3$ and  $l=5$.

%%%%%%%%%%%%%%%%%%%%%%%%%%%%%%%%%%%%%%%%%%%%%%%%%%
%%%%%%%%%%%%%%%%%%%%%%%%%%%%%%%%%%%%%%%%%%%%%%%%%%%%%%

%%%%%%%%%%%%%%%%%%%%%%%%%%%%%%%%%%%%%%%%%%%%%%%%%%%%
%%%%%%%%%%%%%%%%%%%%%%%%%%%%%%%%%%%%%%%%%%%%%%%%%%%%%%%%%%
\section{\underline{\bf Problem ``Sturm Index"}: an homotopy classification of
 patterns via $\e$-regularization (analytic-numerical approach)}
 \label{SectHom}

%%%%%%%%%%%%%%%%%%%%%%%%%%%%%%%%%%%%%%%
\subsection{Sturm index for second-order ODEs: in need to extension
for $m \ge 2$}

As we have mentioned, it is well-known that
 in the second-order case $m=1$, solutions of ODE problems,
 even in the non-Lipschitz case,
 %%% (\ref{4.4}),
 \beq
 \label{4.4}
  \mbox{$
 F''=-F + \big|F \big|^{-\frac n{n+1}}F \inB \re
 $}
  \eeq
 obey Sturm's Theorem on zeros, which is a corollary of the
Maximum Principle. Namely, concerning the problem (\ref{4.4}),
each function $F_l(y)$ from the basic family has precisely $l$
isolated zeros (sign changes) and $l+1$ non-degenerate extremum
points. Therefore, the {\em Sturmian index} $I_l=l$ of $F_l$, as
the number of its ``transversal zeros", uniquely specifies any of
the basic patterns. This is also equivalent to the {\em Morse
index} of the corresponding linearized operator. Moreover, the
Lusternik--Schnirel'man minimax construction of critical points reveals this zero
structure of the minimizers \cite[p.~385]{KrasZ} that is directly
associated with the category \cite[\S~6.6]{Berger}, or genus
\cite[\S~57]{KrasZ}, of the sets involved in the variational
construction.

\ssk

For any $m \ge 2$, the Maximum Principle fails, and such a
rigorous geometric classification of basic patterns is no longer
available in view of existence of oscillatory tails close to both
interfaces. Roughly speaking, each profile obtained by the
Lusternik--Schnirel'man/fibering approach has infinitely many zeros and extremum
points, which  makes it impossible to use the above  simple
geometric characteristics for classification of the patterns as
for $m=1$.

Nevertheless, we claim that a Sturmian-type characterization of
{\em some} (basic) patterns is possible for oscillatory solutions
of higher-order ODEs.  We will
 reveal how to attach the Sturmian
index to solutions from the basic family $\{F_l\}$ in the
higher-order case, and also to other families. We consider the ODE
case (\ref{S2}) for $N=1$, though a similar approach applies to
the radial elliptic setting in (\ref{S2NN}), as well as
non-radial, where though it is not that well-presented and clear.

We begin with description of higher-order equations admitting a
rigorous Sturmian classification of patterns.

%%%%%%%%%%%%%%%%%%%%%%%%%%%%%%%%%%%%%%%%%%%%%%%%%%%%%%%%
\subsection{$2m$th-order equations with Sturmian ordering}

Consider the following ODE with Dirichlet boundary conditions:
 \beq
 \label{m1}
 (-1)^{m+1} F^{(2m)} + F^3 =0 \,\,\, \mbox{in}
  \,\,\, (-R,R), \quad F=F'=...=F^{(m-1)}=0 \,\,\,
  \mbox{at} \,\,\, y= \pm R.
   \eeq
This problem is also variational and admits a countable set of
solutions $\{F_l\}$. Moreover, since the differential operator
 $
 (-1)^{m} D_y^{2m}
   $
is an iteration of the positive operator $-D^2_x$ with the Maximum
Principle,
    according to Elias \cite{Elias} (see also applications to some nonlinear
higher-order eigenvalue problems in \cite{Rynn1, Rynn2}),
 the following result holds:

\begin{proposition}
\label{PrSt33}
 The $l$th-solution $F_l(y)$ of the
problem $(\ref{m1})$ for any $l=0,1,2,...$, has precisely $l$
zeros and $l+1$ extremum points on $(-R,R)$.
 %% For $N>1$, the same
%%%properties hold for solutions $\bar W_{2l}$.
\end{proposition}

\begin{figure}
%\vskip -.3cm
 \centering
\includegraphics[scale=0.8]{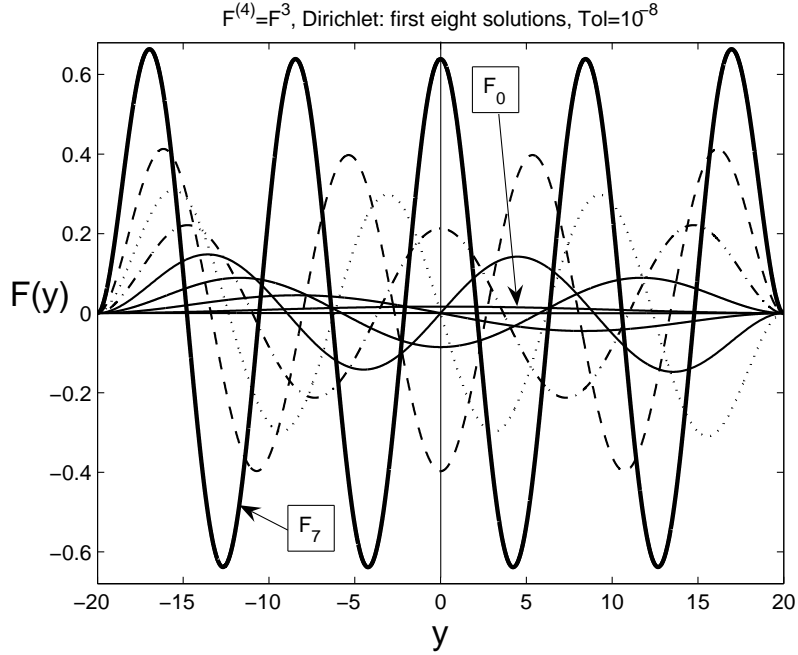}
 \vskip -.4cm
\caption{\rm\small The first eight patterns $F_{l}(y)$ satisfying
(\ref{m1}) for $m=2$.}
   \vskip -.3cm
 \label{Fm2}
\end{figure}

\begin{figure}
%\vskip -.3cm
 \centering
\includegraphics[scale=0.8]{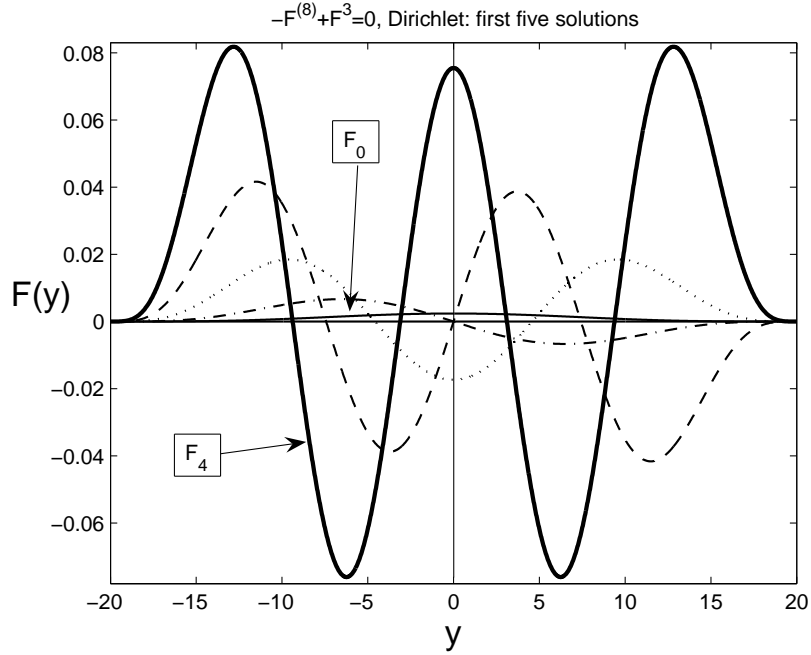}
 \vskip -.4cm
\caption{\rm\small The first five solutions $F_{l}(y)$ of
(\ref{m1}) for $m=4$.}
   \vskip -.3cm
 \label{Fm4}
\end{figure}

In Figures \ref{Fm2} and \ref{Fm4}, we present the first solutions
of (\ref{m1}) for $m=2$ and $m=4$. According to Lusternik--Schnirel'man theory
\cite[pp.~385-387]{KrasZ}, each profile $F_l$ is obtained by the
minimax variational construction (see details  in
\cite[\S~3]{GMPSob})
%%(see (\ref{ck1}) with $k \mapsto l+1$)
 on the sets of the category
$ \rho \ge l+1$.

 %%%\com{To SIP:  how to see that ? A simple justification? }

 %%%%\com{To SIP: how to see this simple? Should be true...}

%%%%%%%%%%%%%%%%%%%%%%%%%%%%%%%%%%%%%%%%%%%%%%%%%%
%%%%%%%%%%%%%%%%%%%%%%%%%%%%%%%%%%%%%%%%%%%%%%%%%%%%%%

%%%%%%%%%%%%%%%%%%%%%%%%%%%%%%%%%%%%%%%%%%%%%%%%%%
%%%%%%%%%%%%%%%%%%%%%%%%%%%%%%%%%%%%%%%%%%%%%%%%%%%%%%

\ssk

 \noi{\bf Remark: convergence to periodic solutions.} It is
clear from two Figures above (cf. the last boldface profiles)
that, for large $l \gg 1$ (actually, already for $l \ge 3$ in both
cases), the solutions of the problem $F_l(y)$ of (\ref{m1}) are
close to a periodic structure.
 Namely, denote by  $F_*(y)$ a $T_*$-periodic solution of the ODE
 (\ref{m1}) in $\re$ normalized so that
 $$
 \sup \, |F_*(y)|=1.
 $$
 Then, by scaling invariance of the equation,
  $$
  \mbox{$
  a^{-m} F_*\bigl(\frac y a\bigr) \quad \mbox{is also a solution in $\re$ for any \,
  $a>0$}.
   $}
  $$
  Therefore, for large $l \gg 1$, the following holds:
   \beq
   \label{a1rr}
    \mbox{$
   F_l(y) \approx a_l^{-m} F_*\bigl(\frac y{a_l}\bigr), \quad \mbox{where}
   \quad a_l= \frac {2R}{lT_*}\, ,
 $}
    \eeq
     and the convergence as $l \to \infty$ is uniform in $y$. Note that
     the periodic solution $F_*$ does not satisfy the Dirichlet
     boundary conditions at $y= \pm R$, and this creates some
     boundary layers. It is easy to see that these are of order $o(1)$, i.e., negligible
 as $l \to \infty$.

%%%%%%%%%%%%%%%%%%%%%%%%%%%%%%%%%%%%%%%%%%%%%%%%%%%%%%%%%%%%%%
\subsection{Homotopic connections to the cubic equation}

We now introduce the {\em basic} one-parametric family of
Dirichlet problems in $(-R,R)$ with the operators
 \beq
 \label{e1}
  \mbox{$
 {\bf A}_\e(F) \equiv (-1)^{m+1} F^{(2m)}  + (1-\e)\bigl(F-
  \bigl|\e^2+F^2\bigr|^{-\frac n {2(n+1)}}F\bigr)+ \e F^3=0,
   %%%\quad \e \in [0,1].
  $}
  \eeq
   where $\e \in [0,1]$.
For $\e=0$, we have the original problem (\ref{S2}), while $\e=1$
gives the above simpler problem (\ref{m1}) with all the solutions
ordered by Sturm's index. Notice that, for all $\e \in (0,1]$,
operators (\ref{e1}) contain analytic nonlinearities, and the
dependence on $\e$ is also analytic. The only problem of concern
is the singular limit $\e \to 0$.
 %% which makes it possible to
%%%define the following concept.

Our further construction  is naturally related to  classic theory
of homotopy of compact continuous vector fields,
\cite[\S~19]{KrasZ}. Denoting by $K(y,\xi)$ the symmetric kernel
of the linear operator $(-1)^{m} F^{(2m)}$ with zero Dirichlet
conditions, the problems (\ref{e1}) can be written in the
equivalent integral form
 \beq
 \label{e1int}
  \mbox{$
 \Phi_\e (F) \equiv F -(1-\e)\displaystyle \int K \bigl(F-
  \bigl|\e^2+F^2\bigr|^{-\frac n {2(n+1)}}F\bigr) - \e \displaystyle \int K
  F^3=0,
   $}
   \eeq
   where each integral Hammerstein operator is compact and continuous in $L^2$
   (or suitable $L^p$ spaces for $p>2$), \cite[p.~83]{KrasZ}. Therefore, the
   function (\ref{e1int}) for $\e \in [0,1]$ establishes a
   deformation of the original vector field
    $$
     \mbox{$
 \Phi_0 (F) = F -\displaystyle \int K \bigl(F-
  \bigl|F\bigr|^{-\frac n {n+1}}F\bigr)
  \quad \mbox{into} \quad  \Phi_1 (F)= F -
 \displaystyle \int K
  F^3;
   $}
 $$
 see \cite[p.~92]{KrasZ}.
If this deformation is non-singular ($0 \not \in
\s(\Phi_\e'(F))$), the two vector fields are homotopic. For
convenience, later on we consider  the differential form of
deformations bearing in mind the necessity to return to the
corresponding compact vector fields for any rigorous
justification.

Thus, we take an arbitrary solution $F(y)$ of (\ref{S2}) with
connected symmetric (always achieved by shifting) compact support
 $$
 {\rm supp} \, F= [-y_0,y_0].
  $$

%%%%%%%%%%%%%%%%%%%%%%%%%%%%%%%%%%%%%%%%%%%%%%5
\begin{definition}
 \label{D.1}
 We say that a solution $F(y)$ of $(\ref{S2})$
has Sturm index $l=I_S(F)$, if there exists its continuous
non-singular deformation, also called a homotopic connection,
  \beq
  \label{epath}
 \{F(y;\e)
\,\, \mbox{on} \,\,[-y_0,y_0]; \, \e \in [0,1]\}
 \eeq
  consisting of
critical points (solutions) of the functional for operators
$(\ref{e1})$ such that $F(y;0)=F(y)$ and $F(y;1)$ coincides with
the solution $F_l(y)$ of $(\ref{m1})$.
 \end{definition}

If, for a given solution $F$ of (\ref{S2}), such a non-singular
homotopic deformation  does not exist, then we say that Sturm's
index $I_S(F)$ cannot be attribute to such a function $F$ in
principle. In what follows, this  nonexistence result can be
associated with the fact that, for these solutions, the homotopic
connections such as (\ref{e1int}) (or  others) become singular at
some saddle-node-type (s-n) bifurcation point $\e_{\rm s-n}<1$, at
which two $\e$-branches of geometrically similar solutions meet
each other.

In general,  Sturm's index can be extended from $\e=1$ (the
ordered cubic problem (\ref{m1})) to $\e=0$ (the original one)
along {\em any} continuous analytic branch that can have arbitrary
even number of turning s-n points for $\e \in (0,1)$, and even
beyond that. Therefore, in fact we ascribe the same Sturm  index
$l$ to all profiles belonging to the same analytic branch started
for $\e=1$ at the point $F_l$ and ended up at $\e=0$. In this
sense, the nonexistence then means that such a branch in principle
is non-extensible to $\e=0$.

%% We will explain what actually means such a nonexistence
%%result.

The possibility of  bifurcation (branching) points is a key
difference between second and higher-order equations. Indeed, for
$m=1$, all the (or most  interesting) solutions have the index by
Sturm's Theorem on zero sets, while, for $m \ge 2$, there are many
others, which principally cannot obey such a  simple
classification, associated with second-order problems only.

Obviously, along any non-singular homotopic path, the critical
points are deformed continuously, which is guaranteed  by the
Inverse Function Theorem ({\em q.v.}  e.g.
  \cite[p.~319]{VainbergTr}).
Therefore,  our strategy is now to use the Lusternik--Schnirel'man/fibering method
for construction of a countable number of branches of different
profiles $\{F_l\}$, which are ordered by the category of the sets
involved. We then apply this  for any $\e \in(0,1)$ in (\ref{e1}).
Theory of compact integral operators \cite{Kras, KrasZ,
VainbergTr} then
 %%%Then the general results on stability of critical values
%%\cite[Sect. 57.6]{KrasZ} would
 suggests  existence of a
countable set of continuous $\e$-curves of critical points that
will continuously attribute  the Sturm index from the regular
problem (\ref{m1}) for $\e=1$ to the non-Lipschitz one (\ref{S2})
for $\e=0$, provided that these branches are extensible and some
of these are not destroyed at saddle-node (or others, even harder)
bifurcations in between. Lusternik--Schnirel'man and fibering  theory guarantee
existence of a countable number of extensible branches. Notice
that, ( see  \cite[p.~387]{KrasZ}), to the authors knowledge,
  $$
\mbox{``It is not known whether the Lusternik--Schnirel'man critical\\
 \,\,values
 %%(\ref{ck1})
are stable."}
 $$

 \noi On the other hand, there are some definitely stable branches.
 %% This is addressed precisely to critical
%%values of even functionals such as (\ref{f4}) constructed by the
 %%minimax S-L method.
 Therefore, in general,  we cannot guarantee that all the
Lusternik--Schnirel'man branches are extensible to $\e=1$.
 In addition, the proof of the fact that the homotopic path (\ref{e1}) (or
suitable others) is non-singular, is also a difficult open
problem.

Nevertheless, we expect that the stability or non-singularity for
(\ref{e1}) takes place for our particular problem, and we end up
this discussion as follows:

 %%%%\com{To SIP: how to get continuous $\e$-branches ?}

\ssk

%%\begin{proposition}
%% \label{Pr.Hom}
\noi {\bf Conjecture \ref{SectHom}.1.} {\em Each function $F_l$
from the basic family $\{F_l, \, l \ge 0\}$ for $(\ref{S2})$ can
be continuously deformed by $(\ref{e1})$ to the corresponding
solutions $F_l$ of $(\ref{m1})$.}
%%%% \end{proposition}

%%%\com{To SIP: any hope of proof ? Ideas ?}

 %%%%%\com{Other patterns cannot... ? Proof ? Any suggestion ???}

\ssk

 During the course of  the inverse $\e$-deformation from $\e=1$ to $\e=0$, the profiles
 $F_l(y;\e)$ get a finite, depending on $\e \in (0,1)$, number of oscillations and zeros close to
 end points $y = \pm R$, and only eventually, at $\e = 0$, this
 number gets infinity, when the nonlinearity becomes
 non-Lipschitz and the solutions become compactly supported in $(-R,R)$.

\smallskip

\begin{figure}
 %%\vskip -.3cm
\centering \subfigure[$\e$-deformation of $F_0$, $R=10$]{
\includegraphics[scale=0.7]{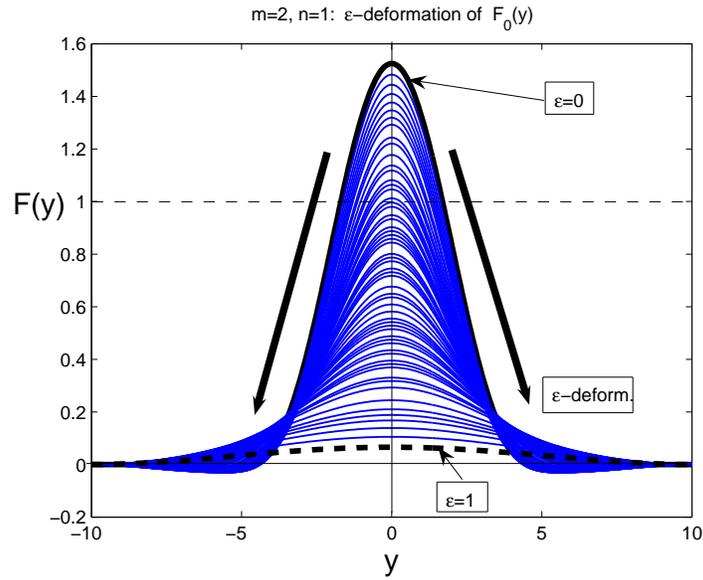} %%%%%%%%%%%%%%%%%%{GF0.eps, OLD}
} \subfigure[$\e$-deformation of $F_1$, $R=10$]{
\includegraphics[scale=0.7]{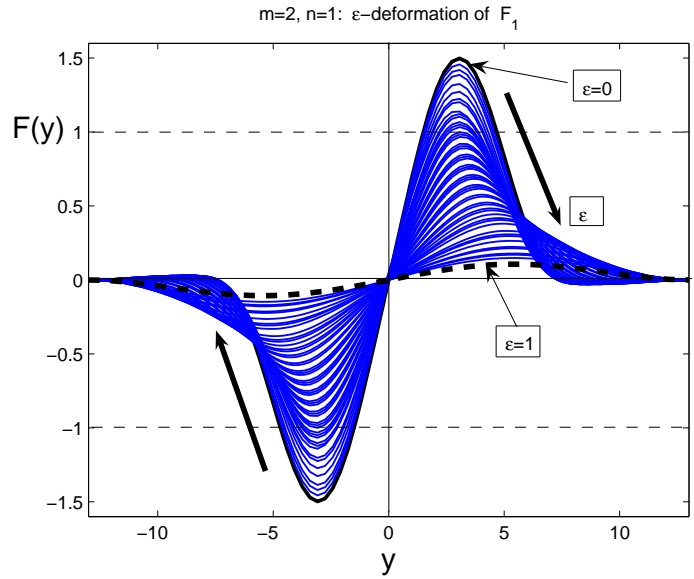}    %%%%%%%%%%%{PHCm236.eps}
}
 \vskip -.2cm
\caption{\rm\small $\e$-deformation via (\ref{e1})  of $F_0(y)$
(a) and  the 1-dipole profile  $F_1(y)$ (b).}
 \label{FFF1}
  %%%%%{FGF.1fig}
\end{figure}

\begin{figure}
 %%\vskip -.3cm
\centering \subfigure[$\e$-deformation of $F_2$, $R=16$]{
\includegraphics[scale=0.7]{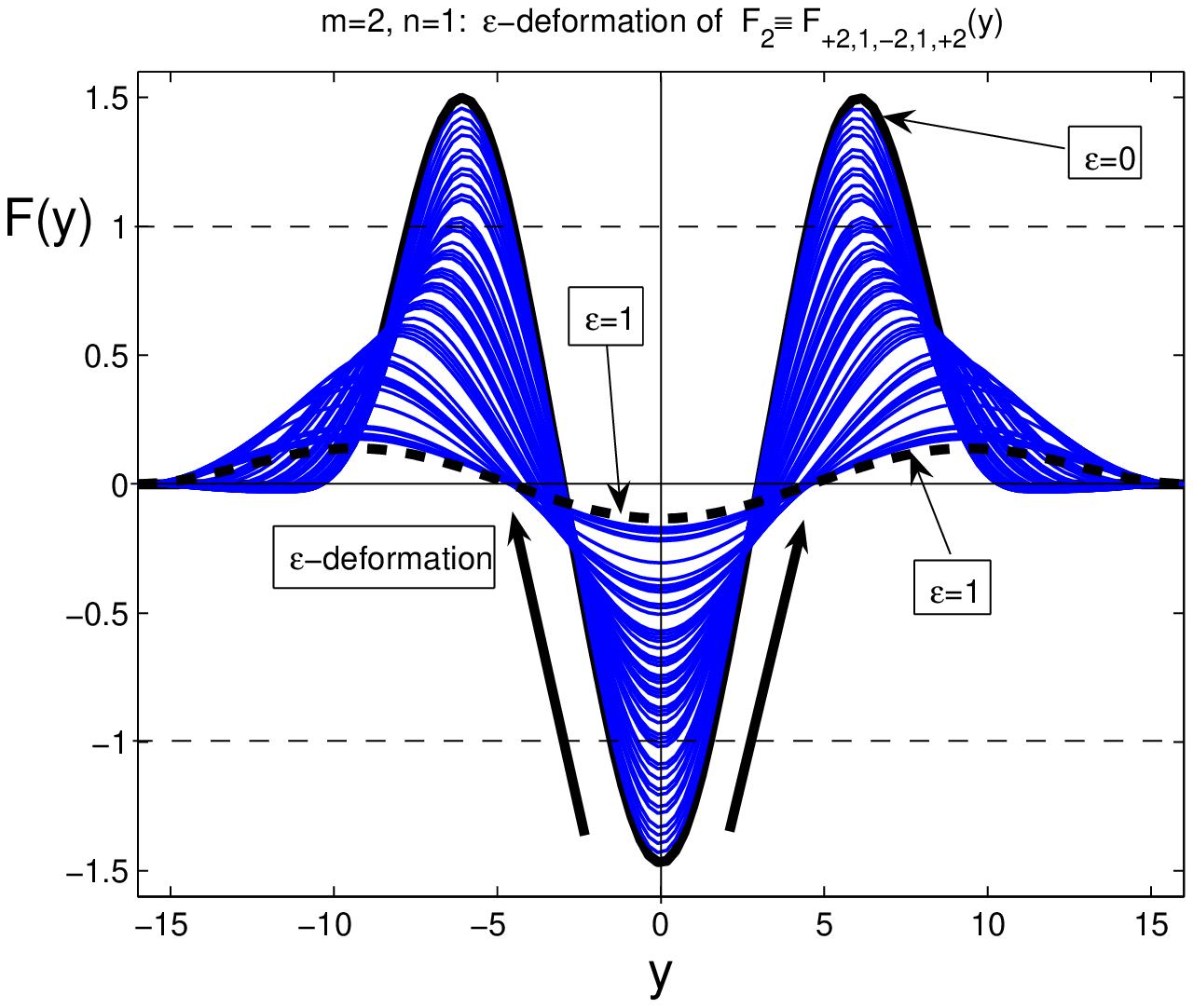}                 %%{F0eps.eps}
} \subfigure[$\e$-branch of $F_2$]{
\includegraphics[scale=0.7]{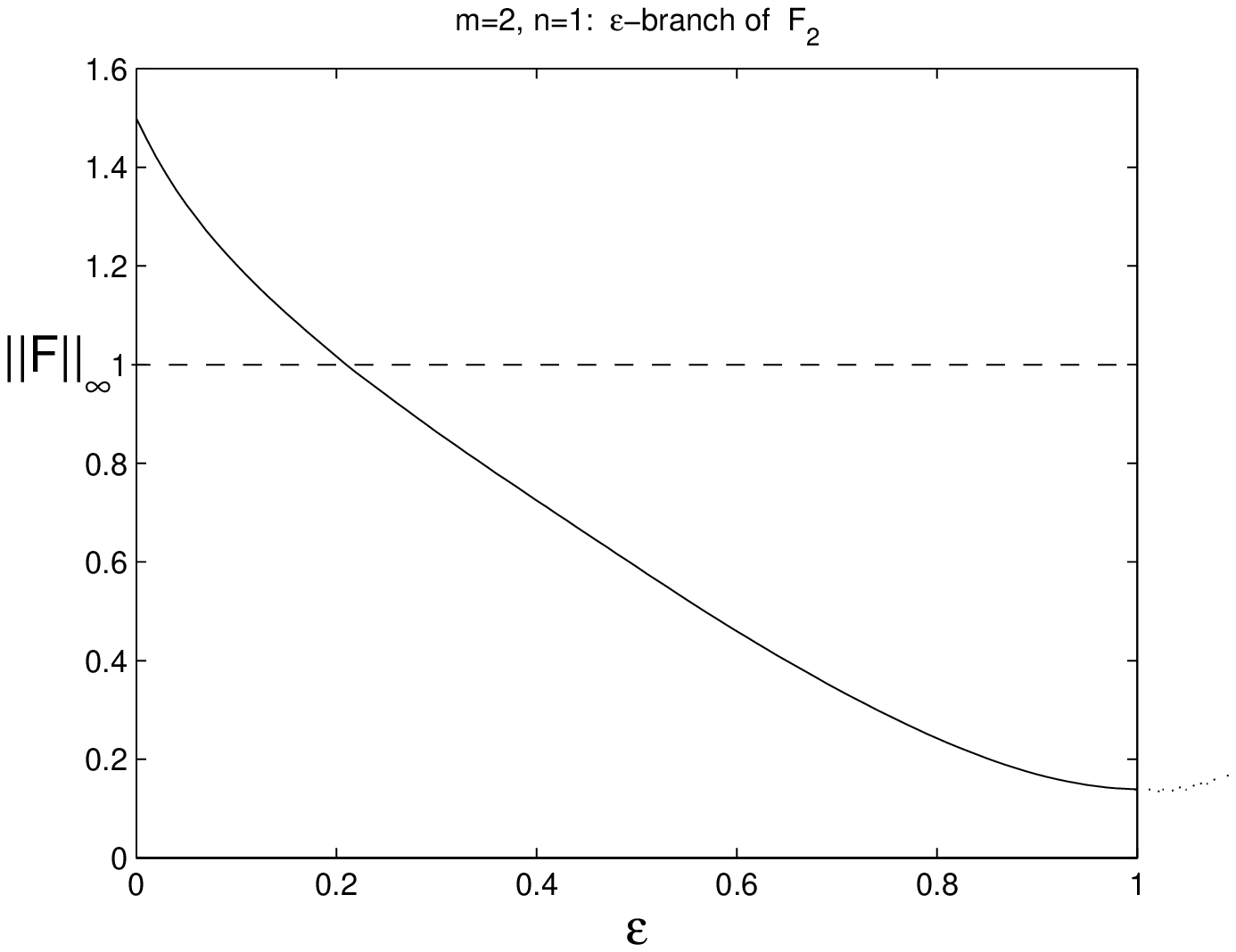}               %%%{FFP01NN.eps}
}
 \vskip -.2cm
\caption{\rm\small The $\e$-deformation of $F_2$ via homotopy
(\ref{e1}) for $n=1$ (a), and the $\e$-branch of $F_2$ (b).}
 %%% for $m=4$ and $n=1$: profiles (a), and zero
%%%%structure (b).}
 \label{FBB1}
  %%%%%{FGF.1fig}
\end{figure}

%%%%%%%%%%%%%%%%%%%%%%%%%%%%%%%%%%%%%%%%%%%%%%%%%%%%%%%%%%%%%%%%%%%%%%%
\noi\underline{\sc Existence of homotopic connections for basic
patterns $F_l$}.
As a typical example,  in Figure \ref{FFF1}, %%% and \ref{FFF2}
we present such $\e$-deformations (\ref{e1}) of  two profiles,
$F_0(y)$ (a) and the dipole $F_1(y)$ (b) for $n=1$ and $R=10$. The
$\e$-deformation of $F_2(y)$ is
 presented in  Figure
\ref{FBB1}(a) for $n=1$ and $R=10$.
 A typical corresponding
$\e$-branch of $F_2$  is shown in Figure \ref{FBB1}(b). All
$\e$-branches of the basic family $\{F_l\}$ look quite similarly.
Note that these branches can be extended beyond $\e=1$, and then
we observe there the absolute minimum of $\|F\|_\iy$ at this value
$\e =1$.

%% Here
%%we have extended these up to $\e =1.5$ and could do it further.
%%Note that at $\e=1$ (i.e., at the equation (\ref{m1})), branches
%%exhibit absolute minimum.

%%%%%%%%%%%%%%%%%%%%%%%%%%%%%%%%%%%%%%%%%%%%%%%%%%%%%%

%%%%%%%%%%%%%%%%%%%%%%%%%%%%%%%%%%%%%%%%%%%%%%%%%%

It is an obvious observation that, by continuity of branches with
respect to small changes of nonlinearities, if the homotopic
connection as in Figure \ref{FBB1}(b) takes place for the basic
deformation (\ref{e1}) and the branch is infinitely extensible for
$\e>1$, a similar connection $\e:0 \to 1$ can be achieved by other
analytic deformations. In this sense, the type of reasonable
homotopic deformations is not that crucial; cf. Proposition
\ref{Pr.NN} below establishing an analogous non-homotopy
conclusion.

\smallskip

%%%%%%%%%%%%%%%%%%%%%%%%%%%%%%%%%%%%%%%%%%%%%%%%%%%%%%%%%%%%%%%%%%%%%%

\noi\underline{\sc Nonexistence of $\e$-connections and
saddle-node bifurcations}.
%%%%%% for non-basic families}.
 We now deal with  other families of non-basic patterns
obtained in Section \ref{SFFh} by  an extra preliminary Cartesian
$h$-approximation.  We then introduce their total {\em generalized
Sturm index}, which should include the number of oscillations
about the non-trivial equilibria $\pm 1$ defined by the structure
of $h(y)$; see Section \ref{Sect8} for an alternative  approach to
the generalized index via a spatial $R$-compression of profiles.

 Then, a homotopic $\e$-deformation of
these patterns to those of the equation (\ref{m1}) with {\em
monotone} nonlinearity {\em is not possible} in principle. We
claim that, on the $\{\e,F\}$-plane of the global bifurcation
diagram, their solution branches appear in standard {\em
saddle-node bifurcations} that occur for $\e<1$, i.e., these
branches do not admit extensions up to the simpler ODEs
(\ref{m1}).

\begin{figure}
%\vskip -.3cm
 \centering
\includegraphics[scale=0.7]{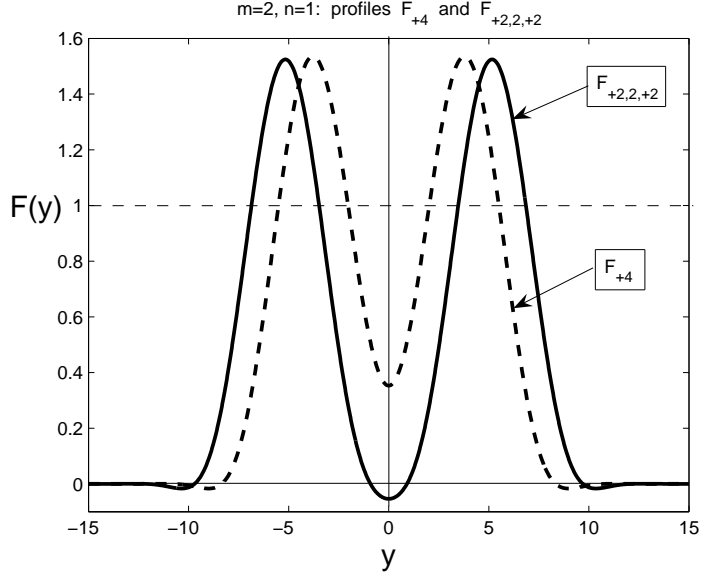}
 \vskip -.4cm
\caption{\rm\small Patterns $F_{+4}$ and $F_{+2,2,+2}$ have
similar geometric shapes and appear simultaneously at the s-n
bifurcation at (\ref{sn1}); $m=2$, $N=n=1$.}
  %%% \vskip -.3cm
 \label{FComp}
\end{figure}

 %%In Figure \ref{HH1},
  Using the enhanced numerics with
Tols=$10^{-4}$ and the step $\D \e=10^{-3}$, we show the
$\e$-deformation of two non-basic profiles given in Figure
\ref{FComp},
 \beq
 \label{bbf1}
 F_{+4}(y) \quad \mbox{and} \quad F_{+2,2,+2}(y),
 \eeq
which have  similar geometric shapes, with the equal numbers of
four intersections with the equilibrium $+1$. This  detailed
$\e$-deformation via  (\ref{e1}) for $n=1$, of $F_{+4}(y)$
%%(a)
 and
$F_2\equiv F_{+2,2,+2}(y)$ is shown in   Figure 15(a) and (b) in
\cite{GMPHomIINA} (these are too big, 5.3 and 2.33 MB, to be
presented in arXiv). The $\e$-deformations of these two profiles
 turns out to stop at the same
saddle-node bifurcation at
 $
 \e=\e_{\rm s-n}=0.709...\, .
 $

 %% see Figure \ref{FComp}
%%%given for convenience.

\begin{figure}
 \centering
 \psfrag{||F||}{$\|F\|_\infty$}
 \psfrag{e}{$\e$}
 \psfrag{t2}{$t_2$}
 \psfrag{F222}{$F_{+2,2,+2}$}
 \psfrag{F+4}{$F_{+4}$}
  \psfrag{F2}{$F_2$}
  \psfrag{final-time profile}{final-time profile}
   \psfrag{F4}{$F_4$}
\psfrag{x}{$x$}
 \psfrag{0<t1<t2<t3<t4}{$0<t_1<t_2<t_3<t_4$}
  \psfrag{0}{$0$}
 \psfrag{l}{$l$}
 \psfrag{-l}{$-l$}
\includegraphics[scale=0.5]{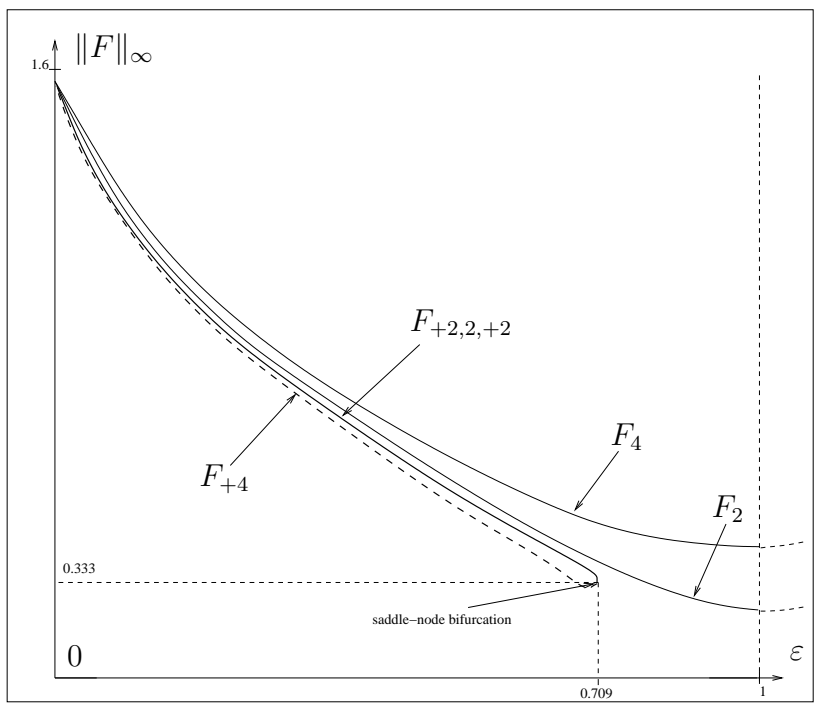}     %%%{F1EE1.pdf}
\caption{\small The $\e$-branches of  $F_{+4}$ and $F_{+2,2,+2}$
cannot be extended from $\e=0$ via (\ref{e1}) to $\e=1$ and meet
at a saddle-node bifurcation at $ \e_{\rm s-n} = 0.709...$\,. The
neighbouring branches of $F_4$ and $F_2$, which are detected in
Figure
%%\ref{HH1}(a) and (b),
 15(a), (b) in \cite{GMPHomIINA}
 are also shown.}
     \vskip -.3cm
 \label{BBF}
\end{figure}

 In Figure \ref{BBF}, we
show the corresponding
  $\e$-bifurcation diagram with a
saddle-node bifurcation at
 \beq
 \label{sn1}
  \fbox{$
 \e_{\rm s-n} = 0.709... \quad (F_{+4}, \,\,\, R=14),
  $}
   \eeq
   at which
the branch of $F_{+4}$
 %%% (see Figure \ref{G8}(a))
  and the branch of
$F_{+2,2,+2}$
 %%(see Figure \ref{G6})
  meet each other.
 For convenience, in Figure \ref{BBF}, we also draw neighbouring
global branches of the basic patterns $F_4$ and $F_2=F_{+2,-2,+2}$
(existing for all $\e \in [0,1]$), to which the  corresponding
branches jump
%%in Figure \ref{HH1} (a) and (b)
 being extended above
the bifurcation values (\ref{sn1}).
 %%%in (a) both profiles are presented.

 %%In Figure \ref{HH1}(a), we show
 It turns out that a neighbouring  branch of the
basic pattern $F_4$ exists for $\e>\e_{\rm s-n}$, while
%%in (b),
the neighbouring basic branch is that of $F_3=F_{+2,1,-2,1,+2}$.
Being extended numerically for $\e>\e_{\rm s-n}$, the
$\e$-branches of profiles (\ref{bbf1}) jump to these basic
$\e$-branches.
%%% as these Figures show.

 Such a branching at
 $\e=\e_{\rm s-n}$ means that the two profiles (\ref{bbf1}) belong to the same
 family,
both having  the generalized Sturm index $\s_{\rm min}=+4$.
%%% In
%%%Figure \ref{BBF2N}, we show the actual deformation of the
%%%profiles, which coincide at $\e=\e_{\rm s-n}$.

%%%%%%%%%%%%%%%%%%%%%%%%%%%%%%%%%%%%%%%%%%%%%%%%%%%%%%

%%\begin{figure}
 %%\vskip -.3cm
%%\centering \subfigure[$\e$-deformation of $F_{+4}$, $R=14$]{
%%\includegraphics[scale=0.52]{VV3.eps}               %%{GF4.eps}                 %%{F0eps.eps}
%%} \subfigure[$\e$-deformation of $F_{+2,2,+2}$, $R=14$]{
%%\includegraphics[scale=0.52]{UU3.eps}               %%%{FFP01NN.eps}
%%}
%% \vskip -.2cm
%%\caption{\rm\small $\e$-deformation  via  (\ref{e1}), $n=1$, of
%%$F_{+4}(y)$ (a) and $F_2\equiv F_{+2,2,+2}(y)$ stops at the same
%%saddle-node bifurcation at $\e=\e_{\rm s-n}=0.709...$\,.}
 %%% for $m=4$ and $n=1$: profiles (a), and zero
%%%%structure (b).}
%% \label{HH1}
  %%%%%{FGF.1fig}
%%\end{figure}
%%%%%%%%%%%%%%%%%%%%%%%%%%%%%%%%%%%%%%%%%%

%%%%%%%%%%%%%%%%%%%%%%%%%%%%%%%%%%%%%%%%%%%%%%%%%%
%%%%%%%%%%%%%%%%%%%%%%%%%%%%%%%%%%%%%%%%%%%%%%%%%%%%%%

It is not difficult to choose other pairs of patterns $F$ with
similar geometries, which have to be originated at saddle-node
bifurcations for $\e < 1$. For instance, these are $F_{+2k}$
%%(Figure \ref{G8}(a))
 and $F_{+2,2,+2,2,...,2,+2}$, %%(Figure \ref{G6}
with $k$ single patterns $ \sim +F_0$ gluing together.

%%%%%%%%%%%%%%%%%%%%%%%%%%%%%%%%%%%%%%%%%%%%%%%%%%%%%%%%
%%%%%%%%%%%%%

%%%%%%%%%%%%%%%%%%%%%%%%%%%%%%%%%%%%%%%%%%%%%%%
%%%\begin{figure}
 %%\vskip -.3cm
%%\centering \subfigure[of $F_{+2,2,+2}$]{
%%\includegraphics[scale=0.52]{F244.eps}
%%} \subfigure[of $F_{+4}$]{
%%\includegraphics[scale=0.52]{F244R.eps}                 %%%{FFP01NN.eps}
%%}
%% \vskip -.2cm
%%\caption{\rm\small $\e$-deformations of $F_{+2,2,+2}$ (a) and
%%$F_{+4}$ (b) via homotopy (\ref{e1}).}
 %%% for $m=4$ and $n=1$: profiles (a), and zero
%%%%structure (b).}
%% \label{FBB2N}
  %%%%%{FGF.1fig}
%%\end{figure}

\begin{figure}
 %%\vskip -.3cm
\centering \subfigure[$\e$-deformation of $F_{+6}$, $R=20$]{
\includegraphics[scale=0.7]{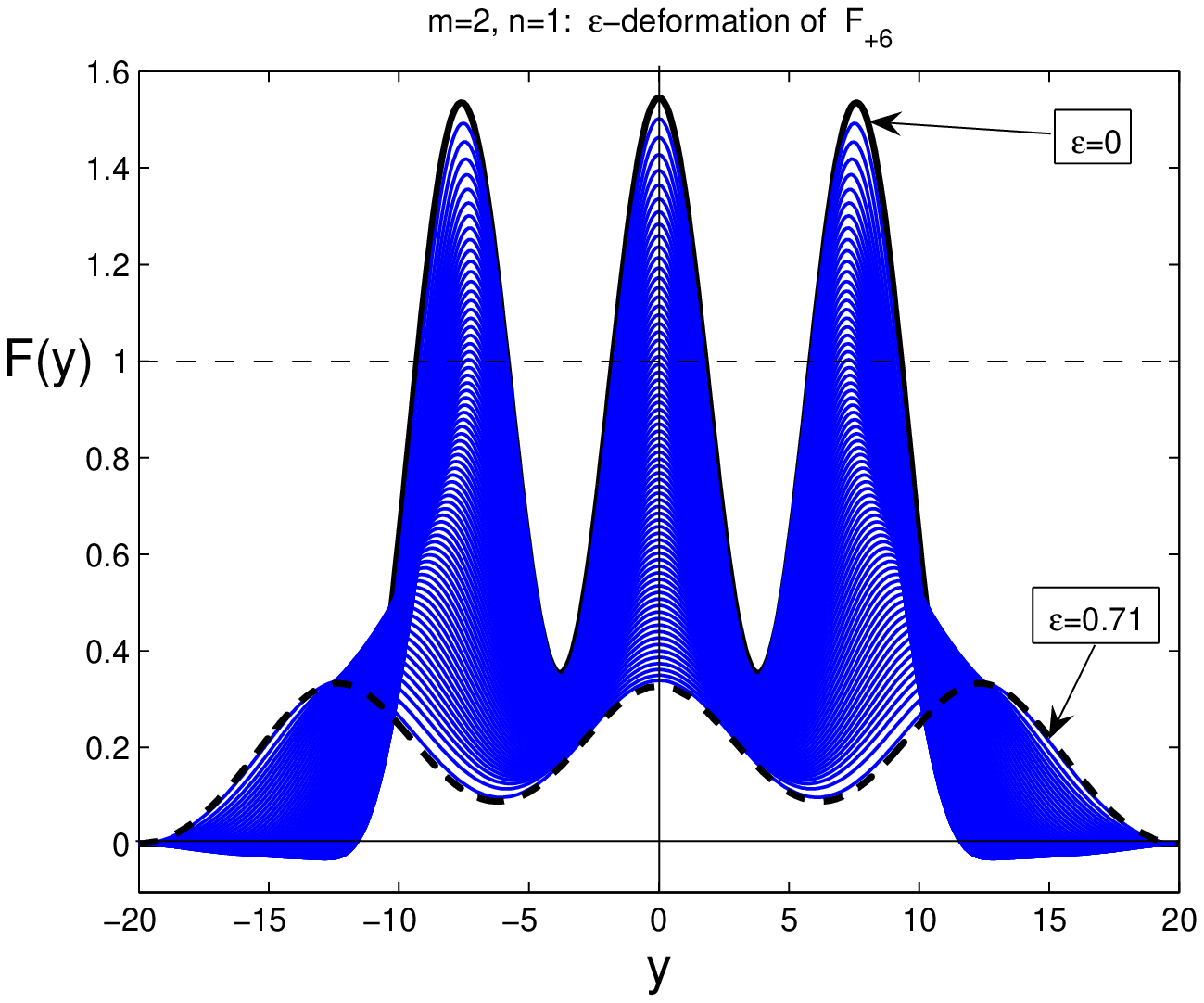}               %%{GF4.eps}                 %%{F0eps.eps}
} \subfigure[$\e$-deformation of $F_{+2,2,+2,2,+2}$, $R=20$]{
\includegraphics[scale=0.7]{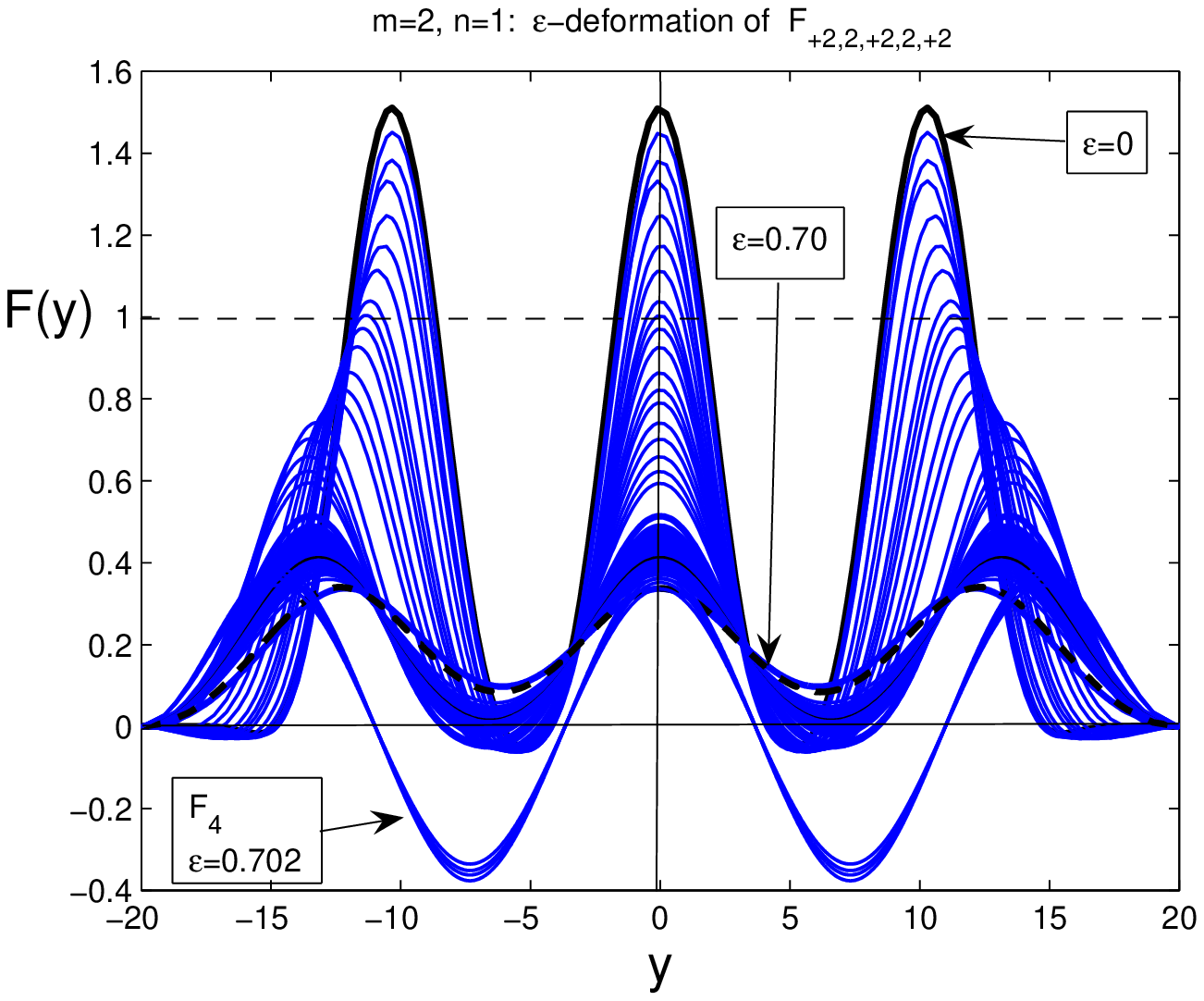}               %%%{FFP01NN.eps}
}
 \vskip -.2cm
\caption{\rm\small $\e$-deformation  via  (\ref{e1}) for $m=2$,
$n=1$, of $F_{+6}(y)$ (a) and $F_{+2,2,+2,2,+2}(y)$ stops at the
same saddle-node bifurcation at $\e=\e_{\rm s-n}=0.700...$\,.}
 %%% for $m=4$ and $n=1$: profiles (a), and zero
%%%%structure (b).}
 \label{HH1N}
  %%%%%{FGF.1fig}
\end{figure}

\begin{figure}
 %%\vskip -.3cm
\centering \subfigure[profiles $F_{+6}$ and $F_{+2,2,+2,2,+2,2}$]{
\includegraphics[scale=0.7]{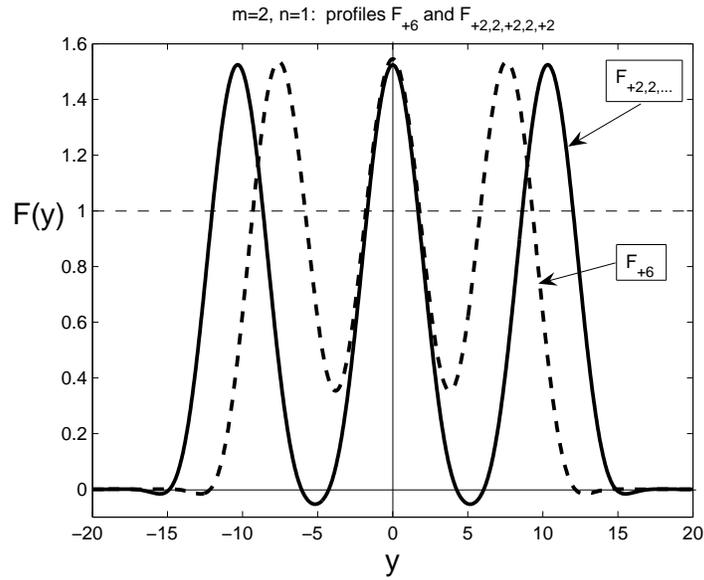}               %%{GF4.eps}                 %%{F0eps.eps}
} \subfigure[$\e$-bifurcation diagram]{
\includegraphics[scale=0.6]{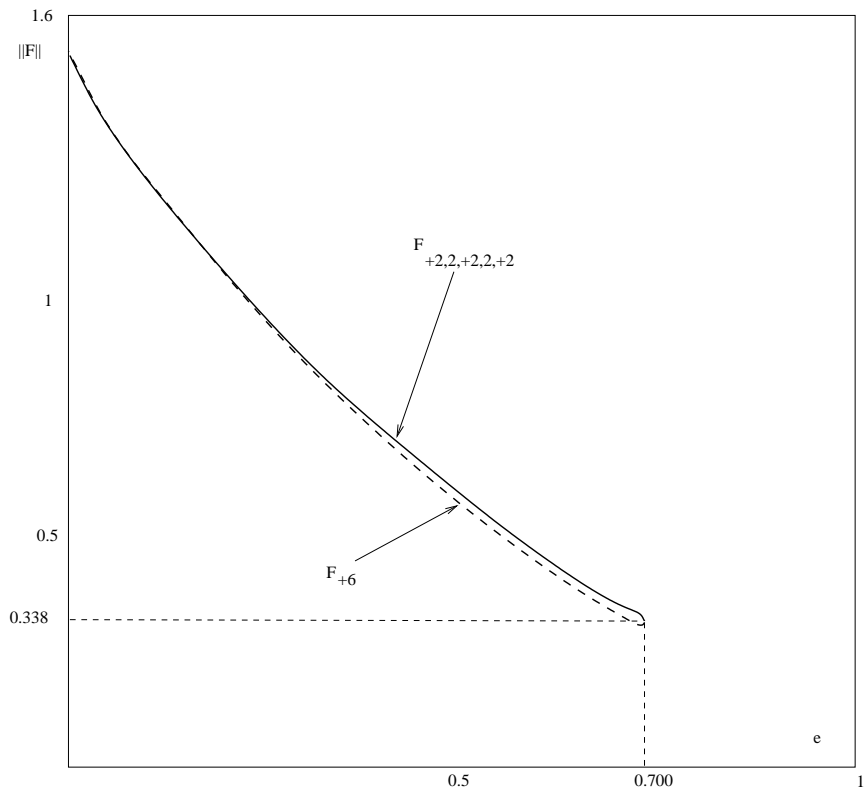}               %%%{FFP01NN.eps}
}
 \vskip -.2cm
\caption{\rm\small $m=2$, $n=1$: profiles $F_{+6}(y)$ and
$F_{+2,2,+2,2,+2}(y)$ (a); bifurcation diagram (b).}
%%%stops at the same saddle-node bifurcation at $\e=\e_{\rm
%%%s-n}=0.700...$\,.}
 %%% for $m=4$ and $n=1$: profiles (a), and zero
%%%%structure (b).}
 \label{FCompN}
  %%%%%{FGF.1fig}
\end{figure}

\begin{figure}
%\vskip -.3cm
 \centering
\includegraphics[scale=0.8]{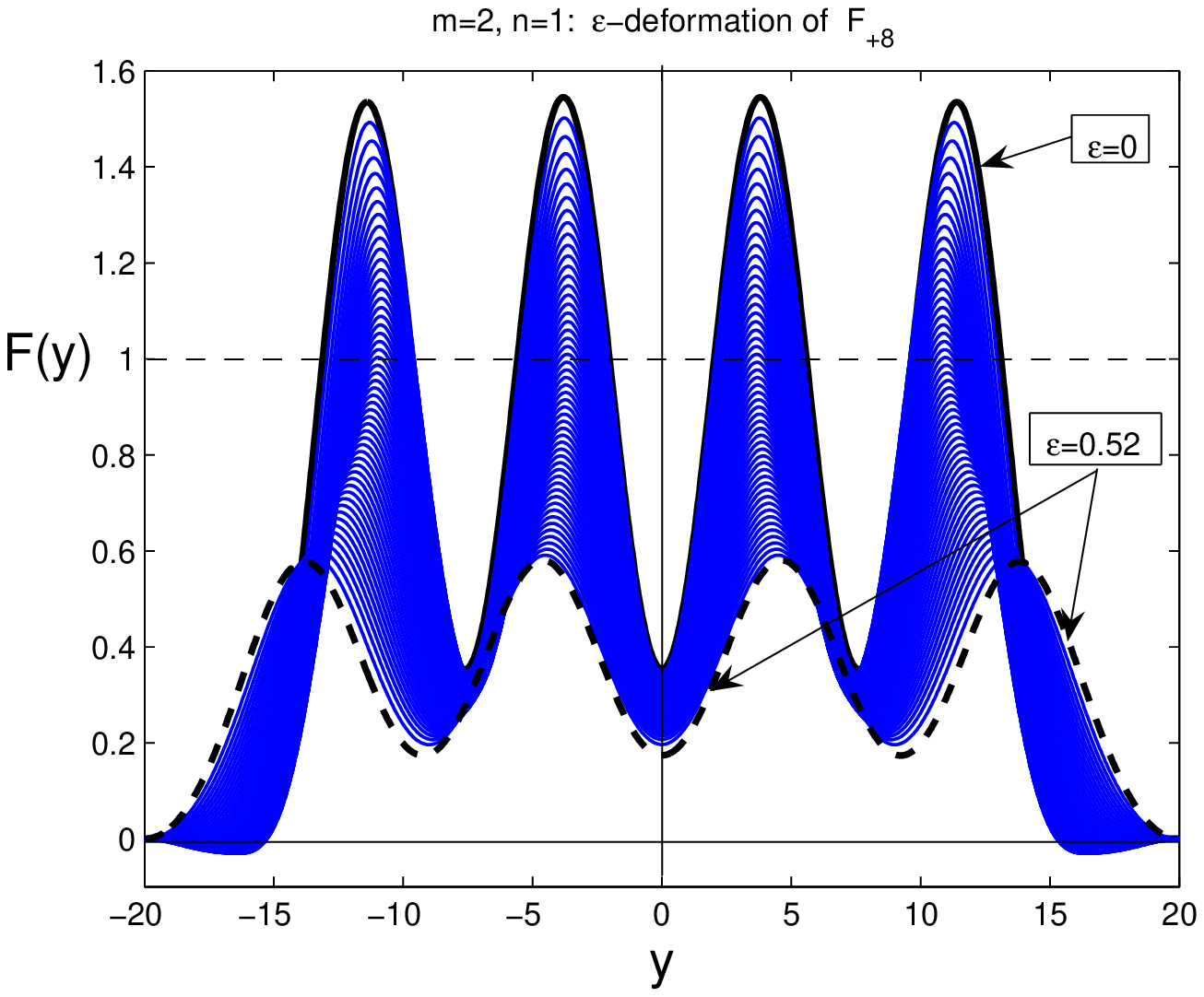}
 \vskip -.4cm
\caption{\rm\small $\e$-deformation  via  (\ref{e1}) for $m=2$,
$n=1$, of $F_{+8}(y)$ is possible until $\e_{\rm s-n}=0.52...$.}
   \vskip -.3cm
 \label{F888}
\end{figure}

For example,  Figure \ref{HH1N} shows
 the
$\e$-deformation of  (see Figure \ref{FCompN}(a))
 \beq
 \label{bbf1N}
 F_{+6}(y) \quad \mbox{and} \quad F_{+2,2,+2,2,+2}(y).
 \eeq
In Figure \ref{FCompN}(b), we show the corresponding bifurcation
diagram. Notice that the corresponding s-n bifurcation point,
 $$
  \fbox{$
 \e_{\rm s-n}=0.700... \quad (F_{+6}, \,\,\, R=20),
  $}
 $$
 is
rather close to (\ref{sn1}) for $F_{+4}$ (notice  different
lengths $R$). For $F_{+8}$ and $F_{+2,2+2,2,+2,2+2,2}$ shown in
Figure \ref{F888} for the same $R=20$, it is different,
 $$
  \fbox{$
 \e_{\rm s-n}=0.52... \quad (F_{+8}, \,\,\, R=20).
  $}
 $$

%%%%%%%%%%%%%%%%%%%%%%%%%%%%%%%%%%%%%%%%%%%%%%%%%%%%%%

%%%%%%%%%%%%%%%%%%%%%%%%%%%%%%%%%%%%%%%%%%

%%%%%%%%%%%%%%%%%%%%%%%%%%%%%%%%%%%%%%%%%%%%%%%%%%%%%%

%%%%%%%%%%%%%%%%%%%%%%%%%%%%%%%%%%%%%%%%%%

%%%%%%%%%%%%%%%%%%%%%%%%%%%%%%%%%%%%%%%%%%%%%%%%%%
%%%%%%%%%%%%%%%%%%%%%%%%%%%%%%%%%%%%%%%%%%%%%%%%%%%%%%

\ssk

Meanwhile, for convenience, we present the following simple
conclusion showing that nonexistence of the homotopic path
(\ref{e1int}) actually means nonexistence of any analytic
non-singular connections. In particular, this indicates that the
geometric type of the branching in Figure \ref{BBF} is generic.

\begin{proposition}
 \label{Pr.NN}
  Let, for a given solution $F$ of $(\ref{S2})$, the basic deformation
   $(\ref{e1int})$ have a singular point at some $\e_{\rm s-n}<1$,
   where two continuous branches of two patterns
    originated at $\e=0$ meet each other and hence cannot be continued up
    to $\e=1$.
   Then any analytic deformation of these patterns generating the functional path
$(\ref{epath})$ ends up at a singular point for some $\e \in
(0,1)$.
 \end{proposition}

In other words, other analytic deformations cannot move the s-n
point into the set $\{\e>1\}$ just by continuity.

\ssk

\noi{\em Proof.} Without loss of generality, we assume that the
corresponding critical values and the points  $\{F_l\}$ of the
cubic problem (\ref{m1}) are non-singular (by changing $R$ if
necessary). Since any continuous deformation of the basic path
will continuously (and analytically) deform the branches, the
existence of a homotopic path would mean that at some instant, the
s-n point of the branches will touch the vertical line $\e=1$. At
this moment, we would create a singular value for the analytic
cubic problem (\ref{m1}), a contradiction. $\qed$

\ssk

Numerically, we have observed a curious phenomenon: in a left-hand
neighbourhood of the saddle-node bifurcation at $\e=\e_{\rm s-n}$,
the profiles keep only essential non-monotonicity features and
loose  all intersections with zero, so become non-oscillatory near
transversal zeros. Therefore, according to such an
$\e$-deformation to saddle-node bifurcations,
 the number of intersections with
the trivial equilibrium 0 at $\e=0$ should  not be taken into
account in the generalized Strum index.
 In this sense, the complicated profiles on Figure \ref{FCmm2} have the
 following generalized Sturm index:
  $$
  \s_{\rm min}=\{-8,+4,-10,+8-12\};
  $$
see Section \ref{Sect8} for details and more mathematics. One can
``split" this index to get equivalent pairs of profiles originated
at some $\e_{\rm s-n}<1$.

\ssk

%%%\com{ALL:  why at the s-n bifurcation, all non-essential zeros
 %%%%disappear?}

%%%%%%%%%%%%%%%%%%%%%%%%%%%%%%%%%%%%%%%%%%%%%%%%%%
 \subsection{Homotopic connection to linear eigenvalue problems}

This is an alternative way to ascribe Sturm's index to basic
patterns $\{F_l$\}.
 Modifying the nonlinearity in the approximation
 (\ref{e1}), we consider the following operator family
 (watch the last term):
\beq
 \label{e1N}
  \mbox{$
 \hat{\bf A}_\e(F)= (-1)^{m+1} F^{(2m)}  + (1-\e)\bigl(F-
  \bigl|\e^2+F^2\bigr|^{-\frac n {2(n+1)}}|F|^\e F\bigr), \quad \e \in [0,1].
  $}
  \eeq
  As usual, by the actual homotopic connection we mean the
  corresponding vector fields composed of compact integral
  operators.
Then, from (\ref{e1N}) at $\e=1$, we obtain the linear operator
 $$
 \hat{\bf A}_1(F)= (-1)^{m+1} F^{(2m)}.
  $$

On the other hand, for any $\e>0$, the linearized operator at 0 is
very simple,
 \beq
 \label{e2N}
 \hat{\bf A}_1'(0)= (-1)^{m+1} D^{2m}_y +(1-\e)I, \quad  \e \in (0,1].
 \eeq
 Denoting by $\s=\{-\l_l >0, \, l \ge 0\}$ the eigenvalues of the negative operator
  $(-1)^{m}
 D^{2m}_y$,
it follows that \cite[Ch.~8]{KrasZ}
 $$
  \e_l=1-\l_l \quad \mbox{for any}  \quad l \ge 0
  $$
  are subcritical bifurcation points, where the necessary $\e$-branches
  appear. Along those branches that are originated at $\e_l>0$,
  for
  $\e=0$ we obtain our nonlinear eigenfunctions $\{F_l\}$, which
  thus inherit the Sturmian structure and the index from the eigenfunction $\psi_l$
  (it has precisely $l$ zeros and $l+1$ extrema points; see Section \ref{SectAn}) that governs
  the pattern for $\e \approx\e_l^-$.
  We continue developing such an $\e$-deformation approach to linear eigenvalue problems
   in Section \ref{SectAn}
devoted to ODEs with analytic nonlinearities.

%%%%%%%%%%%%%%%%%%%%%%%%%%%%%%%%%%%%%%%%%%%%%%%%%%%
\section{\underline{\bf Problem with ``fast diffusion"} : extinction and blow-up  phenomenon in
the  Dirichlet setting}
 %%parabolic and hyperbolic
%%PDEs}
 \label{S.Ext}

Here, using  typical concavity-like techniques,
%% of nonlinear capacity and eigenfunction
%%type,
we prove that {\em finite-time extinction} for the PDE (\ref{S1F})
%% and
%%(\ref{S3F})
is a generic property of their bounded weak solutions. Firstly, we
study this phenomenon in a bounded domain in $\ren$ with
homogeneous Dirichlet boundary conditions. Secondly, for the
Cauchy problem, possible (and rather complicated) types of
extinction patterns will be revealed in the next Section
\ref{SectAn} by using our separate variable similarity patterns.

%%%\com{SIP and EM: please have a look again}

%%\begin{center}
%%  \Large
%%  \bf

 \subsection{Extinction for some nonlinear parabolic problems of higher order: main result}
%%\end{center}

%%\section{Settings of problem}

Let $\O$ be a bounded sufficiently smooth domain in $\ren$. Taking
the original equation (\ref{S1F}) and setting, as usual,
  $$
   \mbox{$
  |u|^nu=v \quad \Longrightarrow \quad u=|v|^{- \frac n{n+1}}\, v,
  \quad \mbox{where} \quad - \frac n{n+1}>0 \,\,\,\mbox{for}
  \,\,\, n \in (-1,0),
  $}
  $$
 we arrive at the following initial boundary value problem:
 \beq
  \label{111}
   \left\{
\begin{aligned}
&\mbox{$
    \frac \partial{\partial t} \psi(v)$} = (-1)^{m+1} \Delta^m v +
    v && \mbox{ in } Q =  \Omega \times \re_+, \\
&    v
=
%%%\mbox{$
 Dv = \ldots = D^{m-1} v = 0 && \mbox{ on }
\partial\Omega \times \re_+,  \\
 &   v(x,0) = v_0(x) && \mbox{ in } \Omega,
\end{aligned}
\right.
 \eeq
 where $v_0$ is an initial function from an appropriate space to
 be specified.
 Here, the only nonlinearity is $\psi(v) = |v|^{-\frac n{n+1}} v$.
%%%% with
%%%$-1<n<0$.
%% and $\Omega$ is a bounded domain in $\R^N$.
We examine the problem (\ref{111}) in the ``native'' energy Sobolev
space.

%%\section{Main result}

Let us introduce the following functionals associated with the
operators in (\ref{111}):
\begin{align*}
 \mbox{$
    \Phi(t)$} &:= \mbox{$\frac 12 \displaystyle \int_\Omega |v(t,x)|^{\frac {n+2}{n+1}} \,{\mathrm d}x,
    $} \ssk
     \\
     \mbox{$
    E(t)$} &:= \mbox{$-\displaystyle \int_\Omega |\tilde D^m v(t,x)|^2 \,{\mathrm d}x - \displaystyle \int_\Omega v^2 \,{\mathrm d}x.
    $}
\end{align*}

\begin{lemma} There holds
\[
 \mbox{$
    \Phi'(t) = \frac {n+2}2 E(t).
    $}
\]
\end{lemma}

%%\begin{proof}

The proof follows from simple calculations. The main result on
extinction in (\ref{111}) is as follows:

%%\end{proof}

\begin{theorem}{(\bf Extinction)}
 \label{Th.Ext}
 For given nontrivial initial data, denote:
\begin{align*}
 \mbox{$
    \Phi(0)$} & \mbox{$= \frac 12 \displaystyle \int_\Omega |v_0|^{\frac {n+2}{n+1}} \,{\mathrm d}x > 0,
    $}\ssk
     \\
    E(0) &=
    \mbox{$-\displaystyle \int_\Omega |\tilde D^m v_0|^2 \,{\mathrm d}x - \displaystyle \int_\Omega v_0^2 \,{\mathrm d}x < 0.
    $}
\end{align*}
Then
 \beq
 \label{Ext66}
 \mbox{$
    \Phi(t) \le \Phi(0) \frac 1{\left(1 - \frac tT \right)^{\frac {n+2}n}} = \left(1 - \frac tT \right)^{\left| \frac {n+2}n \right|} \Phi(0) \to 0 \,\,\, \mbox{as} \,\,\, t \to T^-,
    $}
 \eeq
with
\[
 \mbox{$
    T := \frac {n+2}n \frac {\Phi(0)}{\Phi'(0)} = \frac 2n \frac {\Phi(0)}{E(0)} > 0.
    $}
\]
\end{theorem}

%%\section{The proof}

The proof is divided into several steps.

\subsection{The first energy relation}

Multiplying the equation in (\ref{111})
%%\[
%%    (\psi(v))_t = (-1)^{m+1} \Delta^m v + v
%%%%\]
by $v$ and integrating by parts over $\O \times (0,T)$  by taking
into account the boundary conditions, we obtain
\begin{equation}\label{e:e1}
 \mbox{$
    \frac 1{n+2} \displaystyle \int_\Omega \left. |v|^{\frac {n+2}{n+1}} \right|_0^T \,{\mathrm d}x =
     -\iint_\Omega |\tilde D^m v|^2\, \,{\mathrm d}x\, {\mathrm d}t -
     \iint _\Omega v^2 \, \,{\mathrm d}x\, {\mathrm d}t.
    $}
\end{equation}

\subsection{The second energy relation}

Multiplying the equation in (\ref{111}) by $v_t$ and again
integrating  by parts over $\O \times (0,T)$ by using the boundary
conditions, we obtain
\begin{equation}\label{e:e2}
 \mbox{$
    \displaystyle \iint_\Omega \psi'(v)v_t^2\,{\mathrm d}x\, {\mathrm d}t = \left.
    \left( -\frac 12 \displaystyle \int_\Omega |\tilde D^m v|^2 \,{\mathrm d}x - \frac 12 \displaystyle \int_\Omega v^2 \,{\mathrm d}x \right) \right|_0^T.
$}
\end{equation}

\subsection{Connection (a main point)}

Denote
\[
\mbox{$
    \tilde U(t) := -\displaystyle \int_\Omega |\tilde D^m v|^2 \,{\mathrm d}x - \displaystyle \int_\Omega |v|^2 \,{\mathrm d}x,
$} \quad
%%\]
%%\[
 \mbox{$
    \tilde V(t) := \frac 1{n+2} \displaystyle \int_\Omega |v|^{\frac {n+2}{n+1}} \,{\mathrm d}x.
$}
\]
Then the identity \eqref{e:e1} reads
 \beq
\label{e:1*} %%%%\tag{$1^*$}
 \mbox{$
    \left. \tilde V(t) \right|_0^T = \displaystyle \int_0^T \tilde U(t) {\mathrm d}t,
$}
%%\]
\quad \mbox{or equivalently}, \quad
%%%%\begin{equation}\label{e:1*} \tag{$1^*$}
 \mbox{$\displaystyle
    \frac {d\tilde V}{{\mathrm d}t} = \tilde U(t).
$}
\end{equation}
Analogously, it follows from \eqref{e:e2} that
\begin{equation}\label{e:2*} %%\tag{$2^*$}
 \mbox{$\displaystyle
    \frac {d\tilde U}{{\mathrm d}t} =
     \frac 2{n+1} \displaystyle \int_\Omega |v|^{-\frac n{n+1}} (v_t)^2\, \,{\mathrm d}x.
$}
\end{equation}

\subsection{The main (crucial) relation}

We deduce from \eqref{e:1*} and \eqref{e:2*} that
\begin{equation}\label{e:3}
 \mbox{$
   \displaystyle  \frac {d^2 \tilde V}{{\mathrm d}t^2} = \frac {d \tilde U}{{\mathrm d}t} =
    \frac 2{n+1} \displaystyle \int_\Omega |v|^{-\frac n{n+1}} (v_t)^2 \,{\mathrm d}x.
 $}
\end{equation}
Next, by the definition of $\tilde V$, we have
\[
 \mbox{$
    \tilde V_{tt} =\displaystyle \frac 1{(n+1)^2} \displaystyle \int_\Omega |V|^{-\frac n{n+1}} (v_t)^2
     \,{\mathrm d}x + \frac 1{n+1} \displaystyle \int_\Omega |v|^{-\frac n{n+1}} v v_{tt} \,{\mathrm d}x.
 $}
\]
Thus, equation \eqref{e:3} takes the form
\begin{equation}\label{e:4}
\mbox{$
    \displaystyle \int_\Omega |v|^{-\frac n{n+1}} v v_{tt}\, {\mathrm d}x =
     \frac {2n+1}{n+1} \displaystyle \int_\Omega |v|^{-\frac n{n+1}} (v_t)^2\, \,{\mathrm d}x.
 $}
\end{equation}

\subsection{Replacement}

We introduce a new function $z=z(x,t)$ by the formula
\[
   \displaystyle v = |z|^{\frac n{n+2}} z.
\]
Then equation \eqref{e:4} reads
\begin{equation}\label{e:5}
 \mbox{$
    \displaystyle \int_\Omega zz_{tt}\, \,{\mathrm d}x = c_n \displaystyle \int_\Omega (z_t)^2\, \,{\mathrm d}x,
    \quad \mbox{with \,\, $\displaystyle c_n = \frac {3n+2}{n+2}$}.
$}
\end{equation}
%%%%with $c_n = \frac {3n+2}{n+2}$.

\subsection{Fourier analysis of equation \eqref{e:5}}

Let $(e_k)_{k\in\mathbf N}$ be a complete orthonormal system in $L^2(\Omega)$. Then we have
\begin{align*}
   z(x,t) &= \mbox{ $\sum_{k=1}^\infty z_k(t) e_k(x),$} \\
    z_t(x,t) &= \mbox{ $\sum_{k=1}^\infty z_k'(t) e_k(x),$} \\
    z_{tt}(x,t) &= \mbox{ $ \sum_{k=1}^\infty z_k''(t) e_k(x).$}
\end{align*}
Let
\[
 \mbox{$\displaystyle
    \Phi(t) = \frac 12 \displaystyle \int_\Omega z^2(x,t) \,{\mathrm d}x = \frac 12 \sum_{k=1}^\infty z_k(t)^2,
$}
\]
 so that
\begin{equation}\label{e:6}
 \mbox{$\displaystyle
    \Phi'(t) = \frac 12 \frac {\mathrm d}{{\mathrm d}t} \displaystyle \int_\Omega |z|^2\, \,{\mathrm d}x = \displaystyle \int_\Omega z z_t\, \,{\mathrm d}x = \sum_{k=1}^\infty z_k z_k',
 $}
\end{equation}
and hence, by Holder's inequality,
\[
 \mbox{$\displaystyle
    (\Phi')^2 \le \left( \sum_{k=1}^\infty z_k^2 \right) \left( \sum_{k=1}^\infty z_k'^2 \right)
    = 2\Phi \sum_{k=1}^\infty z_k'^2.
 $}
\]
In Fourier coefficients, equation \eqref{e:5} takes the form
\begin{equation}\label{e:e}
 \mbox{$\displaystyle
    \sum_{k=1}^\infty z z_k'' = c_n \sum_{k=1}^\infty z_k'^2.
 $}
\end{equation}
On the other hand, in terms of the function $\Phi$, we have
\[
 \mbox{$\displaystyle
    \Phi'' = \displaystyle\sum_{k=1}^\infty z_k z_k'' + \sum_{k=1}^\infty z_k'^2.
$}
\]
Consequently, equation \eqref{e:e} takes the form
\[
 \mbox{$\displaystyle
    \Phi'' - \displaystyle\sum_{k=1}^\infty z_k'^2 = c_n \sum_{k=1}^\infty z_k'^2.
$}
\]
Therefore,
\[
 \mbox{$\displaystyle
    \Phi'' = (1+c_n) \sum_{k=1}^\infty z_k'^2 \ge \frac {1+c_n}{2}\, \displaystyle \frac  {(\Phi')^2}{\Phi},
$} \quad \mbox{i.e.,}
\]
%%i.e.,
\begin{equation}\label{e:ee}
   \displaystyle  \Phi\Phi'' \ge k_n \Phi'^2, \quad \mbox{where}
\end{equation}
%%where
\[
 \mbox{$\displaystyle
    k_n := \frac {1+c_n}2 = \frac {2(n+1)}{n+2} < 1 \quad (\mbox{for } n \in (-2,0)).
$}
\]

\subsection{Extinction: the proof}

For the analysis of an ordinary differential inequality appeared, we can use various approaches. For the case where $\Phi'<0$, we apply the standard approach.
Namely, we divide inequality \eqref{e:ee} by $\Phi'$. Then we obtain
\[
 \mbox{$\displaystyle
    \frac {\Phi''}{\Phi'} \le k_n \frac {\Phi'}\Phi.
$}
\]
From here, it then follows that
\[\displaystyle
    \Phi'(t) \le C_1 \Phi^{k_n}(t),
\]
with $\displaystyle C_1 = \frac {\Phi'(0)}{\Phi^{k_n}(0)} < 0$ for $E(0)<0$.

This inequality implies the result of the main theorem.
For the proof, it suffices  to replace our new variables with the original function $v$ by the formula
\[
 \mbox{$\displaystyle
    \Phi(t) = \frac 12 \displaystyle \int_\Omega z^2(x,t) \,{\mathrm d}x = \frac 12 \displaystyle \int_\Omega |v|^{\frac {n+2}{n+1}}\, \,{\mathrm d}x
$}
\]
so that
\begin{align*}
    \mbox{$\Phi'(t) = \displaystyle\frac {d\Phi}{{\mathrm d}t} = \displaystyle \int_\Omega z z_t\, \,{\mathrm d}x $}
     &= \mbox{$\displaystyle \frac 12 \, \frac {n+2}{n+1} \displaystyle \int_\Omega \psi(v) v_t\, \,{\mathrm d}x $} \\
    &= \mbox{$\displaystyle \frac 12 \, \frac {n+2}{n+1} \displaystyle \int_\Omega \frac {\psi(v)}{\psi'(v)} \psi'(v) v_t\, \,{\mathrm d}x $}\\
    &= \mbox{$\displaystyle \frac 12 \, \frac {n+2}{n+1} \displaystyle \int_\Omega (n+1) v \psi(v)_t\, \,{\mathrm d}x $} \\
    &= \mbox{$\displaystyle \frac {n+2}2 \displaystyle \int_\Omega v \psi(v)_t\, \,{\mathrm d}x $}\\
    &= \mbox{$\displaystyle \frac {n+2}2 E(v),$}
\end{align*}
where
\[
 \mbox{$\displaystyle
    E(v) = -\displaystyle \int_\Omega |\tilde D^m v|^2\, \,{\mathrm d}x - \displaystyle \int_\Omega v^2\, \,{\mathrm d}x. $}
\]
This completes the proof of Theorem \ref{Th.Ext}.\\

%%%%%%%%%%%%%%%%%%%%%%%%%%%%%%%%%
\subsection{A diversion to blow-up for $n>0$}

By using the same computations as above we can prove the following blow-up result for the original equation (\ref{S1}).

\begin{theorem}{(\bf Blow-up)}
\label{Th.Blow-up}
Suppose that  $n>0$ in $(\ref{111})$.  Let
\begin{align*}
    \Phi(0) &= \frac 12 \int_\Omega |v_0|^{\frac {n+2}{n+1}} \,{\mathrm d}x > 0, \\
    E(0) &= -\int_\Omega |\tilde D^m v_0|^2 \,{\mathrm d}x + \int_\Omega v_0^2 \,{\mathrm d}x > 0.
\end{align*}
Then
\[
    \Phi(t) \ge \Phi(0) \frac 1{\left(1 - t/T \right)^{\frac {n+2}n}}
\]
with
\[
    T := \frac {n+2}n \frac {\Phi(0)}{\Phi'(0)} = \frac 2n \frac {\Phi(0)}{E(0)}.
\]
\end{theorem}

\subsection{Blow-up: the proof} In this case, following the computations as in  \eqref{e:e}, we obtain,
\[
 \mbox{$\displaystyle
    \Phi'' - \displaystyle\sum_{k=1}^\infty z_k'^2 = c_n \sum_{k=1}^\infty z_k'^2.
$}
\]
Therefore,
\[
    \Phi'' = (1+c_n) \sum_{k=1}^\infty z_k'^2 \ge \frac {1+c_n}{2\Phi} \Phi'^2,
\]
i.e.,
\begin{equation}\label{e:eeN}
    \Phi\Phi'' \ge k_n \Phi'^2,
\end{equation}
where
\[
    k_n := \frac {1+c_n}2 = 2 \frac {n+1}{n+2} > 1 \quad (\mbox{for } n>0).
\]
For the case where $\Phi'>0$, we apply the standard approach.

Namely, we divide inequality \eqref{e:eeN} by $\Phi'$. Then we obtain
\[
    \frac {\Phi''}{\Phi'} \ge k_n \frac {\Phi'}\Phi.
\]
From here, it follows that
\[
    \Phi'(t) \ge C_1 \Phi(t)^{k_n}
\]
with $C_1 = \dfrac {\Phi'(0)}{\Phi(0)^{k_n}} > 0$ for $E(0)>0$.

This inequality implies the result of the main theorem.

Indeed, it is enough to replace our new variables with the original function $v$ by the formula
\[
    \Phi(t) = \frac 12 \int_\Omega z(t,x)^2\, \,{\mathrm d}x = \frac 12 \int_\Omega |v|^{\frac {n+2}{n+1}}\, \,{\mathrm d}x
\]
and
\begin{align*}
    \Phi'(t) = \frac {{\mathrm d}\Phi}{{\mathrm d}t} = \int_\Omega z z_t\, \,{\mathrm d}x &= \frac 12 \cdot \frac {n+2}{n+1} \int_\Omega \psi(v) v_t\, \,{\mathrm d}x \\
    &= \frac 12 \cdot \frac {n+2}{n+1} \int_\Omega \frac {\psi(v)}{\psi'(v)} \psi'(v) v_t\, \,{\mathrm d}x \\
    &= \frac 12 \cdot \frac {n+2}{n+1} \int_\Omega (n+1) v \psi(v)_t\, \,{\mathrm d}x \\
    &= \frac {n+2}2 \int_\Omega v \psi(v)_t\, \,{\mathrm d}x \\
    &= \frac {n+2}2 E(v),
\end{align*}
where
\[
    E(v) = -\int_\Omega \left( |\tilde D^m v|^2 + v^2 \right) \,{\mathrm d}x.
\]

%%%%%%%%%%%%%%%%%%%%%%%%%%%%%%%%%%%%%%%%%%%%%%%%%%%%%
%%%%%%%%%%%%%%%%%%%%%%%%%%%%%%%%%%%%%%%%%%%%%%%%%%%%%%%
\section{\underline{\bf Problem ``fast diffusion":} existence,
multiplicity, etc.}
%% on the ODE
%%%$F^{(4)}=-F+F^3$}
%% from fast diffusion-dispersion PDEs}
 \label{SectAn}

%%%%%%%%%%%%%%%%%%%%%%%%%%%%%%%%%%%%%%%%%%%%%%%%%%%%%%%%
\subsection{Oscillatory ODEs with analytic nonlinearities from fast diffusion}

Here we consider another ODE model (\ref{nn1}), and without loss
of generality, we mainly restrict to $m=2$.
 %% \beq
 %%\label{nn1}
 %%F^{(4)} = -F+ F^3 \quad \mbox{in} \quad \re,
%%\eeq
 We show that (\ref{nn1}) provides us with similar countable families of various
 patterns. Moreover, we claim that the solution set of (\ref{nn1})
 is equivalent to that obtained earlier for non-Lipschitz
 nonlinearities.
 %%Equation (\ref{nn1}) is a non-coercive version of  the
 %%stationary ODE occurred in the extended
 %%Kolmogorov--Petrovskii--Piskunov--Fisher model,
 %%$$
 %%u_t= - u_{xxxx} + \b u_{xx} + u - u^3
  %% \,\, \Longrightarrow  \,\, F^{(4)} = \b F'' + F - F^3.
 %%$$
 %%which leads to (\ref{KPP1});
  %%%There exists a large amount of literature devoted to such ODEs
 %%%with $\b > 0$;
  %%%  see references and results  in \cite{PelTroy, KKVV00, VV02}.
%%%Consider  equation (\ref{nn1}).
  Indeed,  solutions $F(y)$ of (\ref{nn1}) are not compactly
 supported and exhibit {\em oscillatory} exponential decay
at infinity governed by the linearized operator
 \beq
 \label{nn2}
  \mbox{$
F^{(4)} = -F+... \quad \Longrightarrow \quad F(y) \sim {\mathrm
e}^{- \frac {y}{\sqrt 2}} \cos\bigl(\frac {\sqrt 3}2 y +c\bigr),
\quad y \to + \infty.
 $}
 \eeq

%%%%%%%%%%%%%%%%%%%%%%%%%%%%%%%%%%%%%%%%%%%%
\subsection{Patterns}

\begin{figure}
%\vskip -.3cm
 \centering
\includegraphics[scale=0.8]{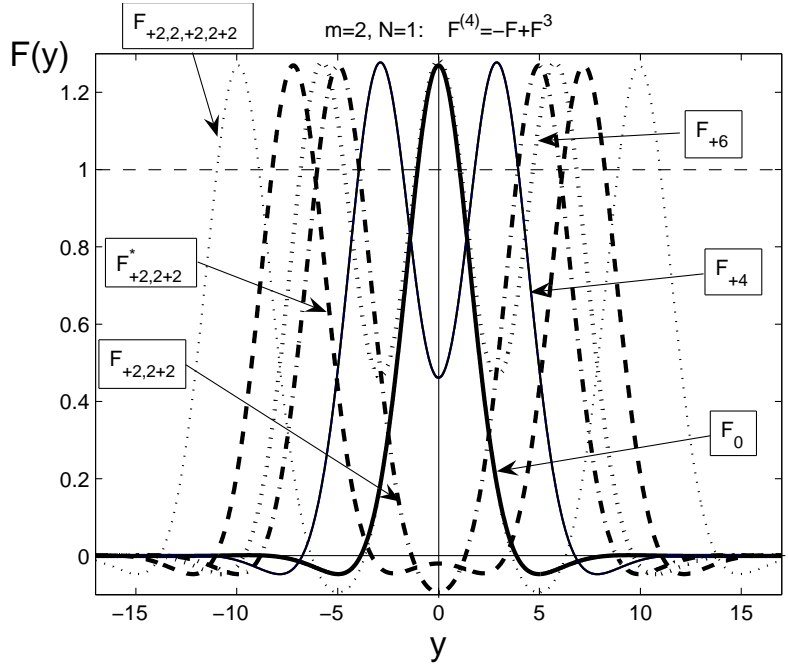}
 \vskip -.4cm
\caption{\rm\small  Various solutions  of the ODE  (\ref{nn1}):
positively dominant patterns.}
  %%% \vskip -.3cm
 \label{A2xNN1}
\end{figure}

\begin{figure}
%\vskip -.3cm
 \centering
\includegraphics[scale=0.8]{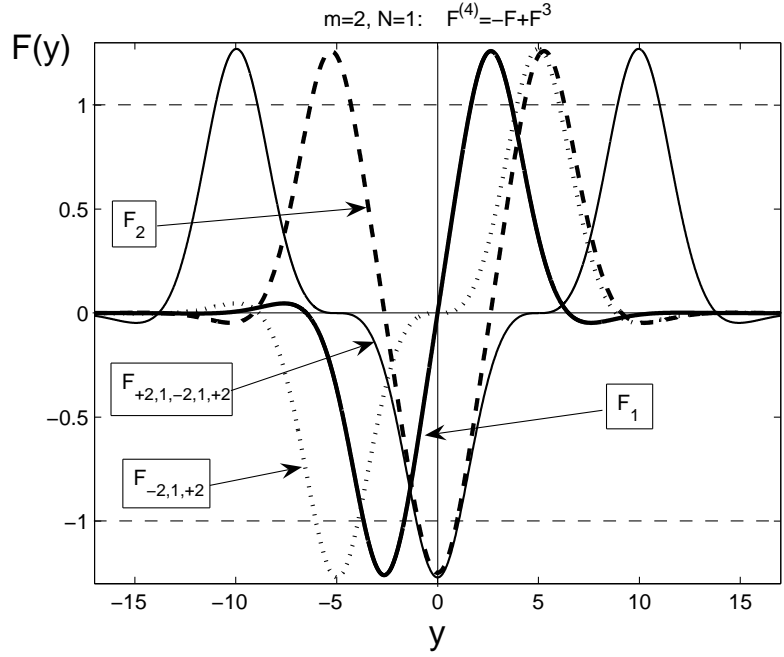}
 \vskip -.4cm
\caption{\rm\small  Various solutions  of the ODE  (\ref{nn1}):
essentially changing sign  patterns.}
  %%% \vskip -.3cm
 \label{A2xNN2}
\end{figure}

 Figures \ref{A2xNN1}
 %%(a)
  and \ref{A2xNN2}
  %%(b)
   show a few  typical patterns,
which we are already familiar with. It is important to notice
that, in Figure \ref{A2xNN1}, by watching the behaviour of small
negative solutions close to the origin $y=0$,
%%(a),
 there are
two different  patterns that can be classified as $F_{+2,2,+2}$,
and the second one is denoted by $ F^*_{+2,2,+2}$. This shows
again that the number of intersections with equilibria $\pm 1$ and
0 are not enough for a complete pattern description (in fact, this
underlines that a homotopy approach using the hodograph plane is
not applicable to equations such as (\ref{nn1E})). It is seen
there that $ F^*_{+2,2,+2}(y)$ exhibit more non-monotone structure
for $y \in (-4,4)$ than $ F_{+2,2,+2}(y)$, so that the derivative
$ F^{*'}_{+2,2,+2}(y)$ has there 3 zeros therein, while $
F'_{+2,2,+2}(y)$ has just one at $y=0$. Thus, these two patterns
can be distinguished by the {\em number of zeros of their
derivatives}, but, by the same reasons, we do not think that
counting {\em internal zeros} of the {\em pairs}
$\{F(y),\,F'(y)\}$ can help to create any rigorous Sturmian-like
classification of such patterns.

\begin{figure}
 %%\vskip -.3cm
\centering \subfigure[$F_{+10,1,-4,+6,4,+8}$]{
\includegraphics[scale=0.65]{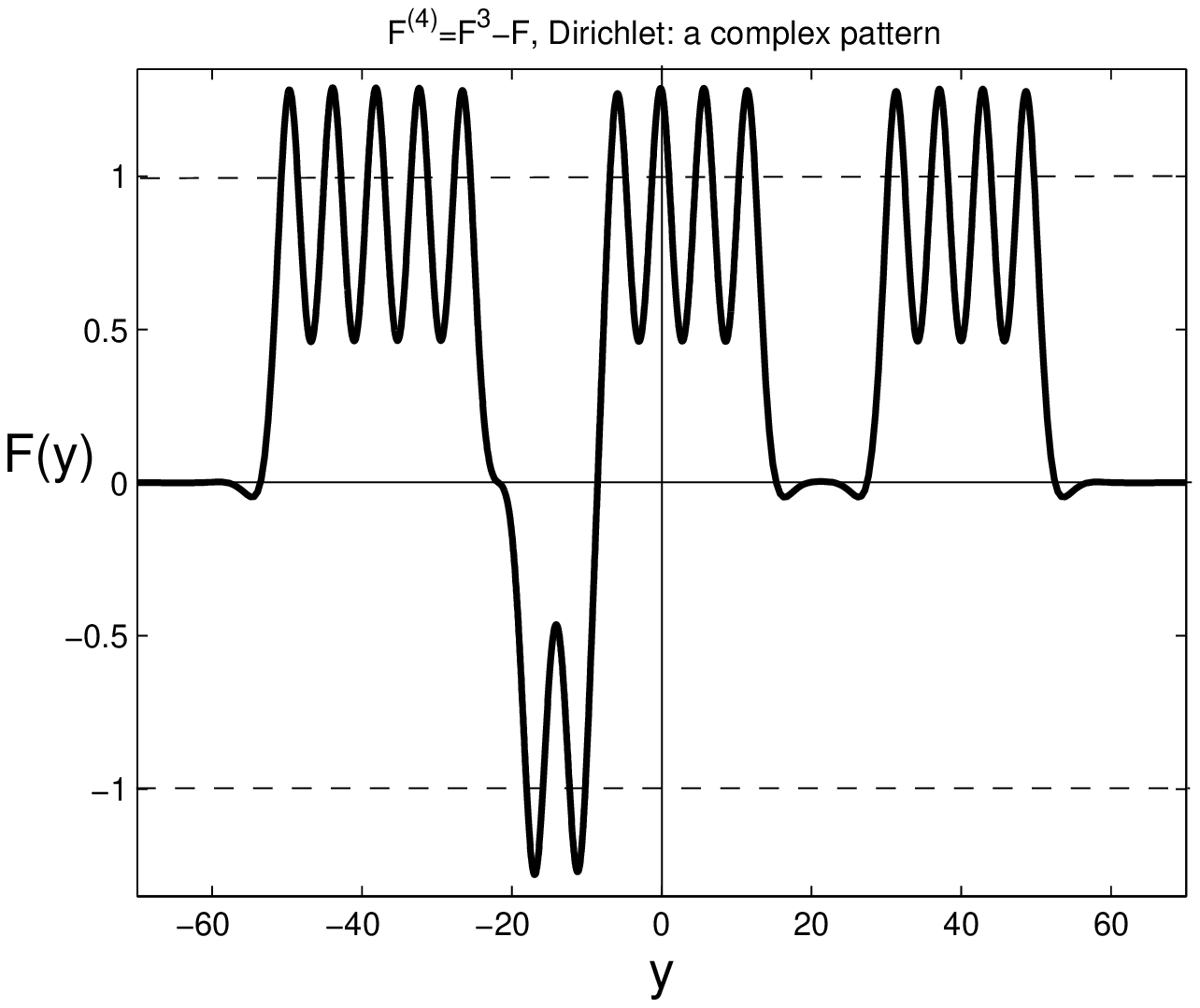}
} \subfigure[$F_{-12,1,+4,1,-12,1,+12,1,-16}$]{
\includegraphics[scale=0.65]{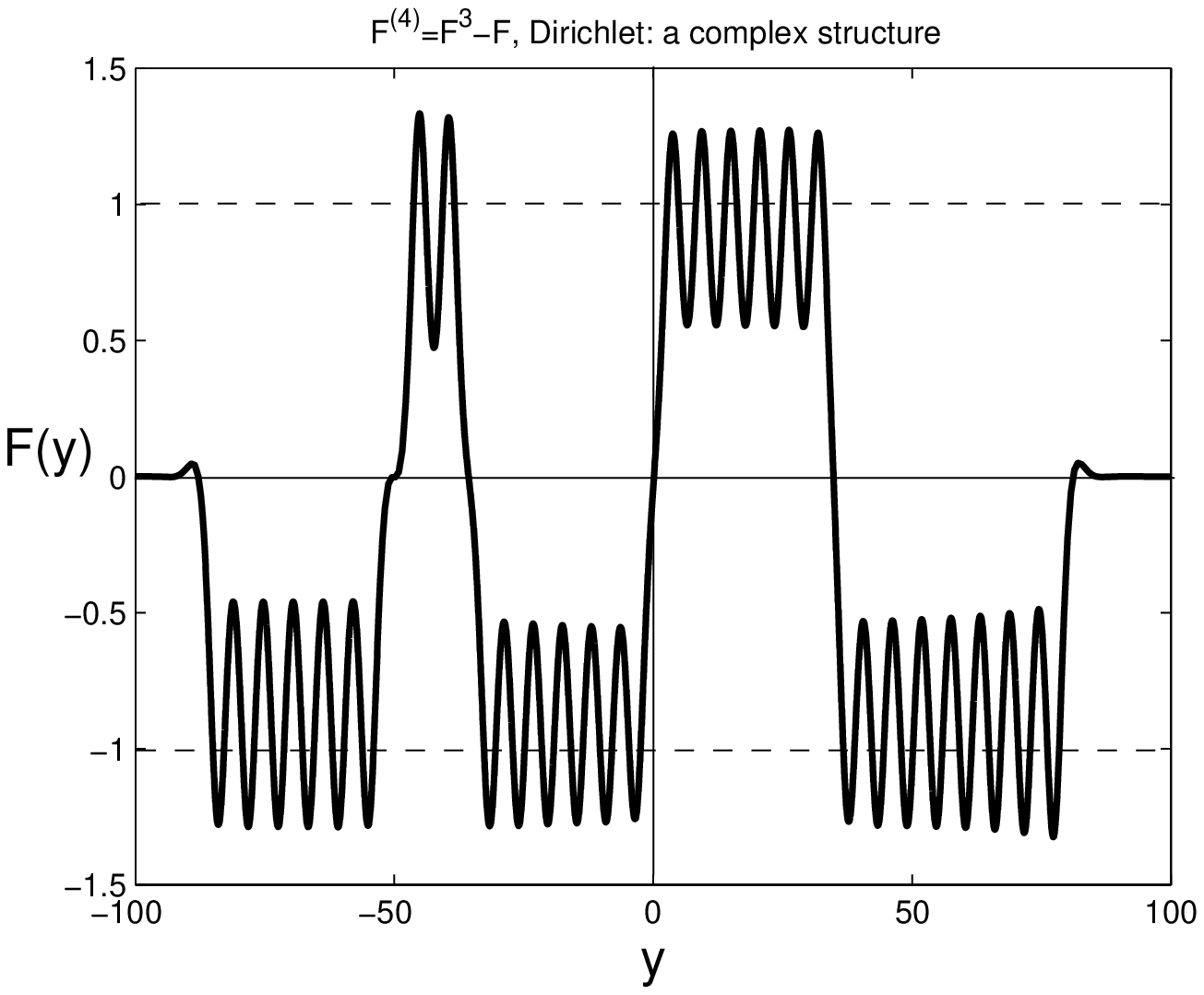}               %%%{FFP01NN.eps}
}
 \vskip -.2cm
\caption{\rm\small Examples of complicated patterns  of the ODE
(\ref{nn1}).}
 %%% for $m=4$ and $n=1$: profiles (a), and zero
%%%%structure (b).}
 \label{A1x}
  %%%%%{FGF.1fig}
\end{figure}

 In Figure \ref{A1x}, we
present typical complicated patterns
 for the problem (\ref{nn1}), which remind similar ``multi-hump"  structures
 obtained above in Figure \ref{FCmm2}.

  Note that, in view of the fast exponential decay (\ref{nn2}),
  it is difficult to observe
 by standard numerical methods  that, unlike the previous problem, the
 profiles in Figure \ref{A1x} {\em are not} compactly supported.
  Note that, for $n>0$, we succeeded
in the logarithmic scale, by taking small regularization parameter
and tolerances $\sim 10^{-12}$, to reveal the difference between
the linearized zero-behaviour (as in (\ref{nn2}))
 and the
nonlinear one at the interface as $y \to y_0^+$ of the type
($\varphi(s)$ is the oscillatory component; see Figure \ref{ZZ1}
and \cite[\S~4]{GMPSob,GMPSobIarX})
%%(\ref{2.2}).
 \beq
   \label{2.2}
    \mbox{$
   F(y) = (y-y_0)^\g \var(s) \whereA s= \ln (y-y_0) \andA
    \g =
    %%\frac{2m}\a \equiv
\frac {2m(n+1)}n>2m.
    $}
    \eeq
For $n<0$, there are no ``nonlinear zeros", but there exists an
{\em infinite number} of linearized ones given by (\ref{nn2}),
first dozens of which can be easily observed numerically.

%%%%%%%%%%%%%%%%%%%%%%%%%%%%%%%%%%%%%%%%%%%%%%%%%%
%%\begin{figure}
%\vskip -.3cm
%% \centering
%%\includegraphics[scale=0.65]{Fan1.eps}
%% \vskip -.4cm
%%\caption{\rm\small Four patterns  of the ODE in (\ref{nn1E}).}
%%   \vskip -.3cm
%% \label{A2x}
%%\end{figure}
%%%%%%%%%%%%%%%%%%%%%%%%%%%%%%%%%%%%%%%%%%%%%%%%%%%%%%
%%%%%%%%%%%%%%%%%%%%%%%%%%%%%%%%%%%%%%%%%%%%%%%%%%
%%%%%%%%%%%%%%%%%%%%%%%%%%%%%%%%%%%%%%%%%%%%%%%%%%%%%%

%%%%%%%%%%%%%%%%%%%%%%%%%%%%%%%%%%%%%%%%%%%%%%%%%%
%%%%%%%%%%%%%%%%%%%%%%%%%%%%%%%%%%%%%%%%%%%%%%%%%%%%%%

%%\begin{figure}
 %%\vskip -.3cm
%%\centering \subfigure[positive dominant patterns]{
%%\includegraphics[scale=0.52]{Anal1.eps}
%%} \subfigure[changing sign patterns]{
%%\includegraphics[scale=0.52]{Anal2.eps}               %%%{FFP01NN.eps}
%%}
%% \vskip -.2cm
%%\caption{\rm\small Various solutions  of the ODE  (\ref{nn1}).}
 %%% for $m=4$ and $n=1$: profiles (a), and zero
%%%%structure (b).}
%% \label{A2xNN}
  %%%%%{FGF.1fig}
%%%\end{figure}

%%%%%%%%%%%%%%%%%%%%%%%%%%%%%%%%%%%%%%%%%%%%%%%%%%%%%%%%%%%%%%%%%%%%%%%%%
 \subsection{Application of Lusternik--Schnirel'man and fibering theory}

Obviously, the ODE  (\ref{nn1}) (or the elliptic problem in
(\ref{nn1E})) possesses the variational setting in $(-R,R)$ (or
$B_R$), to which the same fibering version of Lusternik--Schnirel'man theory applies,
as in \cite[\S~3]{GMPSob,GMPSobIarX}.
%% Section \ref{SectVar}.
This gives a {\em countable family of basic patterns} $\{F_l\}$
for both the  ODE and the elliptic PDE. Introducing the
preliminary $h$-approximation of patterns makes it possible to
reconstruct other families of more complicated geometric
structure. Since $F_l$ are not  compactly supported, we always
assume fixing $R=R(l) \gg 1$ for $l \gg 1$. The solutions in $\re$
($\ren$) are then obtained by passing to the limit $l \to \iy$. We
do not stress the attention to such a  compactness procedure that
assumes deriving some uniform bounds independent of $l$. On the
other hand, in view of the known exponential decay of all the
solutions, the variational statement in the whole space $\ren$ is
also an option; see \cite{PelTroy}.
 We now present brief comments.

%%%%%%\com{To SIP: how to justify the var. method in $\re$?}

 %%%%%\com{Assume that we know the exponential decay at infinity}

The functional is
 \beq
 \label{V1aa}
 %% \mbox{$
 {E}(F)=  \frac 12  \displaystyle \int |\tilde D^m F|^2 + \frac 12 \displaystyle \int
 F^2 -\frac{1}{4} \, \displaystyle \int F^4.
 %%%%%% \quad \mbox{with} \,\,\, \b=\frac {n+2}{n+1}\in (1,2),
  %% \,\,
  %%%\b= \frac 1 n,
 %% $}
  \eeq
 %%% where we use the notation $\tilde D^m= \D^{\frac m2}$ for even
  %%$m$ and $\tilde D^m= \n \D^{\frac {m-1}2}$ for odd $m$.
  %%In general, we have to look for critical points in $W^2_m(\ren)
  %%\cap L^2(\ren) \cap L^{\b}(\ren)$.
 %% Bearing in mind %%% radially symmetric solutions or in 1D,
 %%compactly supported solutions, we choose a sufficiently
  %%large radius $R>0$ of the ball $B_R$ and consider the variational problem for (\ref{V1})
  %%in $W_{m,0}^2(B_R)$, where we assume Dirichlet boundary
  %%conditions on $S_R= \partial B_R$. Then
  %% both spaces $L^2(B_R)$ and $L^{p+1}(B_R)$ are
The  {spherical fibering} \cite[\S~3]{GMPSob,GMPSobIarX}
 %% \beq
 %%\label{f1}
 %%F= r(v) v \quad (r \ge 0),
 %% \eeq
%% (\ref{f1}),
 %%% \beq
 %%%\label{f1aa}
    %%%\quad  %%%(r \ge 0),
 %% \eeq
  with $v$ belonging to the
  %%% a subset
  %% in  $W_{m,0}^2(B_R)$, which is now a
    unit sphere,
   \beq
   \label{f2aa}
    \mbox{$
    F= r(v) v, \quad \mbox{where} \,\,\, v \in
    {\mathcal H}_0=\bigl\{ %%%v \in W_{m,0}^2(B_R): \,\,\,
    H_0(v)
     \equiv    \displaystyle \int |\tilde D^m v|^2 +  \displaystyle \int
 v^2 =1\bigr\},
    $}
    \eeq
    leads to the
 functional
 \beq
 \label{f3aa}
  \mbox{$
H(r,v)= \frac 12 \, r^2 - \frac 1{4}\, r^{4} \displaystyle \int v^4.
 $}
  \eeq
This attains the absolute maximum at
 \beq
 \label{f31aa}
  \mbox{$
 H'_r \equiv r-  r^{3} \displaystyle \int v^{4} =0
  \,\,\Longrightarrow \,\,
 %%$}
 %% \eeq
 %% at
 %% $$
 %%  \mbox{$
   r_0(v)=\frac 1{\sqrt{\displaystyle \int v^4}},
   $}
   \eeq %%$$
  at which
   $
    %%\mbox{$
   H(r_0(v),v)=\frac 1{2\sqrt{\displaystyle \int v^4}}.
   %%% $}
    $
 This defines the positive homogeneous convex functional
 \beq
 \label{f4aa}
 \mbox{$
 \tilde H(v) = \bigl[ \frac 1{2H(r_0(v),v)}\bigr]^2
  \equiv \displaystyle \int v^{4}.
  $}
  \eeq
  Here  Lusternik--Schnirel'man theory applies in its classic form  \cite[p.~387]{KrasZ} giving
  a countable  set $\{c_k\}$ of critical values and points denoted by
   $\{v_k\}$ (see full details in \cite[\S~3]{GMPSob,GMPSobIarX}):
 \beq
 \label{ck1aa}
  \mbox{$
 c_k = \sup_{{\mathcal F} \in {\mathcal M}_k} \,\, \inf_{v \in {\mathcal
 F}} \,\, \displaystyle \int v^4.
  $}
  \eeq
Now the sequence of critical values is decreasing,
 \beq
  \label{decr1}
 c_1 \ge c_2 \ge ... \ge c_k \ge c_{k+1} \ge ... \, .
  \eeq
  For each $v_k$ (a solution $F_k$), the critical values are given by
 \beq
 \label{dd1aa}
  \mbox{$
   v= C F \in {\mathcal H}_0 \quad \Longrightarrow \quad
 c_F \equiv \tilde H(v) = \frac { \displaystyle \int F^4}{\left(\displaystyle \int |\tilde D^m F|^2 + \displaystyle \int
  F^2 \right)^{2}} \, .
   $}
   \eeq
For all the patterns shown in Figures \ref{A2xNN1}
 and \ref{A2xNN2},
%%(a) and  (b),
 these values (Lusternik--Schnirel'man critical or not) are presented in
Table 1. As usual, this table makes it possible to detect the Lusternik--Schnirel'man
critical points that deliver critical values (\ref{ck1aa}) for
each genus. We  then again claim
 that the basic patterns $\{F_l\}$ deliver
all the Lusternik--Schnirel'man critical values (\ref{ck1aa}) with $k=l+1$.

Recall  the {\em Formal Rule of Composition} (FRC) from
\cite[\S~5]{GMPSob,GMPSobIarX}:
 %%\noi{\bf Formal Rule of
 %%Composition} (FRC) {\bf of patterns:}
  {performing  maximization of $\tilde H(v)$ of
 any $(k-1)$-dimensional manifold ${\mathcal F} \in {\mathcal
 M}_k$,}
  \beq
  \label{comp1}
  %% \fbox{$
\mbox{the Lusternik--Schnirel'man point $F_{k-1}(y)$ is obtained by minimizing all
internal tails and zeros,}
  %%%$}
 \eeq
{i.e., making the minimal number of internal transversal zeros
between single structures.
   }
Then
 (\ref{comp1}) also applies, since by the same
reason diminishing a small tail between two $F_0$-structures will
increase the corresponding value $c_F$ in (\ref{dd1aa}).

%%%%%%%%%%%%%%%%%%%%%%%%%%%%%%%%%%%%%%%%%%%%%%%%%%%%%%%%%%%%%%%%%
\begin{table}[h]
\caption{Critical values  of $\tilde H(v)=\frac { \displaystyle \int F^4}{\left(\displaystyle \int
|\tilde D^m F|^2 + \displaystyle \int
  F^2 \right)^{2}}$} %%\lineup
%%\begin{indented}
 %%\item[]
 \begin{tabular}{@{}lll}
%%%%%\vskip -.5cm
%%%%%% \\\hline
 %%\br
 $F$ & $c_F$
 \\\hline
%%% \mr
  $F_0$ & $0.2033...=c_1$\\
  $F_1$ & $0.1080...=c_2$
\\ $F_{+2,2,+2}$ & $0.1019...$\\
 $\tilde F_{+2,2,+2}$ & $0.1017...$\\
$\tilde F_{-2,1,+2}$ & $0.1014...$
\\
$F_{+4}$ & $0.0961...$\\
 $F_{2}$ & $0.0736...=c_3$\\
 $F_{+2,2,+2,2,+2}$ & $0.0680...$\\
$\tilde F_{+2,1,-2,1,+2}$ & $0.0675...$\\
 $F_{+6}$ & $0.0629...$\\
 %% \br
\end{tabular}
%%\end{indented}
\end{table}
%%%%%%%%%%%%%%%%%%%%%%%%%%%%%%%%%%%%%%%%%%%%%%%%%%%%%%%%%%%%%%%%%%%%%%%

%%%%%%%%%%%%%%%%%%%%%%%%%%%%%%%%%%%%%%%%%%%%%%%%%%%%%%%%%%%%%%
\subsection{Homotopic connections with Sturm-ordered linear eigenvalue problems}

An extra advantage of the problem (\ref{nn1}) is that the
homotopic connections of basic patterns can be revealed  easier.
Namely, as above, we consider the ODE  (\ref{nn1}) in $(-R,R)$
with sufficiently large $R >0$ to see the first $l$ patterns which
are still exponentially small for $y \approx \pm R$. Consider the
following homotopic path with the operators:
 \beq
 \label{nn3}
{\bf A}_\e(F) \equiv F^{(4)} + \e F- F^3=0.
 \eeq
 Consider the corresponding linearized operator:
 \beq
 \label{nn4}
{\bf A}_\e'(0)= D_y^4 + \e I.
 \eeq
 Let $\s(D_y^4)=\{\l_l>0, \, l\ge 0\}$ be the discrete spectrum of simple eigenvalues
  of $D_y^4>0$ in
 $L^2((-R,R))$ with the Dirichlet boundary conditions.
 %%It is well known that the corresponding
  The orthonormal eigenfunctions
$\{\psi_l\}$ satisfy Sturm's zero property; we again refer to
\cite{Elias} for most general results.

\begin{figure}
 %%\vskip -.3cm
\centering \subfigure[$\e \sim 1$]{
\includegraphics[scale=0.52]{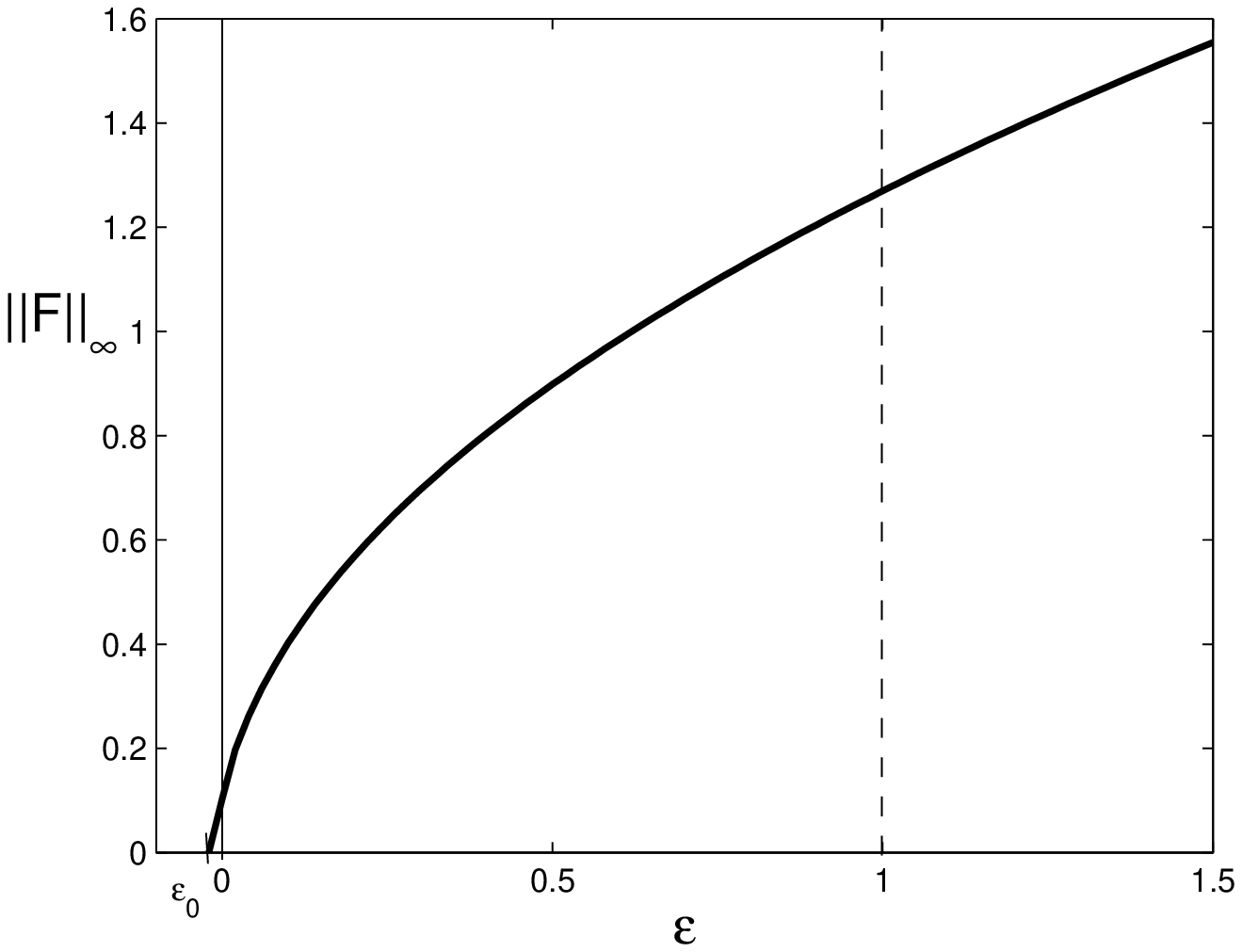}
} \subfigure[$\e \sim 100$]{
\includegraphics[scale=0.52]{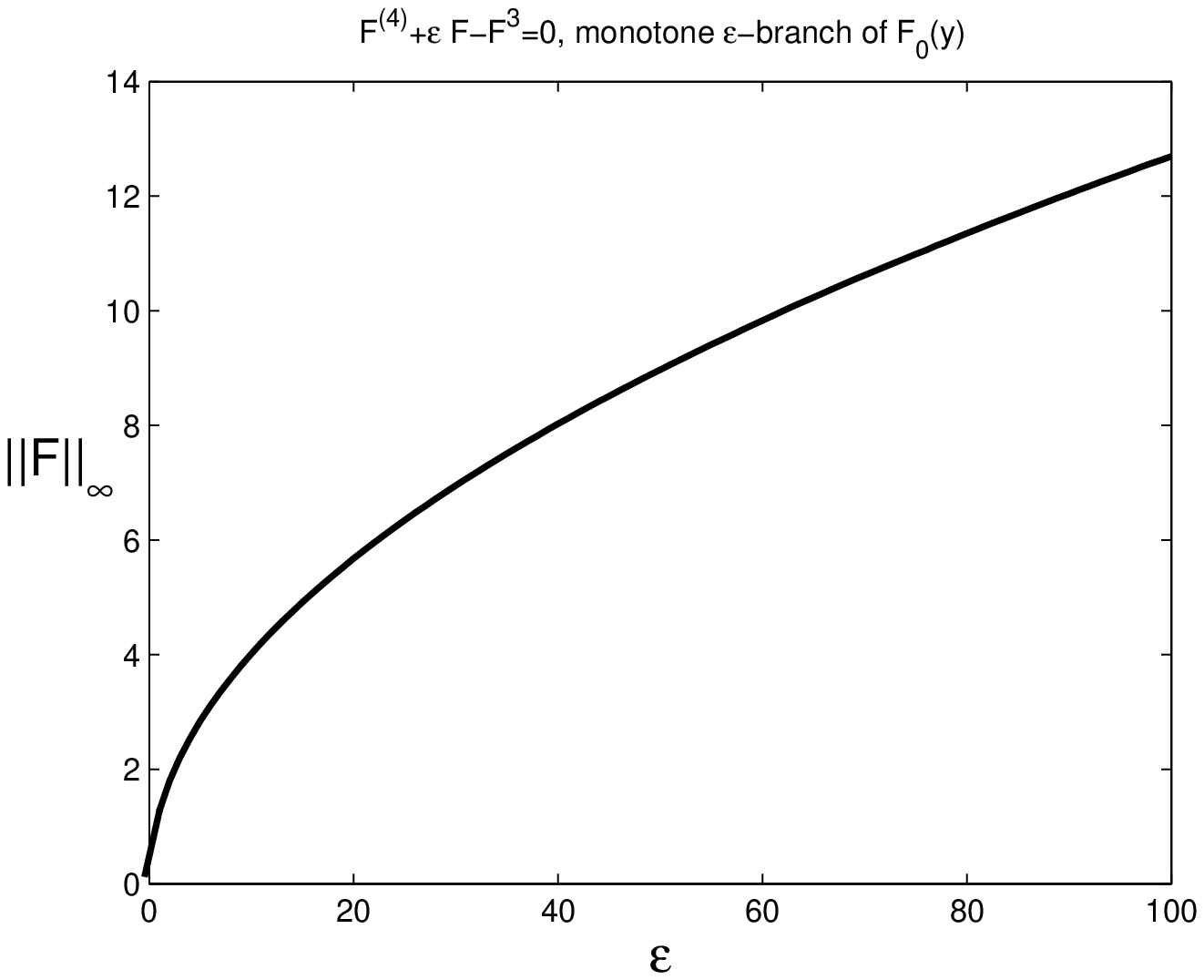}               %%%{FFP01NN.eps}
}
 \vskip -.2cm
\caption{\rm\small The monotone $\e$-branch of profile $F_0$ via
homotopy (\ref{e1}) on different scales; $\e \in(-\e_0,1.5)$ (a)
and $\epsilon \in (0,100)$ (b).}
 %%% for $m=4$ and $n=1$: profiles (a), and zero
%%%%structure (b).}
 \label{FBif}
  %%%%%{FGF.1fig}
\end{figure}

 By classic bifurcation theory \cite[p.~391]{KrasZ}, for
 such variational problems,
  $$
  \e_l=- \l_l < 0 \forA l=0,1,2,...
   $$
    are
 bifurcation points, so there exists a countable number of
 branches emanating from  these points (but we take into account
 the
 first ones).
 In order to identify the type of bifurcations, in a standard manner, setting
 $\e=\e_l+s$ for $|s|$ small and
 $$
 F= C \psi_l+w, \quad \mbox{where} \quad  w \bot \psi_l,
  $$
 substituting into (\ref{nn3}) and multiplying by $\psi_l$
 yields
  $$
   \mbox{$
   s = C^2(s) \displaystyle \int \psi_l^4 +...>0 \quad \Longrightarrow \quad s>0
   \,\,\, \mbox{and}
   \,\,\, C(s)= \pm \sqrt{\frac s{ \displaystyle \int \psi_l^4}} +... \, .
   $}
   $$
Hence, at $\e=\e_l$, there appear two branches with the equations
 $$
  \mbox{$
 F_l(y)= \pm \sqrt{\frac {\e-\e_l}{ \displaystyle \int \psi_l^4}}\,\, \psi_l(y)+... \forA
 \e>\e_l.
  $}
 $$
 Therefore,
   all the bifurcations are {\em pitchfork} and
 are {\em supercritical}, i.e., two symmetric branches are initiated at $\e=
 \e_l^+$; see Figure \ref{FBif}(a), where we construct numerically the first positive $\e$-branch
  of $F_0$ of the
 $\e$-bifurcation diagram and show that this branch is extensible
 up to the necessary value $\e=1$, and even up to $\e \sim 10^2$ (b), and further.

%%%%%%%%%%%%%%%%%%%%%%%%%%%%%%%%%%%%%%%%%%%%%%%%%%%%%%%%%%%%
%%%%%%%%%%%%%%%%%%%%%%%%%%%%%%%%%%%%%%%%%%%%%%%%%%%%%%%

 %%%%%%%%%%%%%%%%%%%%%%%%%%%%%%%%%%%%%%%%%%%%%%%%%%
%%%\begin{figure}
%\vskip -.3cm
%% \centering
%%\includegraphics[scale=0.65]{F0epsN.eps}
%% \vskip -.4cm
%%\caption{\rm\small The $\e$-branches of profile $F_0$  via
%%homotopy (\ref{e1}).}
%%   \vskip -.3cm
%% \label{FBif}
%%\end{figure}
%%%

%%\begin{figure}
 %%\vskip -.3cm
%%\centering \subfigure[$\e$-branch of $F_0$, $R=8$]{
%%\includegraphics[scale=0.52]{F0epsN.eps}
%%} \subfigure[$\e$-branch of $F_1$, $R=8$]{
%%\includegraphics[scale=0.52]{F1epsN.eps}               %%%{FFP01NN.eps}
%%}
%% \vskip -.2cm
%%\caption{\rm\small $\e$-branches of profiles $F_0$ and $F_2$ via
%%homotopy (\ref{e1}).}
 %%% for $m=4$ and $n=1$: profiles (a), and zero
%%%%structure (b).}
%% \label{FBif}
  %%%%%{FGF.1fig}
%%%\end{figure}

It follows %% from (\ref{f56})
that all the branches are originated at $\e_l <0$, so  being
continued up to $\e=1$, give the original equation  (\ref{nn1}).
 The questions of global continuation of branches are classic in
 nonlinear variational theory; see \cite[\S~6.7C]{Berger}.
 The global behaviour of bifurcation branches
for  $2m$th-order ODEs with  analysis of possible types of
end-points is addressed in \cite{Rynn2}. These results hardly
apply to the equation (\ref{e1}) with non-coercive operators
admitting solutions of changing sign near boundary points.

%% To show that branches cannot have {\em
%%turning points} on interval $\e \in (\e_l,1)$, we argue by
%%contradiction.
 %%follows from the straightforward
%%%observation that if, in such circumstances,

The existence of a turning point of the given branch in this real
self-adjoint case, i.e., of a saddle-node bifurcation,  assumes
that there exists an eigenvalue (say, simple)
 $$
 0 \in \s({\bf A}_\e'(F)), \quad \mbox{i.e.,} \quad \exists \,\, \phi_0: \,\,
   \phi_0^{(4)} + \e \phi_0 - 3 F^2 \phi_0=0,
 $$
 where $\phi_0$ is an eigenfunction of ${\bf A}_\e'(F)$
  %%% $$
 %% \phi_0^{(4)} + \e \phi_0 - 3 F^2 \phi_0=0
%%%  $$
  satisfying the Dirichlet conditions $\phi_0=\phi_0'=0$ at $y =
  \pm R$. For a moment, we digress from our difficult ODEs and
  consider simpler models with known bifurcation diagrams.

\ssk

%%%%%%%%%%%%%%%%%%%%%%%%%%%%%%%%%%%%%%%%%%%%%%%%%%%%%%%%%%%%%%%%%%%%%
\noi{\bf Remark 1: on turning saddle-node bifurcation points of
positive solutions.}  Such turning points do  exist for equations
with other nonlinearities and another dependence on $\e$, e.g.,
 \beq
  \label{jj1}
 (-1)^{m} F^{(2m)} = \e {\mathrm e}^F \quad \mbox{or}
  \quad  (-1)^{m} F^{(2m)} = \e (1 + F^3) \,\,\,(F \ge 0),
   \eeq
 with Dirichlet boundary conditions at $\pm R$. Existence of two
 different branches of solutions (they are positive by \cite{Elias})
  for all small $\e>0$ is established by the
 fibering method as in \cite[\S~3]{GMPSob,GMPSobIarX}.%%%% Section \ref{SectVar}.

  On the other hand,
 for  $\e \gg 1$, positive solutions of (\ref{jj1})
 are obviously not possible.
 This is easily seen by multiplying the first equation in
 (\ref{jj1}) by the first eigenfunction $\psi_0>0$ with the eigenvalue $\l_0>0$ of the positive
 operator $(-1)^m D_y^{2m}$. This gives for the first Fourier coefficient $C_0=\displaystyle \int F \psi_0$ the
 following inequality:
 \beq
 \label{J1}
 \mbox{$
 %%\begin{matrix}
 \l_0 C_0 = \e \displaystyle \int {\mathrm e}^F \psi_0 \equiv
\e \|\psi_0\|_1 \displaystyle \int {\mathrm e}^F \frac { \psi_0}{\|\psi_0\|_1}
 %%\ssk\\
 \ge \e \|\psi_0\|_1 \,\, {\mathrm e}^{\displaystyle \int F \frac {
 \psi_0}{\|\psi_0\|_1}}=
\e \|\psi_0\|_1 \,\, {\mathrm e}^{ \frac {
 C_0}{\|\psi_0\|_1}}.
 %% \end{matrix}
  $}
  \eeq
  At the last stage we have used Jensen's inequality for the convex
  function  ${\mathrm e}^F$.
One can see that the resulting  inequality in (\ref{J1}),
 $$
 \mbox{$
 \l_0 C_0 \ge \e \|\psi_0\|_1 \,\, {\mathrm e}^{ \frac {
 C_0}{\|\psi_0\|_1}},
  $}
  $$
 does not
have a solution $C_0 >0$ for all $\e \gg 1$, meaning nonexistence.

By standard theory of compact integral Uryson--Hammerstein
operators \cite{Kras, KrasZ} and branching theory
\cite{VainbergTr}, the solutions of (\ref{jj1}) detected for small
$\e>0$ comprise two continuous branches, which  must end up at a
saddle-node bifurcation point (no other bifurcations are possible)
at some
 $$
 \e=\e_{\rm s-n}>0.
  $$

\begin{figure}
 \centering
 \psfrag{||F||}{$\|F\|_\infty$}
 \psfrag{e}{$\e$}
 \psfrag{t2}{$t_2$}
 \psfrag{t3}{$t_3$}
 \psfrag{t4}{$t_4$}
  \psfrag{v(x,t-)}{$v(x,T^-)$}
  \psfrag{final-time profile}{final-time profile}
   \psfrag{tapp1}{$t \approx 1^-$}
\psfrag{x}{$x$}
 \psfrag{0<t1<t2<t3<t4}{$0<t_1<t_2<t_3<t_4$}
  \psfrag{0}{$0$}
 \psfrag{l}{$l$}
 \psfrag{-l}{$-l$}
\includegraphics[scale=0.5]{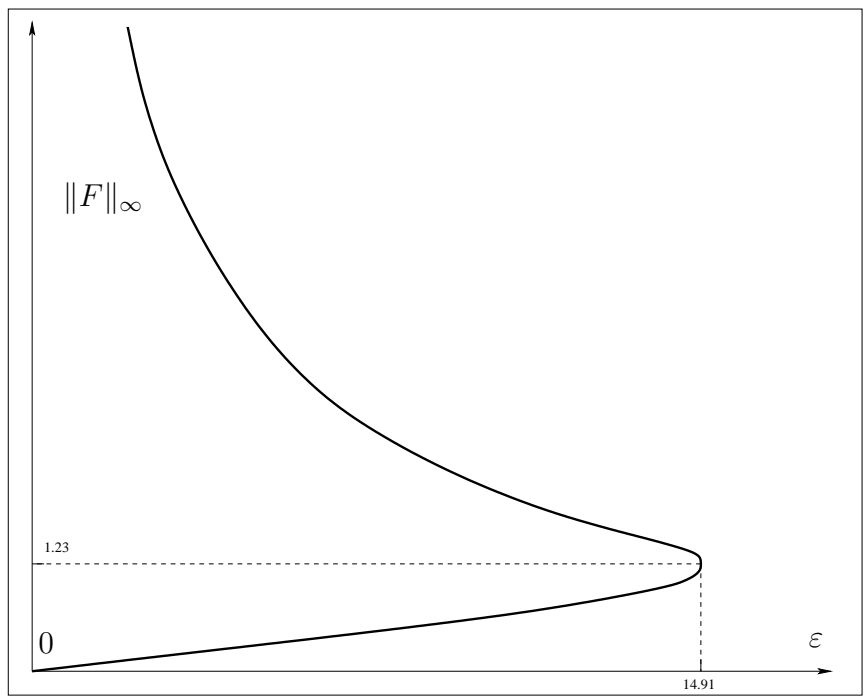}     %%%{F1EE1.pdf}
\caption{\small The bifurcation diagram for (\ref{qq1}).}
     \vskip -.3cm
 \label{SN1}
\end{figure}

In Figure \ref{SN1}, the global bifurcation diagram is presented
for the quadratic equation,
 \beq
 \label{qq1}
  F^{(4)}= \e (1+F^2) \quad \mbox{on \, $(-1,1)$} \whereA \fbox{$\e_{\rm s-n}=14.91...$} \, .
   \eeq
The upper branch blows-up to $+\infty$ as $\e \to 0^+$, while the
lower one vanishes. See \cite{Kor04} for a different approach to
bifurcation analysis for such fourth-order ODEs.

\ssk

%%%%%%%%%%%%%%%%%%%%%%%%%%%%%%%%%%%%%%%%%%%%%%%%%%%%%%%%%%%%%%%%%%%%%%%
%%%%%%%%%%%%%%%%%%%%%%%%%%%%%%%%%%%%%%%%%%%%%%%%%%%%%%%%%%%%%%%%%%%%%%%

%%%%%%%%%%%%%%%%%%%%%%%%%%%%%%%%%%%%%%%%%%%%%%%%%%%%%%%%%%
%%\begin{figure}
 %%\vskip -.3cm
%%\centering \subfigure[lower, stable branch]{
%%\includegraphics[scale=0.52]{LB1.eps}
%%} \subfigure[upper, blowing-up branch]{
%%\includegraphics[scale=0.52]{UB1.eps}               %%%{FFP01NN.eps}
%%}
%% \vskip -.2cm
%%\caption{\rm\small The bifurcation diagram for (\ref{qq1}); the
%%lower branch (a) and the upper one (b).}
%% \label{SN1}
  %%%%%{FGF.1fig}
%%\end{figure}

We return to the  equation (\ref{nn3}).
   Nonexistence of such saddle-node bifurcation points for the ODEs (\ref{nn3}) is a difficult
open problem. Numerically, we have got a strong evidence that the
first branches are strictly monotone without turning points; see
Figure \ref{FBif}(b), where the first branch is  extended up to
$\e =100$.

On the other hand, it is easy to check by the fibering approach
that (\ref{nn3}) has a countable set of solutions for arbitrarily
large $\e>0$, which are also identified as continuous curves by
nonlinear compact operator theory. Hence, there are infinitely
many branches of basic patterns that are unboundently extensible
in $\e$. We then conclude that at $\e=1$ all these branches are
available. Hence,  at $\e=1$, the corresponding solutions $F_l$ of
(\ref{nn1}) inherit
 their Sturm index from eigenfunctions $\psi_l$ that occur at $\e
 = \e_l^+$.

 %%%%\com{To SIP: enough ?  Or in greater detail ?}

  %% Then differentiating (\ref{nn1}) necessarily yields
 %%%$\phi_0 \sim F'(y)$, which does not satisfy boundary conditions,
 %%since $F'''(\pm R) \not =0$ (otherwise, $F(y) \equiv 0$).
%%It also follows that branches cannot be attracted to further
%% supercritical bifurcation point, since, obviously, this would mean
%% appearance a turning point beforehand.

%% \com{To SIP: Correct ???  Should be standard ??? More ???? 0 not simple ? Possible ?}

Of course, we can use the alternative homotopy approach.
  Figure \ref{FBif}(a)  shows
 that the $\e$-branches are well-defined at $\e=0$, where
 we obtain the simpler equation
   $$
 -F^{(4)} + F^3=0
  $$
 that, as in (\ref{m1}), admits Sturm's classification of all the
patterns.

%%Therefore, unlike the  non-Lipschitz problem (\ref{4.3}), the
%%previous $\e$-homotopy approach does not
%% apply to the analytic ODE (\ref{nn1}).

 %%%%that all these three types of problems and solutions have the
 %%%same classifying Sturm index.

%%% These branches being fully extended in $\e$, at $\e=1$
%%generate solutions of the familiar equation (\ref{m1}) with
%%Sturm's property,
%% \beq
%% \label{m12}
%% F^{(4)} + F^3=0,
%%  \eeq
%%  and after that, for $\e=0$,  give (\ref{nn1}).

%%%\com{To Monotonicity in $\e$?  Obviously yes??? Not easy ?}

%%%The fibering method

\ssk

\noi{\bf Remark 2: explicit global monotone  $\e$-bifurcation
branches for a non-local nonlinearity.} Computations of branches
are straightforward for the following non-local equation with the
analytic cubic nonlinearity (cf. (\ref{nn3})): %%%%% and (\ref{V1N})):
 \beq
  \label{non1}
   \mbox{$
F^{(4)}=- \e F + \bigl(\displaystyle \int F^2 \bigr) F \quad \mbox{in} \quad
(-R,R).
 $}
\eeq
 Then the solutions %%% (cf. (\ref{pp2}))
 $$
  \mbox{$
 F_l(y) = \pm \sqrt{\l_l+\e} \,\, \psi_l(y)\quad \mbox{for}
 \quad  l=0,1,2,...
  $}
 $$
are originated at $\e_l=-\l_l^+$, and the branches have  the geometric form as in
Figure \ref{FBif}.

%%%%%%%%%%%%%%%%%%%%%%%%%%%%%%%%%%%%%%%%%%%%%%%
%%%%%%%%%%%%%%%%%%%%%%%%%%%%%%%%%%%%%%%%%%%%%%%%%%%
 \section{\underline{\bf Problem ``Sturm Index":} using $R$-compression}
  \label{Sect8}

  Consider again for a while
   the simplest non-Lipschitz  ODE for $m=2$:
  \beq
 \label{4.3}
  \mbox{$
 F^{(4)}=F- %%%%\frac 1n \,
 \big|F\big|^{-\frac n{n+1}} F \quad \mbox{in} \,\,\,
 \re.
 $}
 \eeq
As we have seen, the oscillatory character of the solutions close
to interfaces, finite or infinite, causes some difficulties in
determining the actual (generalized) Sturm index of the patterns
defined as the number of certain oscillations, which are
``dominant" zeros (or extrema) points about 0. Oscillations about
non-trivial equilibria $\pm 1$ are clear and do not exhibit such
difficulties,
 though sometimes it is rather hard  to distinguish these classes of non-monotonicities. The main
problem of concern is that how to distinguish the ``dominant"
non-monotonicities and ``small" non-dominant ones related to
asymptotic  oscillations close to interfaces, which are almost the
same for all the patterns; cf. \cite[\S~4]{GMPSob,GMPSobIarX} for
the non-Lipschitz nonlinearity and (\ref{nn2}) for the analytic
one (\ref{nn1}).

As we know, in  the CP for the PDEs involved, it is enough to pose
the Dirichlet problem for (\ref{S2}) on sufficiently large
intervals $(-R,R)$, so one does not need to take infinite length
$R=+\infty$. Let is concentrate on the case $m=2$,
 \beq
 \label{4.3R}
 \mbox{$
 F^{(4)}=F- %%%%% \frac 1n \,
 \big|F\big|^{-\frac n{n+1}} F \quad \mbox{in} \,\,\,
 (-R,R),
 $}
 \eeq
 with Dirichlet conditions at the end-points $y= \pm R$. Given a
 compactly supported pattern %%%% in $\re$,
   \beq
   \label{pp1}
    \mbox{$
  F(y)=  F_\s(y), \quad \mbox{with multiindex} \quad
    \s=\{\s_1, \s_2, \s_3,..., \s_l\}
    \quad \big({\rm supp}\, F_\s \subset (-R,R)\big),
     $}
    \eeq
we perform its $R$-{\em compression}, i.e.,  start  to decrease
$R$ observing a {\em continuous} deformation of the corresponding
profile $F(y;R)$ until the minimally possible value
 $$
 R=R_{\rm min} > 2R_*>0
  \quad (\mbox{$R_*$ is defined in \cite[\S~3]{GMPSob,GMPSobIarX}}).
 %%(\ref{4.3NN})}).
 $$
 %%%%%%where $R_*$ is defined in (\ref{4.3NN}).
%% (see Proposition
%%%\ref{Pr.Sm}).

 %%%Here the lower bound $R_0$ is obtained from the linearized problem:  $D_y^4$ has the
 %%%eigenvalue $\l_0=1$ and the eigenfunction $\psi_0>0$ on
 %%%$(-R_0,R_0)$.
Let $\s_{\rm min}$ be the multiindex of $F(y;R_{\rm min})$, if
this profile is bounded. If not, we mean the index $F(y;R_{\rm
min}^+)$ calculated for profiles with $R \approx R_{\rm min}^+$.
As usual, $\s_{\rm min}$ reflects the sequence of intersection
numbers with equilibria $\pm 1$, without taking into account those
with 0, which are, actually, nonexistent; see below.
 There is a direct relation to Definition \ref{D.1} to be explained later
 on:

%%In other words,  if $R_{\rm min} = 0$, by  $\s_{\rm min}$ we mean
%%the multiindex of
%% the steep profile $F(y;R)$ with any sufficiently small $R>0$.

\begin{definition}
 \label{Def.21}
 Given a compactly supported solution $(\ref{pp1})$ of
 $(\ref{4.3R})$, by its generalized Sturm index we mean the multiindex
 $\s_{\rm min}$  obtained by the $R$-compression.
  \end{definition}

Two cases of $R$-compression are distinguished:

\ssk
 (i) As $R \to R_{\rm min}^+$, the profile $F(y;R)$ gets
unbounded and achieves a clear geometric form with $l$ extrema.
Then we say that this $l$ is precisely the {\em Sturm index}
$l=I_{\rm S}$ of the functions $F$ on this $R$-branch. As usual,
we claim that such an index $l$ can be attributed to the basic
family $\{F_l\}$ only.

\ssk

(ii) There exists a finite limit $F(y;R_{\rm min})$. Then the only
{\em generalized Sturm index} can be defined as a characterization
of its geometric structure with fixed numbers of intersections
with equilibria $\pm 1$.

\ssk

\begin{figure}
 %%\vskip -.3cm
\centering \subfigure[$F_{0}$]{
\includegraphics[scale=0.65]{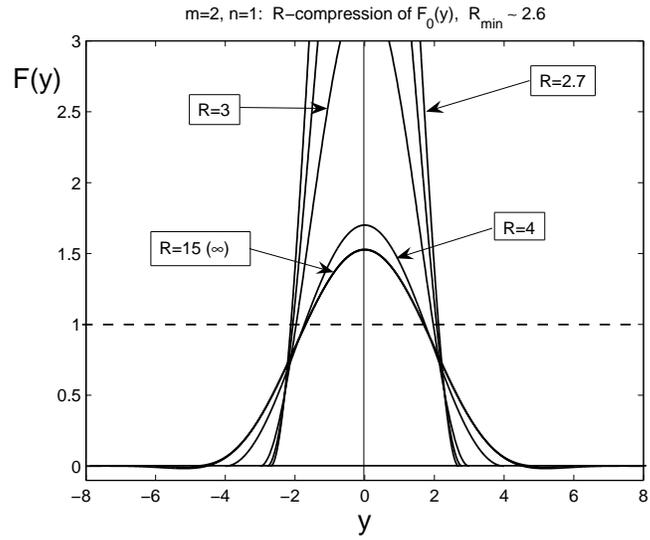}
} \subfigure[$F_{+2,2,+2}$]{
\includegraphics[scale=0.7]{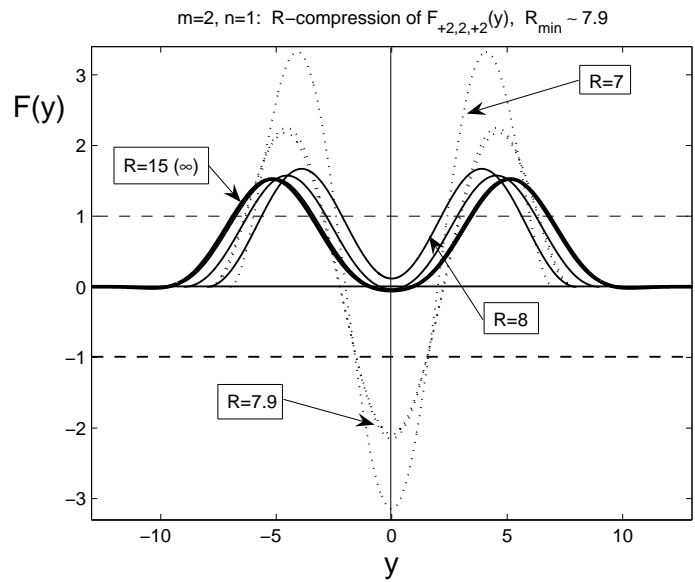}               %%%{FFP01NN.eps}
}
 \vskip -.2cm
\caption{\rm\small The $R$-compression of solutions of the ODE
(\ref{4.3}) for $n=1$;  $F_0$ (a) and $F_{+2,2,+2}$ (b).}
 %%% for $m=4$ and $n=1$: profiles (a), and zero
%%%%structure (b).}
 \label{FCom1}
  %%%%%{FGF.1fig}
\end{figure}

  In Figure \ref{FCom1}, we present some numerical results of the
  $R$-compression of the profile $F_0$ (a) (interfaces become non-oscillatory
  already for $R \sim 5$, $\s_{\rm min}=\{+2\}$
  as expected, i.e., with  no essential ``transversal" zeros and a single maximum) and
  $F_{+2,2,+2}$ (b)
   (non-oscillatory interfaces for $R \sim 10$, $\s_{\rm min}=\{+4\}$, no essential zeros and three
  extrema). In (b), by dotted lines, we denote some other profiles
  that also appear for such $R \approx R_{\rm min}^+$. These belong to other branches of
  solutions of (\ref{4.3R}).

%%%%%%%%%%%%%%%%%%%%%%%%%%%%%%%%%%%%%%%%%%%%%%%%%%%%%%%%%%

\begin{figure}
%\vskip -.3cm
 \centering
\includegraphics[scale=0.75]{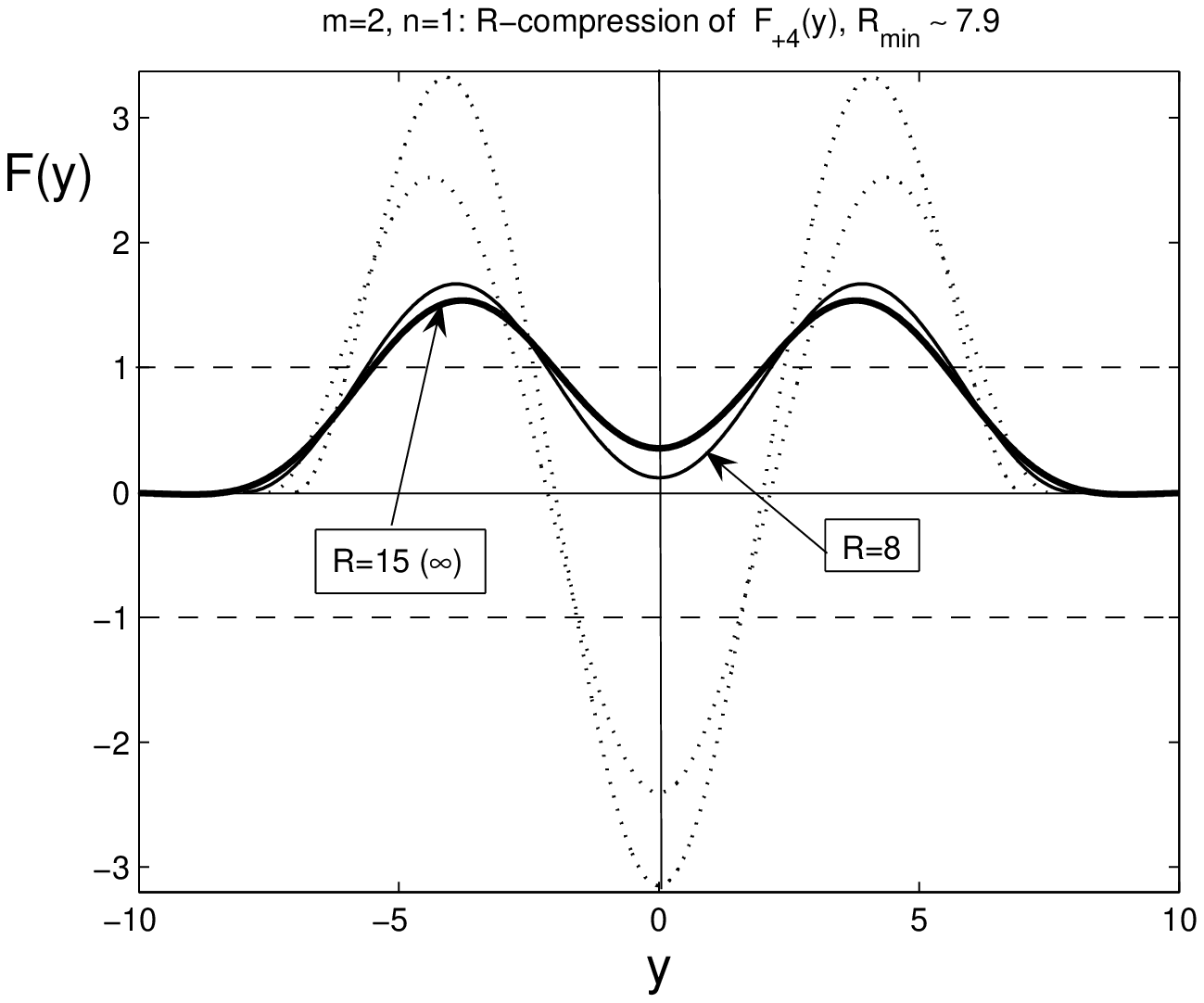}
 \vskip -.4cm
\caption{\rm\small The $R$-compression of the solution $F_{+4}$ of
the ODE (\ref{4.3}) for $n=1$.}
   \vskip -.3cm
 \label{FCom2}
\end{figure}

In Figure \ref{FCom2}, we show the $R$-compression of a different
profile,
 $F_{+4}$, with a similar structure. Note that,
 in Figures \ref{FCom2} and \ref{FCom1}(b), the profiles for $R=8$
  coincide. This again confirms
 (cf. Figure \ref{BBF}  for an $\e$-deformation of those)
   that these profiles belong to two branches originated at
 a supercritical saddle-node $R$-bifurcation at $R_{\rm min} \sim
 7.9$. Observe in Figure \ref{FCom2} other dotted profiles of a similar geometric structure that
 exist close to $R_{\rm min}^+$ and
 indeed correspond to the third basic pattern $F_2=F_{+2,-2,+2}$.
  %%% possibly occurred at saddle-node
 %%bifurcations.
 The intriguing global $R$-bifurcation diagram with saddle-node
 bifurcations
 turns out to be  similar to that observed in Section \ref{SectHom} by the
 $\e$-homotopy approach; see more explanations below.
 %%%remains an open problem to be
 %%discussed elsewhere.

 %%\com{To SIP: should we do that ? In a book, possibly...}

  %%%%%%%%%%%%%%%%%%%%%%%%%%%%%%%%%%%%%%%%%%%%%%%%%%
%%%%%%%%%%%%%%%%%%%%%%%%%%%%%%%%%%%%%%%%%%%%%%%%%%%%%%%%%%%%%

We do not need to develop more consistent theory of
$R$-compression in view of the following simple comment:

%%%\ssk

\noi\underline{\sc $\e$-homotopy and  $R$-compression can be
equivalent}. For simplicity, consider the analytic ODE problem
(\ref{nn3}) on $(-R,R)$, where we  perform the scaling
 $$
 F(y)= \sqrt \e \,
 V(z), \quad y= \e^{\frac 14} z \quad (\e>0).
 $$
 Then $\e$ is scaled out from the ODE and enters the interval, where the problem
 is posed:
  \beq
  \label{V11}
   \mbox{$
  V^{(4)}= - V + V^3 \quad \mbox{on} \quad \bigl(-\e^{\frac 14}R,\e^{\frac
  14}R \bigr).
   $}
  \eeq
  Therefore, in this particular case, the $\e$-deformation of the equation (\ref{nn3}) is
  equivalent to the $R$-compression for (\ref{V11}) as $\e>0$
  decreases. So these two approaches to Sturm's index of highly oscillatory
structures are essentially equivalent and hence lead to the same
results. In particular, this explains the striking phenomenon
(observed in Section \ref{SectHom}) that, at saddle-node
bifurcations, the profiles lose all their oscillations about zero.
%%%%%and hence ????????????????????????  What ??????

%%%%%%%%%%%%%%%%%%%%%%%%%%%%%%%%%%%%%%%%%%%%%%%%%%
\section{On essentially quasilinear extensions: gradient diffusivity}
 \label{SectExt1}

\subsection{Variational problems with $p$-Laplacian operators}

For simplicity, we consider PDEs (\ref{Cl1}) in 1D, where these
have simpler forms,
  \beq
  \label{Cl12}
  \begin{matrix}
  {\bf (I)} \quad u_t = -(|u_{xx}|^n u_{xx})_{xx} + |u|^n u \quad
  (\mbox{parabolic}),\,\,  \ssk\ssk \\
  {\bf (II)} \quad u_{tt} =  -(|u_{xx}|^n u_{xx})_{xx} + |u|^n u \quad
  (\mbox{hyperbolic}),\ssk\ssk \\
{\bf (III)} \quad u_t =  -(|u_{xx}|^n u_{xx})_{xxx} + (|u|^n
u)_{x} \quad
  (\mbox{NDE}). \quad\,\,
   \end{matrix}
   \eeq
 For the reaction-diffusion PDE {\bf (I)}, the blow-up solutions
 are the same, (\ref{RD.31}), with
 %%%%the ODE
 \beq
 \label{o1}
  \mbox{$
  \frac 1n \, f= -(|f''|^n f'')'' + |f|^n f.
  $}
  \eeq
  Using the scaling yields the basic quasilinear  ODE:
  \beq
  \label{o2}
 f= n^{-\frac 1n} F
 \quad \Longrightarrow \quad
 \fbox{$
-(|F''|^n F'')'' + |F|^n F - F=0 \quad \mbox{in} \quad \re.
 $}
  \eeq
  For the hyperbolic PDE {\bf (II)}, we construct the blow-up
  patterns (\ref{RD.31H}), where the same scaling yields (\ref{o2}).
  Finally, for the NDE {\bf (III)}, the TW compacton
  %%(\ref{RD.5})
  %%with
   $$
  u_{\rm c}(x,t)=f(x- \l t), \quad \mbox{with} \quad
   \l=-1,
   $$
   directly leads to the ODE in (\ref{o2}).

   In all the three cases, for the $N$-dimensional PDEs (\ref{Cl1}),
   we arrive at the elliptic PDE,
   \beq
   \label{el111}
   \fbox{$
   -\D(|\D F|^n \D F) + |F|^n F - F=0 \quad \mbox{in}
   \quad \re,
    $}
    \eeq
admitting  variational formulation in  $L^2(\ren)$ with the
potential (cf. (\ref{V1}))
 \beq
 \label{V1NN}
  \mbox{$
 E(F)= - \frac 1{n+2} \, \displaystyle \int_{\ren} |\D F|^{n+2} + \frac 1{n+2} \, \displaystyle \int_{\ren}
 |F|^{n+2} - \frac 12 \, \displaystyle \int_{\ren} F^2.
 $}
  \eeq
  Then the same Lusternik--Schnirel'man/fibering approaches can be applied to get a countable
  basic
   family of compactly supported solutions and next discover other
   countable sets of blow-up patterns, etc.
We claim that the most of principal results obtained above for the
{\em semilinear} elliptic problems can be extended to the {\em
quasilinear} one (\ref{el111}), though some local (e.g.,
oscillatory structures of solutions near finite interfaces) or
global aspects of the behaviour of patterns are more complicated.
 Some local oscillatory properties of solutions of changing sign
 %% near
 %%interfaces
    are discussed in \cite[pp.~246--249]{GSVR} and \cite{GpLap}.
 All these
problems admit more general $2m$th-order extensions along the same
lines as usual.

%%%\ssk

%%%%%%%%%%%%%%%%%%%%%%%%%%%%%%%%%%%%%%%%%%%%%%%%%%%%%
\subsection{Related non-variational problems: branching
``from" potential results}

Finally, let us mention that there exists
 a variety of slightly
changed PDEs {\bf (I)--(II)}, which, on exact blow-up  solutions,
reduce to elliptic or ODE problems that are principally {\em
non-variational}. A systematic study of such problems with
non-coercive, non-monotone, and non-potential operators touches
another wide area of open mathematical problems.

 For instance, consider a reaction-diffusion equation with the
 nonlinearity from (\ref{pn12}),
  \beq
  \label{par1}
  u_t= - \D (|\D u|^n \D u) +
 |u|^{p} u, \quad \mbox{where} \quad p>n>0.
  \eeq
The similarity blow-up solutions are not separable as in
(\ref{RD.31}) and have the form
 \beq
 \label{par2}
  \mbox{$
 u_{\rm S}(x,t) = (T-t)^{-\frac 1p} f(y), \quad y=x/(T-t)^\b, \quad
 \mbox{where} \quad  \b=
 \frac {p-n}{2p(n+2)}>0.
  $}
   \eeq
Substituting (\ref{par2}) into (\ref{par1}) yields the following
elliptic equation:
 \beq
 \label{par3}
  \fbox{$
- \D (|\D f|^n \D f)  - \b y \cdot \n f - \frac 1{p} \, f +
 |f|^{p} f=0 \quad \mbox{in \,\, $\ren$}.
 $}
 \eeq
For $p=n$ we have $\b=0$, so that (\ref{par3}) reduces to the
above variational equation (\ref{el111}).

We claim that, for  $p \not = n$, the differential operator in
 (\ref{par3}) is not variational (even for $n=0$; see
 \cite[\S~7]{GW2}),
 so all above techniques fail. Nevertheless, we also claim
 that the present variational analysis can and  does play a role for
 such problems. Namely, the variational problem (\ref{el111}) for $p=n$
 has the following meaning:
  \beq
  \label{vv1}
  \fbox{$
 \begin{matrix}
  \mbox{potential eq. (\ref{el111}) %%% for $p=n$
    initiates {\em branching}  at $p=n$  of solutions of non-potential} \ssk\\
    \mbox{eq.
 (\ref{par3}) from the Lusternik--Schnirel'man/fibering patterns
    $\{F_l\}$.}
    %%%at $p=n$.}
     %%% for
 %%%$p>n$}.
 \end{matrix}
 $}
  \eeq
  Such ideas well-correspond
 to
    classic branching theory; see Vainberg--Trenogin \cite{VainbergTr}. Therefore, we expect
that still, even not being variational,
 \beq
 \label{vv2}
 \mbox{(\ref{par3}) admits arbitrarily large
 %%(e.g., countable)
 number of  solutions for $p \approx n^+$}.
   \eeq
 Here we should assume that $p$ is sufficiently close to $n^+$,
 since, as usual, global continuation of local bifurcation branches is a
 difficult open problem.  Note that, for
$p>n$, solutions of (\ref{par3}) are not compactly supported in
general (for $p \in (0,n]$ these are). Applications of such a
 branching approach to (\ref{par3}) with $N=1$ and  $p>n$ and
 $p<n$ are given
in \cite{GpLap}.

  Branching phenomena for such nonlinear and
  degenerate operators as in (\ref{par3}) are not standard
  and demand a lot of work; see e.g., \cite{AlvGalTFE}, where thin film operators have been dealt with.
     On the other hand, the
  $p$-branching analysis in the semilinear case $n=0$ for equations such as (\ref{par3}) uses
  spectral theory of non self-adjoint operators and is easier;
 see examples in \cite{BGW1, GHUni, GW2}.

\smallskip

In other words, precisely (\ref{vv1}) is the actual role that
variational problems can play for describing finite or countable
sets
 of solutions of principally non-variational ones, so that
the results such as (\ref{vv2}) can be inherited from a suitable
potential asymptotic setting. Finding such a variational problem
by introducing a parameter as a good approximation of the given
non-potential one can be rather tricky, though for equations like
(\ref{par3}), this looks
 natural and
 straightforward.

%% As a final open problem, we

%%%non-variational ????

%%\ssk \ssk

%%%%%%%%%%%%%%%%%%%%%%%%%%%%%%%%%%%%%%%%%%%%%%%%%%%
%%{\bf Acknowledgements.}  The authors would like to thank  A.E.
%%Shishkov  for a number of useful
%%%  discussions.

%%%%%%%%%%%%%%%%%%%%%%%%%%%%%%%%%%%%%%%%%%%%%%%%%%%%%%%%%%%%%%%%%%%%%%%%%%%

%%%%%%%%%%%%%%%%%%%%%%%%%%%%%%%%%%%%%%%%%%%%%%%

\begin{thebibliography}{10}

\bibitem
{AlvGalTFE}
 P.~\'Alvarez-Caudevilla and V.A.~Galaktionov,
 {\em Local bifurcation-branching analysis of global and
``blow-up" patterns for a fourth-order thin film equation},
submitted (arXiv:1009.5864).



 \bibitem
{Rynn2}
 R.~Bari and B.~Rynne, {\em Solution curves and exact multiplicity results
 for $2m$th order boundary value problems},  J.~Math. Anal. Appl., {\bf 292} (2004), 17--22.




\bibitem %[Berger]
 {Berger}
  M.~Berger, {\rm Nonlinearity and Functional Analysis}, Acad.
  Press, New York, 1977.


 \bibitem
  {BrBra}
  H.~Brezis and F.~Brawder, {\em Partial Differential Equations in the 20th Century},
  Adv. in Math., {\bf 135} (1998), 76--144.




\bibitem
 {BGW1}
C.J.~Budd, V.A.~Galaktionov, and J.F.~Williams, {\em Self-similar
blow-up
  in higher-order semilinear parabolic equations}, SIAM J.~Appl. Math.,
 {\bf 64} (2004), 1775--1809.





 \bibitem
{Coff72}
   C.V.~Coffman, {\em Uniqueness of the ground state for $\D u-u+u^3=0$ and a
   variational characterization of other solutions},
 Arch. Rat. Mech. Anal., {\bf 46} (1972), 81--95.




\bibitem
{Elias} U.~Elias, {\em Eigenvalue problems for the equation $L y +
p(x) y=0$}, J. Differ. Equat., {\bf 29} (1978), 28--57.





\bibitem
 {Galp1}
 V.A.~Galaktionov, {\em On interfaces and oscillatory solutions of
 higher-order semilinear parabolic equations with non-Lipschitz
 nonlinearities}, {Stud. Appl. Math.,} {\bf 117} (2006), 353--389.

\bibitem
 {GpLap}
 V.A.~Galaktionov,
{\em Three types of  self-similar blow-up for  the fourth-order
$p$-Laplacian equation with source}, J.~Comp. Appl. Math., {\bf
223} (2009), 326--355 (arXiv:0903.0981).


\bibitem
 {GHUni}
 V.A.~Galaktinov and P.J.~Harwin,
  {\em Non-uniqueness and  global similarity solutions   for
a higher-order semilinear parabolic equation}, Nonlinearity, {\bf
18} (2005), 717--746.



\bibitem
 {GMPSob}
 V.A.~Galaktionov, E.~Mitidieri,  and S.I.~Pohozaev,
 {\em Variational approach to complicated similarity solutions
of higher-order
 nonlinear evolution equations of parabolic, hyperbolic, and nonlinear dispersion
 types},
%%% parabolic,\\
 %%%hyperbolic, and nonlinear dispersion
%%  PDEs: an analytic-numerical
%%approach}
In: Sobolev Spaces in Mathematics. II, Appl. Anal. and Part.
Differ. Equat., Series: Int. Math. Ser., Vol. {\bf 9}, V.~Maz'ya
Ed., Springer, New York, 2009 (an extended version in arXiv:0902.1425).


\bibitem
 {GMPSobIarX}
 V.A.~Galaktionov, E.~Mitidieri,  and S.I.~Pohozaev,
 {\em Variational approach to complicated similarity solutions
of higher-order
 nonlinear
 PDEs. I},
arXiv:0902.1425.


\bibitem
 {GMPHomIINA}
 V.A.~Galaktionov, E.~Mitidieri,  and S.I.~Pohozaev,
 {\em Variational approach to complicated similarity solutions
of higher-order
 nonlinear
 PDEs. II}, Nonl. Anal: RWA, DOI: 10.1016/j.nonrwa.2011.03.001.
%%%arXiv:0902.1425.






 \bibitem
{GPnde}
 V.A.~Galaktionov and S.I.~Pohozaev, {\em
 Third-order nonlinear dispersive equations: shocks, rarefaction,
 and blow-up waves},
%% On shock,
%%rarefaction waves, and compactons for odd-order nonlinear
%%dispersion PDEs},
 %% Differ. Equat.,
 {Comput. Math. Math. Phys.,} {\bf 48} (2008), 1784--1810
 (arXiv:0902.0253).
%%% Advances Differ. Equat., {\bf 10} (2005), 635--674.
 %%submitted.





\bibitem
 {GSVR} V.A.~Galaktionov and S.R.~Svirshchevskii, Exact Solutions and
 Invariant Subspaces of Nonlinear Partial Differential Equations in Mechanics and Physics,
  Chapman$\,\&\,$Hall/CRC, Boca Raton,
Florida,
 2007.


\bibitem
  {AMGV}
 V.A.~Galaktionov and J.L.~V\'azquez,
 {\rm A Stability Technique  for Evolution Partial Differential Equations.
 A Dynamical Systems Approach}, {\rm Progr. in Nonl. Differ.
 Equat. and their Appl.,} {\bf 56},
Birkh\"auser Boston, Inc., MA, 2004.


\bibitem{GW2}
V.A.~Galaktionov and J.F.~Williams,
% \bysame,
{\em On very singular similarity solutions of a higher-order
  semilinear parabolic equation}, Nonlinearity, {\bf 17} (2004), 1075--1099.
  %({\tt
  %http://www.maths.bath.ac.uk/MATHEMATICS/preprints.html}).



\bibitem
 {KKVV00}
  W.D.~Kalies, J.~Kwapisz, J.B.~VandenBerg, and
  R.C.A.M.~VanderVorst, {\em Homotopy classes for stable periodic
  and chaotic patterns in fourth-order Hamiltonian systems,}
  Commun. Math. Phys., {\bf 214} (2000), 573--592.



\bibitem
 {Kor04}
  P.~Korman, {\em Uniqueness and exact multiplicity of solutions
  for a class of fourth-order semilinear problems}, Proc. Roy.
  Soc. Edinburgh, Sect. A, {\bf 134} (2004), 179--190.


\bibitem %[Kras]
{Kras} M.A.~Krasnosel'skii, {\rm Topological Methods in the Theory
of Nonlinear Integral Equations}, Pergamon Press, Oxford/Paris,
1964.



\bibitem
{KrasZ} M.A.~Krasnosel'skii and P.P.~Zabreiko, {\rm Geometrical
Methods of Nonlinear Analysis}, Springer-Verlag, Berlin/Tokyo,
1984.



 \bibitem %%[MitPoh]
 {Mer05}
F.~Merle  and P.~Raphael, {\em On a sharp lower bound on the
blow-up rate for the $L^2$ critical nonlinear Schr\"odinger
equation}, J.~Amer. Math. Soc., {\bf 19} (2005),  37--90.

\bibitem %%[MitPoh]
 {MerleR05}
F.~Merle  and P.~Raphael, {\em The blow-up dynamics and upper
bound on the blow-up rate for  critical nonlinear Schr\"odinger
equation}, Ann. Math., {\bf 161} (2005), 157--222.

\bibitem %%[MitPoh]
 {MerR052}
F.~Merle  and P.~Raphael, {\em Profiles and quantization of the
blow-up mass for  critical nonlinear Schr\"odinger equation},
Comm. Math. Phys., {\bf 253} (2005), 675--704.



\bibitem %%[MitPoh]
 {MitPoh}
E.~Mitidieri and S.I.~Pohozaev, {\rm  A priori estimates and the absence of solutions of nonlinear partial differential equations and inequalities.},
 Proc. Steklov Inst. Math., Vol. {\bf 234}, Intern.
Acad. Publ. Comp. Nauka/Interperiodica, Moscow, 2001.



\bibitem
{PelTroy} L.A. Peletier and W.C. Troy, {\rm Spatial Patterns:
Higher Order Models in Physics and Mechanics}, Birkh\"auser,
Boston/Berlin, 2001.


\bibitem
{Poh0} S.I.~Pohozaev,  {\em On an approach to nonlinear
equations}, Soviet Math. Dokl.,
 {\bf 20} (1979), 912--916.

 \bibitem
{PohFM} S.I. Pohozaev, \rm {\em The fibering method in nonlinear
variational problems}, Pitman Research Notes in Math., Vol. {\bf
365}, Pitman, 1997, pp. 35--88.







\bibitem
{Rynn1}
 B.~Rynn, {\em Global bifurcation for $2m$th-order boundary
value problems and infinitely many solutions of superlinear
problems},  J. Differ. Equat., {\bf 188} (2003), 461--472.




\bibitem %[SGK]
{SGKM}  A.A.~Samarskii, V.A.~Galaktionov, S.P.~Kurdyumov, and
A.P.~Mikhailov, {Blow-up in Quasilinear Parabolic Equations}, \rm
Walter de Gruyter, Berlin/New York, 1995.

\bibitem %[SZKM2]
{SZKM2} A.A.~Samarskii, N.V.~Zmitrenko, S.P.~Kurdyumov, and
A.P.~Mikhailov, {\em Thermal structures and fundamental length in
a medium with non-linear heat conduction and volumetric heat
sources}, {Soviet Phys. Dokl.,} {\bf 21} (1976), 141--143.



%%\bibitem
%%{SW} J. Shi and J. Wang, {\em Morse indices and exact multiplicity
%%of solutions to semilinear elliptic problems}, Proc. Amer. Math.
%%Soc.,  {\bf 127} (1999), 3685--3695.

\bibitem
{Shi2}
 A.E.~Shishkov, {\em Dead cores and instantaneous
compactification of the supports of energy solutions of
quasilinear parabolic
equations of arbitrary order}, {Sbornik: Math.}, {\bf 190}    %:12
 (1999), 1843--1869.


\bibitem %[VainbergTr]
{VainbergTr} M.A. Vainberg and V.A. Trenogin, {\rm Theory of
Branching of Solutions of Non-Linear Equations}, Noordhoff Int.
Publ., Leiden, 1974.



\bibitem
 {VV02}
   J.B.~Van Den Berg and
  R.C.~Vandervorst, {\em Stable
   patterns for fourth-order parabolic equations,}
  Duke Math.~J., {\bf 115} (2002), 513--558.




\end{thebibliography}
\end{document}